\documentclass[11pt]{article}

\usepackage[utf8]{inputenc}

\usepackage{graphicx}
\usepackage{xcolor}
\usepackage{dsfont}
\usepackage{amsmath,amsthm,amssymb}
\usepackage[margin=1in]{geometry}
\usepackage[pdftex,colorlinks=true,linkcolor=blue,citecolor=blue,urlcolor=blue]{hyperref}
\usepackage{comment}
\usepackage{bbm}
\usepackage{tikz}
\usetikzlibrary{decorations.markings}
\usetikzlibrary{decorations.markings,decorations.pathmorphing}

\renewcommand{\d}{{\mathrm d}}

\makeatletter
\newcommand\xleftrightarrow[2][]{%
  \ext@arrow 9999{\longleftrightarrowfill@}{#1}{#2}}
\newcommand\longleftrightarrowfill@{%
  \arrowfill@\leftarrow\relbar\rightarrow}
\makeatother

\usepackage[normalem]{ulem}
\usepackage{cite}
\usepackage{tikz}
\usepackage{tcolorbox}

\def\ba{\begin{array}}
\def\ea{\end{array}}

\def\0{{\bf 0}}

\newcommand{\e}{\mathrm{e}}

\newcommand{\subjclass}[1]{%
  \par\noindent
  \textbf{2020 Mathematics Subject Classification.}
  #1\par
}

\DeclareMathOperator{\Var}{Var}

\renewcommand{\Re}{\operatorname{Re}}

\newcommand{\be}{\begin{equation}}
\newcommand{\ee}{\end{equation}}
\newcommand{\bea}{\begin{eqnarray}}
\newcommand{\eea}{\end{eqnarray}}
\newcommand{\bes}{\begin{equation*}}
\newcommand{\ees}{\end{equation*}}
\newcommand{\beas}{\begin{eqnarray*}}
\newcommand{\eeas}{\end{eqnarray*}}

\makeatletter
\newtheorem*{rep@theorem}{\rep@title}
\newcommand{\newreptheorem}[2]{%
\newenvironment{rep#1}[1]{%
 \def\rep@title{#2 \ref{##1} (restated)}%
 \begin{rep@theorem}}%
 {\end{rep@theorem}}}
\makeatother

\allowdisplaybreaks

\newtheorem{thm}{Theorem}[section]
\newtheorem*{thm*}{Theorem}
\newtheorem{cor}[thm]{Corollary}
\newtheorem{con}[thm]{Conjecture}
\newtheorem{lem}[thm]{Lemma}
\newtheorem*{lem*}{Lemma}

\newtheorem{prop}[thm]{Proposition}

\newreptheorem{thm}{Theorem}
\newreptheorem{lem}{Lemma}
\newreptheorem{cor}{Corollary}

\newtheorem{theo}{Theorem}[section]

\newtheorem{remark}[theo]{Remark}
\newtheorem{definition}[theo]{Definition}

\newenvironment{keywords}{}{}

\def\i{{\textnormal{i}}}

\usepackage[title]{appendix}
\usepackage{pgfplots}
\usepackage{thmtools}
\usepackage{thm-restate}

\providecommand{\keywords}[1]
{ \textbf{Keywords:} #1 }

\usepackage{authblk}

\date{}

\begin{document}
\title{ \large  \bf Complex Moments of the Derivative of the Riemann Zeta Function}
\author{
Joseph Najnudel$^{1}$\thanks{$^1$School of Mathematics, University of Bristol, Fry Building,
Woodland Road, Bristol, BS8 1UG, United Kingdom.
\texttt{joseph.najnudel@bristol.ac.uk}}
~and Fei Wei$^{2}$\thanks{$^2$Key Laboratory of Mathematics and Information Networks (Beijing University of Posts and Telecommunications), Ministry of Education, Science Building, No. 1 Nanfeng Road, Beijing 102206, China.
\texttt{weif0831@gmail.com}}
}
\date{}
\maketitle
\vspace{-1.5em}
\begin{abstract}
We investigate the complex moments of
the derivative of the Riemann zeta function
through a random Euler product model.
For each fixed $\sigma>1/2$ and every
$s\in\mathbb{C}$ with $\operatorname{Re}(s)>0$,
we show, under the Lindel\"of hypothesis,
that the $2s$-th moment of
$|\zeta'(\sigma+\mathrm{i}t)|$ over $t\in[1,T]$
converges, as $T\to\infty$,
to the corresponding moment of the derivative
of the random Euler product.
We then establish, unconditionally, an
asymptotic formula for the latter moment
as $\sigma\downarrow 1/2$ throughout the
half-plane $\operatorname{Re}(s)>0$.
For positive real moments, our approach
is based on a change-of-measure argument
and a Gaussian limit theorem for the
logarithmic derivative under a tilted
probability measure.
For general complex moments, we combine
Fourier analysis of weighted complex
measures with regularized Riesz
distributions. As an application, under the Lindel\"of
hypothesis, we resolve the conjecture
of Simm and Wei
[J. Lond. Math. Soc. 113 (2026)]
on the complex moments of the derivative
of the Riemann zeta function for every
$s\in\mathbb C$ with $\operatorname{Re}(s)>0$.
We establish the predicted asymptotic
formula as $\sigma\downarrow1/2$, with
the leading order $s^2+2s$ and an explicit
leading coefficient involving an arithmetic
factor and the confluent hypergeometric
function of the first kind.
The conjecture also holds unconditionally
when $\operatorname{Re}(s)\in (0,2)$.

~\\
\keywords{{\bf Keywords:} Riemann zeta function, complex moments,
random Euler products, Gaussian limit theorem,
Riesz distributions}~\\

\subjclass{11K65, 11M06, 60E10}
\end{abstract}

\setcounter{tocdepth}{2}
\tableofcontents

\section{Introduction}
\label{sec:introduction}

The study of moments of the Riemann zeta
function has a long history in analytic
number theory, dating back to the classical
work of Hardy and Littlewood
\cite{hardy1923lindelofs}.
In particular, the Lindel\"of hypothesis
is equivalent to the assertion that,
for every positive integer $s$ and
every fixed $\sigma>1/2$,
\begin{equation}
\label{eq:intro-classical-moments}
\lim_{T\to\infty}
\frac{1}{T}\int_1^T
|\zeta(\sigma+\mathrm{i}t)|^{2s}
\,\mathrm{d}t
=
\sum_{n=1}^{\infty}
\frac{d_s(n)^2}{n^{2\sigma}},
\end{equation}
where $d_s(n)$ denotes the $s$-fold
divisor function; see
\cite[Theorem~13.2]{titchmarsh1986theory}. The case $s=1$ of
\eqref{eq:intro-classical-moments}
is known unconditionally by the
classical mean-value theorem for
Dirichlet series; see also
\cite{HardyLittlewood1916}.
The case $s=2$ was established
by Hardy and Littlewood
\cite{hardy1923approximate}.
However, for $s=3$, the identity
remains unproved throughout the
entire range $\sigma>1/2$. 
From \eqref{eq:intro-classical-moments},
we see that, under the Lindel\"of
hypothesis, for every positive integer
$s$,
\begin{equation}
\label{eq:intro-classical-moment-asymptotics}
\lim_{T\to\infty}
\frac{1}{T}\int_1^T
|\zeta(\sigma+\mathrm{i}t)|^{2s}
\,\mathrm{d}t
\sim
\frac{a_s}{(2\sigma-1)^{s^2}},
\qquad \sigma\downarrow1/2,
\end{equation}
where $a_s$ is the arithmetic factor
defined by
\begin{equation}
\label{definitionofthearithmeticfactor}
a_s
=
\prod_p
(1-p^{-1})^{s^2}
\sum_{k=0}^{\infty}
\left(\frac{(s)_k}{k!}\right)^2
p^{-k},
\end{equation}
with $(s)_k=s(s+1)\cdots(s+k-1)$ for $k\geq1$
and $(s)_0=1$ denoting the rising factorial.
In particular, the asymptotic formula
\eqref{eq:intro-classical-moment-asymptotics}
holds unconditionally for $s=1$ and $s=2$.

Recently, Simm and Wei \cite{simm2026moments}
studied the analogous problem for the
moments of the derivative of the Riemann
zeta function, $\zeta'$.
More precisely, they proposed the following
conjecture, analogous to
\eqref{eq:intro-classical-moment-asymptotics},
concerning the asymptotic behaviour of
the moments of $\zeta'(\sigma+\mathrm{i}t)$
as $\sigma\downarrow1/2$.
\begin{con}
\label{conjectureonthemomentsoffthecriticalline}
For every $s\in\mathbb C$ with
$\operatorname{Re}(s)>0$, one has
\begin{equation}
\label{thelimtintherandommatrixside}
\lim_{T\to\infty}
\frac{1}{T}\int_1^T
|\zeta'(\sigma+\mathrm{i}t)|^{2s}
\,\mathrm{d}t
\sim a_{s}
\frac{\e^{-s^2}\Gamma(s+1)
{}_1F_1(s+1;1;s^2)}
{(2\sigma-1)^{s^2+2s}},
\qquad \sigma\downarrow1/2,
\end{equation}
where $a_s$ is the arithmetic factor
defined by (\ref{definitionofthearithmeticfactor}).
\end{con}

Compared with
\eqref{eq:intro-classical-moment-asymptotics},
this conjecture predicts that
differentiation changes the exponent
of $(2\sigma-1)^{-1}$ from $s^2$
to $s^2+2s$ and introduces an additional
factor involving the confluent
hypergeometric function arising from
random matrix theory.
The conjecture is known to hold for
all positive integers $s$.
The case $s=1$ follows from the
mean-value theorem for Dirichlet series,
while the case $s=2$ was proved
unconditionally in
\cite[Theorem~1.7]{simm2026moments}.
For every integer $s\geq3$, it was
established under the Lindel\"of hypothesis
in \cite[Theorem~1.8]{simm2026moments}. A difficulty in studying these derivative
moments is that differentiation introduces
logarithmic factors into the coefficients
of the associated Dirichlet series,
which are no longer multiplicative.
Consequently, the Euler product approach
used for the moments of the Riemann
zeta function cannot be applied directly.

The conjecture above is motivated by
the work of Simm and Wei \cite{simm2026moments},
who studied complex moments of the
derivatives of CUE characteristic polynomials
at points inside the unit disc.
By analyzing the asymptotic behaviour of
these moments as the spectral parameter
approaches the unit circle, they obtained
an asymptotic formula with the same
leading order $s^2+2s$ as in
Conjecture~\ref{conjectureonthemomentsoffthecriticalline}.
This suggests a corresponding asymptotic
formula for the moments of the derivative
of the Riemann zeta function, with an
additional arithmetic factor $a_s$. The idea of using characteristic polynomials
of random unitary matrices to model moments
of the Riemann zeta function goes back to
Keating and Snaith \cite{KeatingSnaith2000},
who proposed conjectures for the moments
of the Riemann zeta function on the
critical line based on calculations in
the CUE.
This approach was subsequently extended
to moments of the derivative of the
Riemann zeta function; see, for e.g.,
\cite{hughes2001theis, HughesKeatingOConnell2000,
ConreyRubinsteinSnaith2006, forrester2006boundary,Winn2012,  Basor2019, AssiotisKeatingWarren2022,BothnerWei2026,assiotis2026painlev}.
For related results on moments of
higher-order derivatives, see, e.g.,
\cite{KeatingWei2024I,KeatingWei2024II,AssiotisGunesKeatingWei2026,
GroverMezzadriSimm2026}.

In the present paper, we establish the
asymptotic formula predicted by
Conjecture~\ref{conjectureonthemomentsoffthecriticalline}
under the Lindel\"of hypothesis for every
$s\in\mathbb C$ with $\operatorname{Re}(s)>0$.
Moreover, the result holds unconditionally
when $0<\operatorname{Re}(s)<2$.
Our approach is based on probabilistic
models of the Riemann zeta function. The classical work of Bohr and Jessen
\cite{BohrJessen1930,BohrJessen1932}
established the existence of limiting
value distributions for
$\log\zeta(\sigma+\mathrm{i}t)$
for fixed $\sigma>1/2$. This result was further developed by
Jessen and Wintner \cite{JessenWintner1935},
who investigated the associated limiting
distributions using infinite convolutions
and Fourier-analytic methods. Bagchi \cite{Bagchi1981} subsequently
established a functional limit theorem
for the vertical translates
$\zeta(\,\cdot+\mathrm{i}t)$,
with the limiting distribution
described by a random Euler product. 
The random Euler product arising in this
limiting distribution is constructed by
replacing the prime-dependent phases
$p^{-\mathrm{i}t}$ with independent
random variables $X(p)$ uniformly
distributed on the unit circle.
Extending these variables completely
multiplicatively gives a Steinhaus
random multiplicative function $X(n)$
and the associated random Euler product
\[
\zeta(w,X)
=
\prod_p(1-X(p)p^{-w})^{-1},
\qquad \operatorname{Re}(w)>1/2.
\]
The study of random multiplicative functions
has a long history, dating back to the
early work of Wintner \cite{Wintner1944}.
Various aspects of the partial sums of
random multiplicative functions, including
their almost sure behavior and moment
estimates, have been extensively studied;
see, for example,
\cite{LauTenenbaumWu2013,
Harper2019HighMoments,Harper2020LowMoments,
AymoneHeapZhao2021,SoundararajanXu2023, Atherfold2025}
and the references therein.
Random Euler products have also played
an important role in the study of value
distributions and moments of various
$L$-functions; see, for example,
\cite{CogdellMichel2004,LiuRoyerWu2008,
Lamzouri2010,Lamzouri2011}
and the references therein.
In the context of the Riemann zeta function,
Arguin, Belius, and Harper
\cite{ArguinBeliusHarper2017}
studied the extreme values of a randomized
model based on truncated Euler products;
see also the subsequent work of Arguin,
Bourgade, and Radziwi{\l}{\l}
\cite{ArguinBourgadeRadziwill2020,
ArguinBourgadeRadziwill2023}
on the Riemann zeta function itself.
These developments illustrate the
importance of random Euler products
in the probabilistic study of the
Riemann zeta function and other
$L$-functions.
\subsection{Main results}
In the present paper, we investigate
the complex moments of the derivative
$\zeta'(\sigma,X)$ of the random Euler
product and their asymptotic behaviour
as $\sigma\downarrow1/2$. For comparison, the corresponding
moments without differentiation admit
a direct evaluation using the independence
of the prime-indexed random variables.
Indeed, the Euler product representation
implies that, for every $s\in\mathbb C$,
\begin{equation}
\label{eq:random-zeta-moment-asymptotics}
\lim_{\sigma\downarrow1/2}
(2\sigma-1)^{s^2}
\mathbb E\left[
|\zeta(\sigma,X)|^{2s}
\right]
=a_s,
\end{equation}
where $a_s$ is the arithmetic factor
defined in
\eqref{definitionofthearithmeticfactor}.
Here and throughout the paper, $\mathbb E$
denotes expectation with respect to the
probability measure $\mathbb P$ defined
in \eqref{eq:probability-measure-definition}.

In contrast, differentiation introduces
a factor of $\log n$ into the Dirichlet
coefficients, which are no longer
multiplicative. Consequently, the moments
of $\zeta'(\sigma,X)$ require a different
analysis. Our first main result, proved in
Subsection~\ref{thm41}, determines
the limiting behaviour of the
normalized complex moments of
$\zeta'(\sigma,X)$ as
$\sigma\downarrow1/2$.

\begin{thm}
\label{thm:complex-random-derivative-moments}
For every $s\in\mathbb C$ with
$\operatorname{Re}(s)>0$, one has
\begin{equation}
\label{eq:complex-random-derivative-main}
\lim_{\sigma\downarrow1/2}
(2\sigma-1)^{s^2+2s}
\mathbb E\left[
|\zeta'(\sigma,X)|^{2s}
\right]
=
a_s \mathrm e^{-s^2}\Gamma(s+1)
{}_1F_1(s+1;1;s^2),
\end{equation}
where $a_s$ is the arithmetic factor
defined in
\eqref{definitionofthearithmeticfactor},
and ${}_1F_1(a;b;z)$ denotes the confluent
hypergeometric function of the first kind.
Moreover, the convergence is locally uniform
in $s$ throughout the half-plane
$\operatorname{Re}(s)>0$.
\end{thm}
We next establish a connection between the
complex moments of the derivative of the
random Euler product and the corresponding
mean values of the derivative of the
Riemann zeta function.
More precisely, we show that, under the
Lindel\"of hypothesis, these moments
coincide for every fixed $\sigma>1/2$.
We also establish this identity
unconditionally in a restricted range
of the complex parameter $s$.
The following result is proved in
Section~\ref{section2}.

\begin{thm}
\label{mainresult}
Let $\sigma>1/2$ be fixed, and let
$\zeta(w,X)$ be defined by
\eqref{defofrandomzeta} for
$\operatorname{Re}(w)>1/2$.
Assume the Lindel\"of hypothesis.
Then, for every $s\in\mathbb C$ with
$\operatorname{Re}(s)>0$, one has
\begin{equation}
\label{eq:complex-moment-main4}
\lim_{T\to\infty}
\frac{1}{T}
\int_1^T
|\zeta'(\sigma+\mathrm{i}t)|^{2s}
\,\mathrm{d}t
=
\mathbb E\left[
|\zeta'(\sigma,X)|^{2s}
\right].
\end{equation}
Moreover, \eqref{eq:complex-moment-main4}
holds unconditionally for every
$s\in\mathbb C$ with
$0<\operatorname{Re}(s)<2$.
\end{thm}

Combining
Theorems~\ref{thm:complex-random-derivative-moments}
and~\ref{mainresult}, we obtain the following
consequence concerning Conjecture \ref{conjectureonthemomentsoffthecriticalline}.

\begin{cor}
\label{onemainresult}
Assume the Lindel\"of hypothesis.
Then, for every $s>0$, Conjecture \ref{conjectureonthemomentsoffthecriticalline} holds. Moreover, Conjecture \ref{conjectureonthemomentsoffthecriticalline} holds unconditionally for every
$s\in \mathbb{C}$ with $\textnormal{Re}(s)\in (0,2)$ and $s=2$.
\end{cor}

We now outline the main steps in the proof
of Theorem~\ref{thm:complex-random-derivative-moments}.
We first consider real parameters $s>0$,
for which a change-of-measure argument
allows us to analyze the moments of the
logarithmic derivative of the random
Euler product under a tilted probability
measure. For complex parameters
$s$ with $\operatorname{Re}(s)>0$,
we employ Fourier--Riesz methods
to establish the corresponding asymptotic
formula.

\subsection{Real moments via a tilted probability} \label{realsinintroduction}

The value distribution of the logarithmic
derivative of the Riemann zeta function
has been investigated by several authors.
In particular, Lester \cite{Lester2014}
established a central limit theorem for
$\zeta'/\zeta(\sigma+\mathrm{i}t)$
in a range where $\sigma$ approaches
the critical line from the right.
In a related direction, Fazzari
\cite{Fazzari2021,Fazzari2021WeightedValues}
investigated weighted central limit
theorems for the logarithm of the
Riemann zeta function on the critical
line, under the Riemann hypothesis.

The logarithmic derivative
$L'(\sigma,X)$ of the random Euler
product provides a probabilistic
counterpart to
$(\zeta'/\zeta)(\sigma+\mathrm{i}t)$.
For every fixed $\sigma>1/2$,
the limiting value distribution of
the latter as $T\to\infty$ is
described by the distribution of
$L'(\sigma,X)$.
In the present work, for real $s>0$,
we study $L'(\sigma,X)$ under a
tilted probability measure obtained
by weighting the original probability
measure with $|\zeta(\sigma,X)|^{2s}$.
Under this new measure, the
prime-indexed random variables
$\{X(p)\}_p$ remain independent,
although their marginal distributions
are modified.
This allows us to apply the
multivariate Lindeberg--Feller
central limit theorem to establish
a Gaussian limit with a nonzero
mean arising from the change of
measure. More precisely, we obtain
the following result.

\begin{prop}
\label{prop:tilted-gaussian-limit}
Let $s>0$ be fixed. As $\sigma\downarrow1/2$, under
the tilted probability measure
$\widetilde{\mathbb P}_{\sigma,s}$,
\begin{equation}
\label{eq:tilted-gaussian-limit-introduction}
-(2\sigma-1)L'(\sigma,X)
\xrightarrow{\mathrm d}s+G,
\end{equation}
where $L'(\sigma,X)$ is defined by the prime series
in \eqref{eq:derivative-random-log-series}, which
converges in
$L^2(\Omega,\widetilde{\mathbb P}_{\sigma,s})$,
and $G$ is a standard circularly symmetric complex
Gaussian random variable satisfying
\[
\mathbb E[G]=0,\qquad
\mathbb E[|G|^2]=1,\qquad
\mathbb E[G^2]=0.
\]
Equivalently, $G=G_1+\mathrm{i}G_2$, where
$G_1$ and $G_2$ are independent $N(0,1/2)$
random variables.
\end{prop}
Proof of Proposition~\ref{prop:tilted-gaussian-limit} is given in Subsection \ref{33}. To pass from convergence in distribution to convergence of moments, we establish uniform integrability using Rosenthal's inequality for sums of independent, centered random variables indexed by the primes. More precisely, for some \(q>2s\) and \(\sigma_0>1/2\), we prove that
\bea\label{uniformity}
\sup_{1/2<\sigma\leq\sigma_0}
\widetilde{\mathbb E}_{\sigma,s}
\left[
|(2\sigma-1)L'(\sigma,X)|^q
\right]<\infty.
\eea
Consequently, the Gaussian limit also holds at the level of moments:
\[
\lim_{\sigma\downarrow1/2}
(2\sigma-1)^{2s}
\widetilde{\mathbb E}_{\sigma,s}
\left[|L'(\sigma,X)|^{2s}\right]
=
\mathbb E\left[|s+G|^{2s}\right].
\]
Combining this convergence with the asymptotic formula for \(\mathbb E\left[|\zeta(\sigma,X)|^{2s}\right]\), obtained from the Euler product, yields the desired asymptotics for the real moments of \(\zeta'(\sigma,X)\). The resulting Gaussian factor can be evaluated explicitly
in terms of a confluent hypergeometric function:
\begin{equation}
\label{eq:Gaussian-factor-introduction2}
\mathbb E\left[|s+G|^{2s}\right]
=
\e^{-s^2}\Gamma(s+1)
{}_1F_1(s+1;1;s^2),
\end{equation}
where ${}_1F_1$ denotes the confluent hypergeometric
function of the first kind. This agrees precisely with the formula in
\eqref{thelimtintherandommatrixside}.
Both \eqref{uniformity} and
\eqref{eq:Gaussian-factor-introduction2}
are established in
Proposition~\ref{prop:random-zeta-derivative-moment-asymptotics}.

\subsection{Complex moments via Fourier--Riesz analysis}

The preceding argument for real $s>0$ in Subsection \ref{realsinintroduction} relies
essentially on the positivity of the weight
$|\zeta(\sigma,X)|^{2s}$, which allows us to
introduce a tilted probability measure and
apply probabilistic limit theorems.
For non-real $s$ with $\operatorname{Re}(s)>0$,
the same weighting procedure gives rise to
a complex measure rather than a probability measure.
Although the probabilistic argument does not
extend directly to this setting, the resulting
complex measure provides a natural framework
for studying complex moments through Fourier
analysis and the theory of Riesz distributions.

\subsubsection{Weighted complex measures and Fourier asymptotics}

Write
\bea\label{defofTsigma}
T_\sigma=(2\sigma-1)L'(\sigma,X).
\eea
For $\operatorname{Re}(s)>0$, we introduce
the weighted complex Borel measure
\begin{equation}
\label{defofthemeasure}
\nu_{\sigma,s}(B)
:=
(2\sigma-1)^{s^2}
\mathbb E\left[
|\zeta(\sigma,X)|^{2s}
\mathbf 1_{\{T_\sigma\in B\}}
\right],
\qquad
B\in\mathcal B(\mathbb R^2).
\end{equation}
By the fixed-$\sigma$ moment estimates
established in
Lemma~\ref{lem:basic-random-Euler-product},
this is a finite complex Borel measure.
For real $s>0$, the normalized measure
\[
\frac{\nu_{\sigma,s}}
{\nu_{\sigma,s}(\mathbb R^2)}
\]
coincides with the distribution of $T_\sigma$
under the tilted probability measure
$\widetilde{\mathbb P}_{\sigma,s}$.
For non-real $s$, however, $\nu_{\sigma,s}$
is generally a complex measure and does not
admit such a probabilistic interpretation. Instead, the complex measure $\nu_{\sigma,s}$
allows us to express the weighted complex moments
as integrals with respect to a finite complex
Borel measure. By studying its associated
Fourier transform $H_{\sigma,s}$, defined in (\ref{eq:Hsigma-def-final}) below, we establish
the asymptotic behaviour of the weighted
joint Fourier transform, together with
suitable decay and regularity estimates. These properties allow us to apply the Fourier
transform identity for the homogeneous
distribution $|x|^{2s}$ established in
\eqref{1023evening}, together with the
meromorphic continuation of Riesz distributions,
to derive a Fourier--Riesz representation
of the complex moments. More precisely, we consider the associated
weighted joint Fourier transform
\begin{equation}
\label{eq:Hsigma-def-final}
\begin{aligned}
H_{\sigma,s}(\lambda,u)
&:=
\int_{\mathbb R^2}
\e^{\i(\lambda x_1+u x_2)}
\,\d\nu_{\sigma,s}(x)\\
&=
(2\sigma-1)^{s^2}
\mathbb E\left[
|\zeta(\sigma,X)|^{2s}
\e^{\i\lambda\operatorname{Re}T_\sigma
+\i u\operatorname{Im}T_\sigma}
\right].
\end{aligned}
\end{equation}
Our first key result describes the limiting
behaviour of this transform as
$\sigma\downarrow1/2$.

\begin{prop}
\label{prop:compact-joint-transform-rigorous}
As $\sigma\downarrow1/2$,
\begin{equation}
\label{eq:compact-joint-transform-rigorous}
H_{\sigma,s}(\lambda,u)
\longrightarrow
a_s
\exp\left(
-\i s\lambda-\frac{\lambda^2+u^2}{4}
\right)
\end{equation}
locally uniformly for
\[
\operatorname{Re}(s)>0,
\qquad
(\lambda,u)\in\mathbb R^2.
\]
Here \(H_{\sigma,s}(\lambda,\mu)\) and \(a_s\) are defined in \eqref{eq:Hsigma-def-final} and \eqref{definitionofthearithmeticfactor}, respectively.
Moreover, on every bounded complex polydisc in $(\lambda,u)$, the
Euler-product expression for $H_{\sigma,s}$ admits a holomorphic continuation, and the
same convergence holds locally uniformly.
\end{prop}
Proof of Proposition \ref{prop:compact-joint-transform-rigorous} is given in Subsection \ref{43}. The limiting function has a Gaussian-type
structure. Unlike the real case, however,
for non-real $s$ it is not the characteristic
function of a probability distribution. To use this convergence in the analysis
of complex moments, we also require
uniform estimates for $H_{\sigma,s}$
over the entire Fourier space.
The following proposition provides
the necessary global decay.

\begin{prop}
\label{prop:global-Fourier-decay-rigorous}
Let \(K\) be a compact subset of the open right half-plane
$
\{s\in\mathbb C:\operatorname{Re}(s)>0\}.
$
For every $N>0$, there exist $\sigma_0>1/2$ and $C_{K,N}>0$ such that
\begin{equation}
\label{eq:global-Fourier-decay-rigorous}
|H_{\sigma,s}(\lambda,u)|
\leq
C_{K,N}(1+\lambda^2+u^2)^{-N}
\end{equation}
for all
\[
\frac12<\sigma\leq\sigma_0,
\qquad
s\in K,
\qquad
(\lambda,u)\in\mathbb R^2.
\]
Here \(H_{\sigma,s}(\lambda,u)\) is defined by
\eqref{eq:Hsigma-def-final}.
\end{prop}
Proof of Proposition \ref{prop:global-Fourier-decay-rigorous} is given in Subsection \ref{441}.
Together with the local derivative estimates
established in
Lemma~\ref{lem:local-Fourier-derivatives-rigorous},
these results provide the analytic control
needed to establish a Fourier--Riesz
representation of the complex moments.

\subsubsection{Removal of Gaussian regularization via Riesz distributions}

To recover the complex moments from
$H_{\sigma,s}$, we first introduce,
for $t>0$, the Gaussian-regularized moment
\begin{equation}
\label{eq:Gaussian-moment-introduction}
\mathcal M_{\sigma,t}(s)
:=
(2\sigma-1)^{s^2}
\mathbb E\left[
|\zeta(\sigma,X)|^{2s}
|T_\sigma|^{2s}
\e^{-t|T_\sigma|^2}
\right].
\end{equation}
By Lemma~\ref{lem:smoothed-moment-fourier},
this quantity admits a Fourier representation
in terms of $H_{\sigma,s}$.
For each fixed $t>0$,
Lemma~\ref{lem:fixed-t-smoothed-moment-limit}
establishes the convergence of
$\mathcal M_{\sigma,t}(s)$ as
$\sigma\downarrow1/2$.
Furthermore,
Lemma~\ref{lem:evaluation-fixed-t-limit}
gives the explicit formula
\begin{equation}
\label{eq:Mt-introduction}
\begin{aligned}
\mathcal M_t(s)
&:=
\lim_{\sigma\downarrow1/2}
\mathcal M_{\sigma,t}(s)=
a_s\e^{-s^2}\Gamma(s+1)
(1+t)^{-s-1}
{}_1F_1\left(
s+1;1;\frac{s^2}{1+t}
\right).
\end{aligned}
\end{equation}

We now explain how the Gaussian-regularized
moments are related to the complex moments
in Theorem~\ref{thm:complex-random-derivative-moments}.
Since
$\zeta'(\sigma,X)=\zeta(\sigma,X)L'(\sigma,X)$
and
$T_\sigma=(2\sigma-1)L'(\sigma,X)$,
we have
\begin{align}
\label{eq:original-moment-as-M0}
(2\sigma-1)^{s^2+2s}
\mathbb E\left[
|\zeta'(\sigma,X)|^{2s}
\right]=
(2\sigma-1)^{s^2}
\mathbb E\left[
|\zeta(\sigma,X)|^{2s}
|T_\sigma|^{2s}
\right]=:\mathcal M_{\sigma,0}(s).
\end{align}
Thus, the quantity appearing on the left-hand
side of our main theorem is precisely the
unregularized moment $\mathcal M_{\sigma,0}(s)$. On the other hand, letting $t\downarrow0$
in \eqref{eq:Mt-introduction} gives
\begin{align}
\label{eq:smoothed-limit-as-t-zero}
\lim_{t\downarrow0}
\lim_{\sigma\downarrow1/2}
\mathcal M_{\sigma,t}(s)
&=
a_s\e^{-s^2}\Gamma(s+1)
{}_1F_1(s+1;1;s^2),
\end{align}
which agrees with the right-hand side of
Theorem~\ref{thm:complex-random-derivative-moments}. To establish Theorem~\ref{thm:complex-random-derivative-moments}, however, we need
to prove that
\[
\lim_{\sigma\downarrow1/2}
\mathcal M_{\sigma,0}(s)=\lim_{\sigma\downarrow1/2}
\lim_{t\downarrow0}
\mathcal M_{\sigma,t}(s)
\]
exists and equals the expression in
\eqref{eq:smoothed-limit-as-t-zero}.
The main difficulty is therefore to remove
the Gaussian regularization while controlling
the limit as $\sigma\downarrow1/2$.

For real $s>0$, the uniform integrability
established in \eqref{uniformity}
provides the necessary control and permits
the interchange of the limits
$t\downarrow0$ and $\sigma\downarrow1/2$.
For non-real $s$, however, the same argument
does not apply directly.
Writing $s=a+\i b$ with $a>0$, we have
\[
|(2\sigma-1)^{s^2}|
=
(2\sigma-1)^{a^2-b^2}.
\]
Taking absolute values in the difference
between the regularized and unregularized
moments introduces an additional factor
$(2\sigma-1)^{-b^2}$ compared with the
corresponding estimates for real parameters. Consequently, for $b\neq0$, the available
uniform moment bounds do not directly justify
the interchange of these limits.

To overcome this difficulty, we use a
Fourier--Riesz approach that avoids
interchanging the two limits.
We first fix $\sigma>1/2$ and remove
the Gaussian regularization by letting
$t\downarrow0$.
We then represent the resulting complex
moment in terms of a regularized Riesz
functional of $H_{\sigma,s}$.
Finally, the Fourier asymptotics and decay
estimates established above allow us to
pass to the limit $\sigma\downarrow1/2$
in this representation. 
More precisely, we introduce the circular
average and a regularized Riesz functional
arising from the classical theory of
homogeneous distributions.

The circular average is the two-dimensional
specialization of the spherical mean used
by Gel'fand and Shilov in their study of
the homogeneous distribution $r^\lambda$;
see \cite[Chapter I, Section 3.9,
Eq.~(2), p.~71]{GelfandShilov}.
For a Schwartz function
$H\in\mathcal S(\mathbb R^2)$,
define its circular average by
\begin{equation}
\label{defofcirmean}
\mathcal A H(r)
:=
\frac{1}{2\pi}
\int_0^{2\pi}
H(r\cos\theta,r\sin\theta)
\,\d\theta,
\qquad r\geq0.
\end{equation}
To handle singular radial integrals at the origin,
we use the analytic continuation method
based on Taylor subtraction, developed
in \cite[Chapter I, Section 3.8,
Eq.~(2), pp.~66--67]{GelfandShilov},
together with the expansion of spherical
means given by Pizzetti's formula;
see \cite[Chapter I, Section 3.9,
Eq.~(7), p.~74]{GelfandShilov}. Let
\[
\Delta
:=
\frac{\partial^2}{\partial\lambda^2}
+
\frac{\partial^2}{\partial u^2}
\]
denote the Laplacian in the Fourier
variables $(\lambda,u)$. For $H\in\mathcal S(\mathbb R^2)$, $M\in\mathbb Z_{\geq0}$, and
\[
-1<\operatorname{Re}(w)<M+1,
\qquad
w\notin\{0,1,\ldots,M\},
\]
we define the regularized Riesz functional
\begin{align}
\mathcal R_{M,w}[H]
:={}&
\int_0^1
\left[
\mathcal A H(r)
-
\sum_{j=0}^M
\frac{r^{2j}}{4^j(j!)^2}
\Delta^jH(0)
\right]
r^{-2w-1}\,\d r
\nonumber\\
&+
\int_1^\infty
\mathcal A H(r)r^{-2w-1}\,\d r+
\sum_{j=0}^M
\frac{\Delta^jH(0)}
{4^j(j!)^2(2j-2w)}.
\label{eq:Riesz-functional-abstract}
\end{align}
Here $0$ denotes the origin in Fourier
space $\mathbb R^2$.
The Taylor subtraction removes the
singularity at the origin, while the
finite correction terms provide the
meromorphic continuation of the radial
integral. Thus, this functional gives
an explicit two-dimensional realization
of the regularization procedure for
homogeneous distributions described
in \cite[Chapter I, Sections 3.8--3.9]
{GelfandShilov}.

The Fourier--Riesz representation of the
complex moments is based on the classical
Fourier transform formula for the homogeneous
distribution $|x|^{2w}$; see
\cite[Chapter II, Section 3.3]
{GelfandShilov}.
Under the Fourier transform convention
\[
\widehat{f}(\xi)
:=
\int_{\mathbb R^2}
\e^{-\i\langle x,\xi\rangle}
f(x)\,\d x,
\]
this formula takes the form
\begin{equation}
\label{1023evening}
\widehat{|x|^{2w}}(\xi)
=
2^{2w+2}\pi
\frac{\Gamma(w+1)}{\Gamma(-w)}
\operatorname{Reg}|\xi|^{-2w-2},
\end{equation}
in the sense of tempered distributions,
where the right-hand side is understood
through meromorphic continuation.
For Schwartz functions $H$, the action
of the regularized distribution is given by
\[
\left\langle
\operatorname{Reg}|\xi|^{-2w-2},H
\right\rangle
=
2\pi\mathcal R_{M,w}[H].
\]

The Fourier--Riesz identity
\eqref{1023evening} provides the basis for
expressing the original complex moments
in terms of the regularized Riesz functional
$\mathcal R_{M,s}[H_{\sigma,s}]$.
Although the Fourier--Riesz identity is
initially formulated for Schwartz functions,
we will extend the corresponding
representation to a broader class of smooth
and integrable functions, allowing us to
apply it to the weighted Fourier transform
$H_{\sigma,s}$.

The following lemma establishes an exact
Fourier--Riesz representation of the complex
moments after removing the Gaussian
regularization for fixed $\sigma>1/2$.

\begin{lem}[Removal of Gaussian regularization]
\label{lem:moment-Riesz-Hsigma}
Let $M\in\mathbb Z_{\geq0}$, and let
$K$ be a compact subset of the open
right half-plane
$\{s\in\mathbb C:\operatorname{Re}(s)>0\}$.
There exists $\sigma_0>1/2$ such that,
for every
\[
\frac12<\sigma\leq\sigma_0,
\qquad s\in K\cap\Omega_M,
\]
where
\[
\Omega_M
:=
\left\{
w\in\mathbb C:
-1<\operatorname{Re}(w)<M+1
\right\}
\setminus\{0,1,\ldots,M\},
\]
one has
\begin{equation}
\label{eq:moment-Riesz-Hsigma}
\begin{aligned}
\lim_{t\downarrow0}
\mathcal M_{\sigma,t}(s)=
(2\sigma-1)^{s^2}
\mathbb E\left[
|\zeta(\sigma,X)|^{2s}
|T_\sigma|^{2s}
\right]=
2^{2s+1}
\frac{\Gamma(s+1)}{\Gamma(-s)}
\mathcal R_{M,s}[H_{\sigma,s}].
\end{aligned}
\end{equation}
Here $\mathcal M_{\sigma,t}(s)$,
$T_\sigma$, $H_{\sigma,s}$, and
$\mathcal R_{M,s}[H_{\sigma,s}]$
are defined in
\eqref{eq:Gaussian-moment-introduction},
\eqref{defofTsigma},
\eqref{eq:Hsigma-def-final}, and
\eqref{eq:Riesz-functional-abstract},
respectively.
\end{lem}

Lemma~\ref{lem:moment-Riesz-Hsigma}
is proved in Subsection~\ref{471}.
The first equality follows from dominated
convergence for fixed $\sigma>1/2$,
while the second is obtained using the
Fourier transform identity for homogeneous
distributions and the analytic continuation
of Riesz distributions. In particular, Lemma~\ref{lem:moment-Riesz-Hsigma}
expresses the original complex moments
exactly in terms of the regularized Riesz
functional $\mathcal R_{M,s}[H_{\sigma,s}]$
for each fixed $\sigma>1/2$.
This allows us to study their asymptotic
behaviour as $\sigma\downarrow1/2$
through the convergence of the corresponding
Riesz functionals, without having to justify
the interchange of the limits
$t\downarrow0$ and $\sigma\downarrow1/2$.
We then study the limit
$\sigma\downarrow1/2$ through this
Fourier--Riesz representation.
Lemma~\ref{lem:Riesz-functional-convergence}
establishes that
\[
\mathcal R_{M,s}[H_{\sigma,s}]
\longrightarrow
\mathcal R_{M,s}[H_s],
\]
locally uniformly away from the positive
integers. Here $H_s$ denotes the limiting function
of $H_{\sigma,s}$ as $\sigma\downarrow1/2$,
given explicitly in
\eqref{eq:compact-joint-transform-rigorous}.
Finally,
Lemma~\ref{lem:limiting-Riesz-evaluation}
gives the explicit evaluation
\[
\begin{aligned}
&2^{2s+1}
\frac{\Gamma(s+1)}{\Gamma(-s)}
\mathcal R_{M,s}[H_s]=
a_s\e^{-s^2}\Gamma(s+1)
{}_1F_1(s+1;1;s^2).
\end{aligned}
\]
The positive integer values of $s$ are
then treated by holomorphicity and the
maximum modulus principle, completing
the proof of
Theorem~\ref{thm:complex-random-derivative-moments}.

\medskip
\noindent\textbf{Use of AI.}
ChatGPT was used to assist in improving the clarity and
presentation of the manuscript, including reviewing the
presentation of mathematical arguments, checking references,
and identifying typographical errors.

\section{A probabilistic model for the derivative of the Riemann zeta function}

In this section, we construct a probabilistic
model for the derivative of the Riemann zeta
function using independent random variables
indexed by the primes. We then establish
the corresponding limiting distribution
and moment identities, and identify the
random model with the derivative of a
random Euler product.

\subsection{Construction of the random model}\label{section2}
\label{subsec:construction-random-model}
The classical Bohr--Jessen theorem describes the limiting value
distribution of $\log \zeta(\sigma+\mathrm{i}t)$ for fixed
$\sigma>1/2$. A natural probabilistic model for this limiting
distribution is obtained by replacing the phases $p^{-\mathrm{i}t}$
associated with distinct primes by independent random variables on the
unit circle. More precisely, let
\[
\Omega:=\prod_p\mathbb T,
\qquad
\mathbb T:=\{z\in\mathbb C:|z|=1\},
\]
and equip $\Omega$ with the product topology and the
corresponding Borel $\sigma$-algebra
\[
\mathcal F:=\bigotimes_p\mathcal B(\mathbb T),
\]
where $\mathcal B(\mathbb T)$ denotes the Borel
$\sigma$-algebra on $\mathbb T$. We endow $(\Omega,\mathcal F)$ with the probability measure
\begin{equation}
\label{eq:probability-measure-definition}
\mathbb P:=\bigotimes_p m_{\mathbb T},
\end{equation}
where $m_{\mathbb T}$ denotes the normalized Haar
measure on $\mathbb T$.
Then $\Omega$ is a compact abelian group and $\mathbb P$ is its
normalized Haar probability measure.
Recall that $\mathbb E$ denotes expectation with respect to $\mathbb P$.
For each prime $p$, define
\[
X(p)(\omega):=\omega_p,
\qquad
\omega=(\omega_p)_p\in\Omega.
\]
Then the random variables $\{X(p)\}_p$ are independent and
Haar-uniform on $\mathbb T$. We extend $X$ completely
multiplicatively to $\mathbb N$ by
\begin{equation}\label{defofX}
X(n):=\prod_{p^\nu\Vert n}X(p)^\nu.
\end{equation}
For fixed $\sigma>1/2$, we define
\begin{equation}\label{defofzetaderi}
Z(\sigma,X)
=
\lim_{N\to\infty}
\left(
-\sum_{n\le N}
\frac{\log n}{n^\sigma}X(n)
\right)
\qquad\text{in }L^2(\Omega,\mathbb P).
\end{equation}
This limit is well-defined since
\[
\mathbb{E}\Bigg[\Big|\sum_{M<n\leq N}\frac{\log n}{n^\sigma}X(n)\Big|^2\Bigg]
=
\sum_{M<n\leq N}
\frac{(\log n)^2}{n^{2\sigma}}
\leq
\sum_{n>M}
\frac{(\log n)^2}{n^{2\sigma}}
\longrightarrow0
\]
as $M\to\infty$, uniformly for $N>M$.

\subsection{Limiting value distribution and moment identity at fixed $\sigma$}
\label{subsec:random-model-complex-moments}
We shall show below that $Z(\sigma,X)$ is the limiting random
variable associated with $\zeta'(\sigma+\mathrm{i}t)$, and that,
under the Lindel\"of hypothesis, the corresponding complex moments
converge.

\begin{prop}
\label{prop:complex-moments-zeta-prime}
Assume the Lindel\"of hypothesis and let $\sigma>1/2$ be fixed. Let $Z(\sigma,X)$ be
defined as in \eqref{defofzetaderi}.
Then, for every $s\in\mathbb{C}$ with $\Re(s)>0$,
\begin{equation}
\label{eq:complex-moment-main}
\lim_{T\to\infty}
\frac{1}{T}
\int_1^T
\left|\zeta'(\sigma+\mathrm{i}t)\right|^{2s}\,dt
=
\mathbb{E}
[\left|Z(\sigma,X)\right|^{2s}].
\end{equation}
In particular, the limit on the left-hand side exists and is finite. Here, for $z\neq0$,
\[
|z|^{2s}:=\exp(2s\log|z|),
\]
and we set $0^{2s}:=0$.
\end{prop}

\begin{proof}
We first establish the limiting distribution of
$\zeta'(\sigma+\mathrm{i}t)$. For every fixed integer $n\geq1$, we claim that
\begin{equation}
\label{eq:zeta-fourier}
\lim_{T\to\infty}
\frac{1}{T}
\int_1^T\zeta(\sigma+\mathrm{i}t)n^{\mathrm{i}t}\d t
=
\frac{1}{n^\sigma}.
\end{equation}
If $\sigma>1$, this follows immediately from the absolutely
convergent Dirichlet series
\[
\zeta(\sigma+\mathrm{i}t)
=
\sum_{m=1}^{\infty}m^{-\sigma-\mathrm{i}t}.
\]
We therefore only need to consider $1/2<\sigma\leq 1$. Applying Cauchy's theorem to
\[
F(w):=\zeta(w)n^w
\]
on the rectangle with vertices
\[
\sigma+\mathrm{i},\quad
2+\mathrm{i},\quad
2+\mathrm{i}T,\quad
\sigma+\mathrm{i}T,
\]
and using the standard estimate (see e.g.,\cite[Chapter~V, section 5.1]{titchmarsh1986theory}),
\begin{equation}
\label{eq:zeta-growth-strip}
\sup_{\sigma\leq u\leq2}
|\zeta(u+\mathrm{i}T)|
=o(T),
\qquad T\to\infty,
\end{equation}
we may shift the vertical line of integration from $\Re(w)=\sigma$
to $\Re(w)=2$.  The two horizontal integrals, after division by $T$,
tend to zero, and hence
\begin{equation}
\label{eq:vertical_shift_average}
\frac{1}{T}\int_{1}^{T}
\zeta(\sigma+\mathrm{i}t)n^{\mathrm{i}t}\d t
-
n^{2-\sigma}\frac{1}{T}
\int_{1}^{T}
\zeta(2+\mathrm{i}t)n^{\mathrm{i}t}
\d t
\longrightarrow 0.
\end{equation}
The argument above on the line $\Re(w)=2$ gives
\[
\lim_{T\rightarrow \infty}n^{2-\sigma}
\int_{1}^{T}
\zeta(2+\mathrm{i}t)n^{\mathrm{i}t}\d t
=
n^{-\sigma}.
\]
and therefore \eqref{eq:zeta-fourier} holds for every $\sigma>1/2$. An integration by parts now gives
\begin{equation}\label{eq:zeta-prime-fourier}
\lim_{T\to\infty}
\frac{1}{T}\int_{1}^{T}
\zeta'(\sigma+\mathrm{i}t)n^{\mathrm{i}t}\d t
=
-\frac{\log n}{n^\sigma}.
\end{equation}
Indeed, if $g(t):=\zeta(\sigma+\mathrm{i}t)$ and
$\lambda:=\log n$, then
\[
g'(t)=\mathrm{i}\zeta'(\sigma+\mathrm{i}t),
\]
and hence
\[
\int_1^T
\zeta'(\sigma+\mathrm{i}t)e^{\mathrm{i}\lambda t}\,dt
=
-\mathrm{i}
\left[g(t)e^{\mathrm{i}\lambda t}\right]_1^T
-
\lambda
\int_1^T
g(t)\e^{\mathrm{i}\lambda t}\,dt.
\]
The boundary term divided by $T$ tends to zero by
\eqref{eq:zeta-growth-strip}, and \eqref{eq:zeta-prime-fourier}
follows from \eqref{eq:zeta-fourier}. Set
\[
a_n:=\frac{\log n}{n^\sigma},
\qquad
P_N(t):=
-\sum_{n\leq N}a_n n^{-\mathrm{i}t}.
\]
The case $k=1$ of the integer moment formula follows from the mean-value theorem for Dirichlet series, e.g., see \cite[Theorem~A and Equation (2)]{Ingham1928},
\begin{equation}
\label{eq:zeta_prime_second_moment}
\lim_{T\to\infty}\frac{1}{T}\int_{1}^{T}
|\zeta'(\sigma+\mathrm{i}t)|^2
\d t
=
\sum_{n=1}^{\infty}a_n^2.
\end{equation}
Notice that
\[
\sum_{n=1}^{\infty}a_n^2
=
\sum_{n=1}^{\infty}
\frac{(\log n)^2}{n^{2\sigma}}
<\infty,
\]
since $\sigma>1/2$. Since $P_N$ is a finite Dirichlet polynomial,
\begin{equation}
\label{eq:PN_second_moment}
\lim_{T\to\infty}\frac{1}{T}\int_{1}^{T}
|P_N(t)|^2\d t
=
\sum_{n\leq N}a_n^2.
\end{equation}
By \eqref{eq:zeta-prime-fourier}, 
\begin{align}
\lim_{T\to\infty}\frac{1}{T}
\int_{1}^{T}\zeta'(\sigma+\mathrm{i}t)
\overline{P_N(t)}\d t
=
\sum_{n\leq N}a_n^2.
\label{eq:cross_term}
\end{align}
Expanding the square and using
\eqref{eq:zeta_prime_second_moment}--\eqref{eq:cross_term},
we obtain
\begin{equation}
\label{eq:L2_approximation}
\lim_{T\to\infty}
\frac{1}{T}\int_{1}^{T}
\left|
\zeta'(\sigma+\mathrm{i}t)-P_N(t)
\right|^2\d t
=
\sum_{n>N}a_n^2.
\end{equation}
In particular,
\begin{equation}
\label{eq:L2_approximation_limit}
\lim_{N\to\infty}
\lim_{T\to\infty}\int_{1}^{T}
\left|
\zeta'(\sigma+\mathrm{i}t)-P_N(t)
\right|^2\d t
=0.
\end{equation} Define
\[
Y_{\sigma,N}
:=
-\sum_{n\leq N}a_n X(n),
\]
with $X(n)$ defined as in (\ref{defofX}).
Note that, for $M<N$,
\[
\mathbb{E}
|Y_{\sigma,N}-Y_{\sigma,M}|^2
=
\sum_{M<n\leq N}a_n^2.
\]
Since $\sum_n a_n^2<\infty$, the sequence
$\{Y_{\sigma,N}\}_{N\geq1}$ is Cauchy in $L^2$. By the definition of $Z(\sigma,X)$, so
\begin{equation}
\label{eq:Y_L2_convergence}
Y_{\sigma,N}\longrightarrow Z(\sigma,X)
\qquad\text{in }L^2.
\end{equation}
We next identify $Z(\sigma,X)$ as the limiting distribution of
$\zeta'(\sigma+\mathrm{i}t)$. Fix $N$, and let
$p_1,\ldots,p_r$ denote all primes not exceeding $N$.
Every $n\leq N$ can be written uniquely in the form
\[
n=p_1^{\nu_1(n)}\cdots p_r^{\nu_r(n)}.
\]
The numbers
\[
\log p_1,\ldots,\log p_r
\]
are linearly independent over $\mathbb{Q}$. Indeed, if
\[
m_1\log p_1+\cdots+m_r\log p_r=0,
\qquad
m_j\in\mathbb{Z},
\]
then
\[
p_1^{m_1}\cdots p_r^{m_r}=1,
\]
and unique factorization implies
$m_1=\cdots=m_r=0$. It follows from the continuous Kronecker-Weyl theorem that
\[
\left(
\e^{-\mathrm{i}t\log p_1},
\ldots,
\e^{-\mathrm{i}t\log p_r}
\right)
\]
is equidistributed on $\mathbb{T}^r$ with respect to Haar
measure. Hence, for every bounded continuous function
$\Phi:\mathbb{C}\to\mathbb{C}$,
\begin{equation}
\label{eq:PN_distribution}
\lim_{T\to\infty}\frac{1}{T}\int_{1}^{T}
\Phi(P_N(t))\d t
=
\mathbb{E}[\Phi(Y_{\sigma,N})].
\end{equation}
We now pass from $P_N$ to $\zeta'$. Let $\Phi$ be bounded
and Lipschitz, with Lipschitz constant
$\operatorname{Lip}(\Phi)$. By Cauchy--Schwarz and
\eqref{eq:L2_approximation},
\begin{align}
&\limsup_{T\to\infty}
\left|
\frac{1}{T}\int_{1}^{T}
\Phi(\zeta'(\sigma+\mathrm{i}t))\d t
-
\frac{1}{T}\int_{1}^{T}
\Phi(P_N(t))\d t\right|
\leq
\operatorname{Lip}(\Phi)
\left(
\sum_{n>N}a_n^2
\right)^{1/2}.
\label{eq:test_function_approximation}
\end{align}
On the other hand, by \eqref{eq:Y_L2_convergence},
\begin{align}
\left|
\mathbb{E}\left[\Phi(Y_{\sigma,N})\right]
-
\mathbb{E}\left[\Phi(Z(\sigma,X))\right]
\right|
&
\leq
\operatorname{Lip}(\Phi)
\|Y_{\sigma,N}-Y_\sigma\|_{L^2}
=
\operatorname{Lip}(\Phi)
\left(
\sum_{n>N}a_n^2
\right)^{1/2}.
\label{eq:random_model_test_function}
\end{align}
Combining
\eqref{eq:PN_distribution}--\eqref{eq:random_model_test_function},
first letting $T\to\infty$ and then $N\to\infty$, gives
\begin{equation}
\label{eq:weak_convergence_test}
\lim_{T\to\infty}
\frac{1}{T}\int_{1}^{T}
\Phi(\zeta'(\sigma+\mathrm{i}t))
\d t
=
\mathbb{E}\left[\Phi(Z(\sigma,X))\right]
\end{equation}
for every bounded Lipschitz function $\Phi$. Hence
\begin{equation}
\label{eq:weak_convergence}
\zeta'(\sigma+\mathrm{i}U_T)
\xrightarrow[T\to\infty]{d}
Z(\sigma,X)
\end{equation}
where $U_T$ is uniformly distributed on $[1,T]$. We now pass from bounded test functions to the complex moment.
Write
\[
a:=\Re(s)>0
\]
and choose an integer
\[
k>\max\{a,2\}.
\]
Under the Lindel\"of hypothesis, the integer moment estimate
\cite[Lemma~3.2]{simm2026moments} gives
\begin{equation}
\label{eq:high-integer-moment}
\lim_{T\to\infty}\frac{1}{T}\int_{1}^{T}
|\zeta'(\sigma+\mathrm{i}t)|^{2k}\d t
=
A_{k,\sigma}
<\infty.
\end{equation}
In particular, for all sufficiently large $T$,
\[
\frac{1}{T}\int_{1}^{T}
|\zeta'(\sigma+\mathrm{i}t)|^{2k}\d t
\leq C
\]
for some $C=C(k,\sigma)$. Since
\[
\left||z|^{2s}\right|=|z|^{2a},
\]
we have, on $|z|>R$,
\[
\left||z|^{2s}\right|
\leq
R^{-2(k-a)}|z|^{2k}.
\]
Hence, for all sufficiently large $T$,
\begin{equation}
\label{eq:tail-zeta-prime}
\frac{1}{T}\int_{1}^{T}\left||\zeta'(\sigma+\mathrm{i}t)|^{2s}\right|
\mathbf{1}_{\{|\zeta'(\sigma+\mathrm{i}t)|>R\}}
\leq
C R^{-2(k-a)}.
\end{equation}
We also need the corresponding bound for the limiting random
variable.  For $M>0$, set
\[
G_M(z):=\min\{|z|^{2k},M\}.
\]
By \eqref{eq:weak_convergence},
\[
\mathbb{E}\left[G_M(Z(\sigma,X))\right]
=
\lim_{T\to\infty}
\frac{1}{T}\int_{1}^{T}
G_M(\zeta'(\sigma+\mathrm{i}t))\d t
\leq A_{k,\sigma}.
\]
Letting $M\to\infty$ and using monotone convergence, we obtain
\begin{equation}
\label{eq:random-high-moment}
\mathbb{E}\left[|Z(\sigma,X)|^{2k}\right]
\leq A_{k,\sigma}<\infty.
\end{equation}
Consequently,
\begin{equation}
\label{eq:tail-random-zeta-prime}
\mathbb{E}
\left[
\left|
|Z(\sigma,X)|^{2s}
\right|
\mathbf{1}_{\{|Z(\sigma,X)|>R\}}
\right]
\leq
A_{k,\sigma}R^{-2(k-a)}.
\end{equation}
Choose a continuous function
$\chi:[0,\infty)\to[0,1]$ such that
\[
\chi(r)=1\quad (0\leq r\leq1),
\qquad
\chi(r)=0\quad (r\geq2),
\]
and define
\[
F_R(z)
:=
\chi(|z|/R)|z|^{2s}.
\]
Since $\Re(s)>0$, the function $z\mapsto |z|^{2s}$, with the value
$0$ assigned at $z=0$, is continuous on $\mathbb{C}$. So $F_R$
is bounded and continuous.  Hence, by \eqref{eq:weak_convergence},
\begin{equation}
\label{eq:cutoff_convergence}
\lim_{T\to\infty}\frac{1}{T}
\int_{1}^{T}
F_R(\zeta'(\sigma+\mathrm{i}t))\d t
=
\mathbb{E}\left[F_R(Z(\sigma,X))\right].
\end{equation}
It follows from (\ref{eq:tail-zeta-prime}) that
\begin{align}
&\limsup_{T\to\infty}
\left|
\frac{1}{T}\int_{1}^{T}
|\zeta'(\sigma+\mathrm{i}t)|^{2s}\d t
-
\frac{1}{T}\int_{1}^{T}
F_R(\zeta'(\sigma+\mathrm{i}t))\d t
\right|\leq
C R^{-2(k-a)}.
\label{eq:cutoff_error_deterministic}
\end{align}
Similarly, by \eqref{eq:random-high-moment},
\begin{equation}
\label{eq:cutoff_error_random}
\left|
\mathbb{E}\left[|Z(\sigma,X)|^{2s}\right]
-
\mathbb{E}\left[F_R(Z(\sigma,X))\right]
\right|
\leq
A_{k,\sigma}
R^{-2(k-a)}.
\end{equation}
Using \eqref{eq:cutoff_convergence}-
\eqref{eq:cutoff_error_random}, we obtain
\[
\begin{aligned}
&\limsup_{T\to\infty}
\left|
\frac{1}{T}\int_{1}^{T}
|\zeta'(\sigma+\mathrm{i}t)|^{2s}\d t
-
\mathbb{E}\left[|Z(\sigma,X)|^{2s}\right]
\right|\leq
(C+A_{k,\sigma})
R^{-2(k-a)}.
\end{aligned}
\]
Since $k>a$, letting $R\to\infty$ gives
\[
\lim_{T\to\infty}
\int_{1}^{T}
|\zeta'(\sigma+\mathrm{i}t)|^{2s}\d t
=
\mathbb{E}\left[|Z(\sigma,X)|^{2s}\right].
\]
By  \eqref{eq:random-high-moment}, the right-hand side is finite. 
\end{proof}

\begin{remark}
\label{rem:unconditional-complex-moments}
The Lindel\"of hypothesis in Proposition~\ref{prop:complex-moments-zeta-prime}
is only used to obtain sufficiently high integer moments of
$\zeta'(\sigma+\mathrm{i}t)$. In particular, the conclusion of
Proposition~\ref{prop:complex-moments-zeta-prime} holds unconditionally
for every fixed $\sigma>1/2$ and every $s\in\mathbb C$ satisfying
\[
0<\Re(s)<2.
\]
Indeed, write $a:=\Re(s)$. Since $a<2$, we may take $k=2$ in the
truncation argument.  By \cite[Lemma 3.1]{simm2026moments},  the fourth moment formula
\[
\lim_{T\to\infty}
\frac1T\int_1^T
|\zeta'(\sigma+\mathrm{i}t)|^4\,\d t
=
A_{2,\sigma}
<\infty
\]
holds unconditionally for every $\sigma>1/2$. Hence, for $R>0$,
\[
\frac1T\int_1^T
|\zeta'(\sigma+\mathrm{i}t)|^{2a}
\mathbf 1_{\{|\zeta'(\sigma+\mathrm{i}t)|>R\}}
\,\d t
\ll_{\sigma}
R^{-2(2-a)}
\]
for all sufficiently large $T$. Moreover, using the weak convergence established in the proof of
Proposition~\ref{prop:complex-moments-zeta-prime} and the bounded
continuous functions
\[
G_M(z):=\min\{|z|^4,M\},
\]
we obtain
\[
\mathbb E\left[|Z(\sigma,X)|^4\right]
\leq A_{2,\sigma}<\infty.
\]
Consequently,
\[
\mathbb E\left[
|Z(\sigma,X)|^{2a}
\mathbf 1_{\{|Z(\sigma,X)|>R\}}
\right]
\leq
A_{2,\sigma}R^{-2(2-a)}.
\]
The same cutoff argument as in the proof of
Proposition~\ref{prop:complex-moments-zeta-prime} therefore yields
\[
\lim_{T\to\infty}
\frac1T\int_1^T
|\zeta'(\sigma+\mathrm{i}t)|^{2s}\,\d t
=
\mathbb E\left[|Z(\sigma,X)|^{2s}\right]
\]
without assuming the Lindel\"of hypothesis.
\end{remark}
\subsection{Identification with the derivative of the random Euler product}
\label{subsec:identification-random-derivative}

Proposition~\ref{prop:complex-moments-zeta-prime}
expresses the limiting complex moments of
$\zeta'(\sigma+\mathrm{i}t)$ in terms of
the random variable $Z(\sigma,X)$
defined in \eqref{defofzetaderi}.
We now identify this random variable with
the derivative of the random Euler product
$\zeta(\sigma,X)$.
More precisely, in the following lemma \ref{lem:random-zeta-prime-log-derivative}, we prove that
\[
Z(\sigma,X)
=
\zeta'(\sigma,X)
=
\zeta(\sigma,X)L'(\sigma,X)
\qquad\text{almost surely}.
\]
This identification allows us to express
the limiting complex moments in terms of
the derivative of the random Euler product,  as stated in Theorem \ref{mainresult}.
\begin{lem}
\label{lem:random-zeta-prime-log-derivative}
Let $\{X(p)\}_{p}$ be independent random variables uniformly
distributed on $\mathbb T$, and extend $X$ completely multiplicatively
to the positive integers. For $\Re(w)>1/2$, let
\begin{equation}\label{defofrandomzeta}
\zeta(w,X)
:=
\prod_p
\left(1-\frac{X(p)}{p^w}\right)^{-1}
\end{equation}
denote the corresponding random Euler product, and define
\bea\label{defofL}
L(w,X)
:=
\sum_p
\log\left(1-\frac{X(p)}{p^w}\right)^{-1},
\eea
where the logarithm is defined by
\[
\log(1-z)^{-1}
=
\sum_{m=1}^{\infty}\frac{z^m}{m},
\qquad |z|<1.
\]
Then, almost surely, $L(\,\cdot\,,X)$ and $\zeta(\,\cdot\,,X)$ are
holomorphic on the half-plane $\Re(w)>1/2$ and
\[
\zeta(w,X)=\e^{L(w,X)}.
\]
Moreover, for every fixed $\sigma>1/2$, let $Z(\sigma,X)$ be
defined as in \eqref{defofzetaderi}. Then
\begin{equation}
\label{eq:random-zeta-prime-log-derivative}
Z(\sigma,X)
=
\zeta'(\sigma,X)
=\zeta(\sigma,X)\frac{\partial}{\partial\sigma}L(\sigma,X)
\qquad\text{almost surely}.
\end{equation}
\end{lem}
\begin{proof}
Fix $\sigma>1/2$, and choose $\tau$ such that
\[
\frac12<\tau<\sigma.
\]
By \cite[Proposition~3.2.9]{KowalskiPNT}, almost
surely, both
\[
\zeta(w,X)
=
\prod_p
\left(1-\frac{X(p)}{p^w}\right)^{-1}
\]
and
\[
\widetilde \zeta(w,X)
:=
\sum_{n=1}^{\infty}\frac{X(n)}{n^w}
\]
converge on the half-plane
\[
U_\tau:=\{w\in\mathbb C:\Re(w)>\tau\},
\]
locally uniformly on compact subsets, define holomorphic functions
there, and satisfy
\begin{equation}
\label{eq:random-euler-dirichlet-equality}
\zeta(w,X)=\widetilde \zeta(w,X)
\qquad (w\in U_\tau)
\end{equation}
almost surely. More precisely, if
\[
L_N(w,X)
:=
\sum_{p\leq N}
\log\left(1-\frac{X(p)}{p^w}\right)^{-1},
\]
then $L_N(\,\cdot\,,X)$ converges almost surely, locally uniformly on
$U_\tau$, to a holomorphic random function $L(\,\cdot\,,X)$.
Consequently,
\[
\exp(L_N(w,X))
=
\prod_{p\leq N}
\left(1-\frac{X(p)}{p^w}\right)^{-1}
\]
converges locally uniformly to $\e^{L(w,X)}$. Hence
\begin{equation}
\label{eq:Z-exp-L}
\zeta(w,X)=\e^{L(w,X)}
\end{equation}
almost surely on $U_\tau$. Since both sides of \eqref{eq:Z-exp-L} are holomorphic, the ordinary
chain rule gives
\begin{equation}
\label{eq:random-chain-rule}
\zeta'(w,X)
=
\e^{L(w,X)}L'(w,X)
\end{equation}
almost surely for every $w\in U_\tau$.

It remains to identify $\zeta'(\sigma,X)$ with the random variable
$Z(\sigma,X)$ defined through the $L^2$ limit. For $N\geq1$, set
\[
\zeta_N(w,X)
:=
\sum_{n\leq N}\frac{X(n)}{n^w}.
\]
By the almost sure local uniform convergence of the random Dirichlet
series,
\[
\zeta_N(\,\cdot\,,X)
\longrightarrow
\zeta(\,\cdot\,,X)
\]
locally uniformly on $U_\tau$ almost surely. Since the $\zeta_N$ are
holomorphic, the Weierstrass theorem, or equivalently Cauchy's
integral formula, implies that their derivatives converge locally
uniformly to the derivative of the limit. Note that $\sigma>\tau$. Thus
\[
\zeta_N'(\sigma,X)
\longrightarrow
\zeta'(\sigma,X)
\qquad\text{almost surely}.
\]But
\[
\zeta_N'(\sigma,X)
=
-\sum_{n\leq N}
\frac{\log n}{n^\sigma}X(n).
\]
Therefore
\begin{equation}
\label{eq:derivative-series-as}
-\sum_{n\leq N}
\frac{\log n}{n^\sigma}X(n)
\longrightarrow
\zeta'(\sigma,X)
\qquad\text{almost surely}.
\end{equation} On the other hand, using
\[
\mathbb E\!\left[X(n)\overline{X(m)}\right]
=
\delta_{n,m},
\]
we obtain, for $M<N$,
\[
\begin{aligned}
\mathbb E
\left[\left|
\zeta_N'(\sigma,X)-\zeta_M'(\sigma,X)
\right|^2\right]
&=
\sum_{M<n\leq N}
\frac{(\log n)^2}{n^{2\sigma}}.
\end{aligned}
\]
Since $\sigma>1/2$,
\[
\sum_{n=1}^{\infty}
\frac{(\log n)^2}{n^{2\sigma}}
<\infty.
\]
Hence $\{\zeta_N'(\sigma,X)\}$ is Cauchy in $L^2$, and by definition of $Z(\sigma,X)$, 
\[
\zeta_{N}'(\sigma,X)
\longrightarrow
Z(\sigma,X)
\qquad\text{in }L^2.
\]
In particular,
\begin{equation}
\label{eq:derivative-series-probability}
\zeta_N'(\sigma,X)
\longrightarrow
Z(\sigma,X)
\qquad\text{in probability}.
\end{equation}
By \eqref{eq:derivative-series-as}, 
\[
\zeta_N'(\sigma,X)
\longrightarrow
\zeta'(\sigma,X)
\qquad\text{almost surely},
\]
and hence also in probability. Since limits in probability are unique
up to almost sure equality, we conclude that
\begin{equation}
\label{eq:identify-random-derivative}
\zeta'(\sigma,X)
=
Z(\sigma,X)
\qquad\text{almost surely}.
\end{equation}
Combining \eqref{eq:identify-random-derivative} with
\eqref{eq:random-chain-rule}, evaluated at $w=\sigma$, gives
\[
Z(\sigma,X)
=
\zeta(\sigma,X)L'(\sigma,X)
\qquad\text{almost surely}.
\]
Finally, since $L(\,\cdot\,,X)$ is holomorphic, its complex derivative
at the real point $\sigma$ agrees with its derivative along the real
axis:
\[
L'(\sigma,X)
=
\frac{\partial}{\partial\sigma}L(\sigma,X).
\]
\end{proof}



We are now ready to prove
Theorem~\ref{mainresult}.

\begin{proof}[Proof of Theorem~\ref{mainresult}]
The assertion under the Lindel\"of hypothesis
follows directly from
Proposition~\ref{prop:complex-moments-zeta-prime}
and Lemma~\ref{lem:random-zeta-prime-log-derivative}.
The unconditional assertion for
$0<\operatorname{Re}(s)<2$
then follows from
Remark~\ref{rem:unconditional-complex-moments}.
\end{proof}

For the subsequent analysis of the moments
of $\zeta'(\sigma,X)$, we shall use
the prime-sum representation of the
logarithmic derivative $L'(\sigma,X)$
established in Lemma~\ref{lem:derivative-random-log} below.
\begin{lem}
\label{lem:derivative-random-log}
Let $\sigma>1/2$ be fixed. Then, almost surely,
\begin{equation}
\label{eq:derivative-random-log}
L'(\sigma,X)
=
-\lim_{N\to\infty}
\sum_{p\leq N}
\frac{(\log p)X(p)p^{-\sigma}}
     {1-X(p)p^{-\sigma}}.
\end{equation}
Equivalently,
\begin{equation}
\label{eq:derivative-random-log-series}
L'(\sigma,X)
=
-\sum_p
\frac{(\log p)X(p)}
     {p^\sigma-X(p)}
\qquad\text{almost surely},
\end{equation}
where the series on the right-hand side is understood as the almost surely convergent limit of its partial sums over the primes.
\end{lem}

\begin{proof}
Fix $\tau$ such that
\[
\frac12<\tau<\sigma,
\]
and set
\[
U_\tau:=\{w\in\mathbb C:\Re(w)>\tau\}.
\]
For $N\geq2$, define
\[
L_N(w,X)
:=
\sum_{p\leq N}
\log\left(1-\frac{X(p)}{p^w}\right)^{-1}.
\]
By the proof of \cite[Proposition~3.2.9]{KowalskiPNT}, there exists
an event $\Omega_0\subset\Omega$ with $\mathbb P(\Omega_0)=1$ such
that, for every $\omega\in\Omega_0$,
\[
L_N(\,\cdot\,,X(\omega))
\longrightarrow
L(\,\cdot\,,X(\omega))
\qquad\text{as }N\to\infty,
\]
locally uniformly on $U_\tau$. Fix $\omega\in\Omega_0$. Since each
$L_N(\,\cdot\,,X(\omega))$ is holomorphic on $U_\tau$ and the
convergence is locally uniform, the standard theorem on locally
uniform limits of holomorphic functions implies that
$L(\,\cdot\,,X(\omega))$ is holomorphic and
\[
L_N'(\,\cdot\,,X(\omega))
\longrightarrow
L'(\,\cdot\,,X(\omega)) \qquad\text{as }N\to\infty
\]
locally uniformly on $U_\tau$. In particular, since $\sigma\in U_\tau$,
\begin{equation}
\label{eq:LN-prime-limit}
L_N'(\sigma,X(\omega))
\longrightarrow
L'(\sigma,X(\omega)), \qquad\text{as }N\to\infty.
\end{equation}
For every fixed $N$, the sum defining $L_N$ is finite, and hence it
may be differentiated term by term. For each prime $p$, since
\[
\left|\frac{X(p)}{p^w}\right|<1,
\qquad w\in U_\tau,
\]
we have
\[
\log\left(1-\frac{X(p)}{p^w}\right)^{-1}
=
\sum_{m=1}^{\infty}
\frac{X(p)^m}{m p^{mw}}.
\]
Therefore,
\begin{align*}
\frac{\d}{\d w}
\log\left(1-\frac{X(p)}{p^w}\right)^{-1}
&=
-\sum_{m=1}^{\infty}
\frac{(\log p)X(p)^m}{p^{mw}}=
-(\log p)
\frac{X(p)p^{-w}}
     {1-X(p)p^{-w}}.
\end{align*}
Consequently,
\[
L_N'(\sigma,X(\omega))
=
-\sum_{p\leq N}
\frac{(\log p)X(p)p^{-\sigma}}
     {1-X(p)p^{-\sigma}}.
\]
Combining this identity with \eqref{eq:LN-prime-limit} yields
\[
L'(\sigma,X(\omega))
=
-\lim_{N\to\infty}
\sum_{p\leq N}
\frac{(\log p)X(p)p^{-\sigma}}
     {1-X(p)p^{-\sigma}}.
\]
Since $\omega\in\Omega_0$ was arbitrary and
$\mathbb P(\Omega_0)=1$, \eqref{eq:derivative-random-log} holds almost
surely.
\end{proof}

\section{Real moments via tilted measures}

In Section~\ref{section2}, we established
Theorem~\ref{mainresult}, which identifies,
under the Lindel\"of hypothesis, the complex
moments of the derivative of the random Euler
product with the corresponding limiting
moments of the derivative of the Riemann zeta
function for fixed $\sigma>1/2$.
We now investigate the asymptotic behaviour
of the former as $\sigma\downarrow1/2$.
We begin with real parameters $s>0$,
for which the positivity of the weight
$|\zeta(\sigma,X)|^{2s}$ allows us to
introduce a tilted probability measure.
.

\subsection{Construction of the tilted probability measure}
\label{sec:tilted-expectation}

In this subsection, we introduce a tilted probability measure associated
with the random Euler product. Throughout this section, let
\[
\sigma>\frac12,
\qquad
s>0,
\]
and recall that, almost surely,
\[
\zeta(\sigma,X)
=
\e^{L(\sigma,X)}
=
\prod_p
\left(1-\frac{X(p)}{p^\sigma}\right)^{-1}.
\]
Since $L(\sigma,X)$ is finite almost surely, it follows that
\begin{equation}
\label{eq:random-zeta-positive-finite}
0<|\zeta(\sigma,X)|<\infty
\qquad\text{almost surely}.
\end{equation}
We first verify that the moments required for the change of measure
are finite. For $N\geq2$, let
\[
\zeta_N(\sigma,X)
:=
\prod_{p\leq N}
\left(1-\frac{X(p)}{p^\sigma}\right)^{-1}.
\]
For $0<r<1$, set
\begin{equation}\label{defofmp2}
m_s(r)
:=
\mathbb E\left[
\left|1-rX(p)\right|^{-2s}
\right].
\end{equation}
Since $X(p)$ is Haar-uniform on $\mathbb T$, we have
\[
m_s(r)
=
\frac{1}{2\pi}
\int_0^{2\pi}
\left|1-r\e^{\mathrm{i}\theta}\right|^{-2s}
\,\d\theta.
\]
Using
\[
(1-z)^{-s}
=
\sum_{k=0}^{\infty}
\frac{(s)_k}{k!}z^k,
\qquad |z|<1,
\]
together with the orthogonality of the characters on $\mathbb T$,
we obtain
\begin{equation}
\label{eq:ms-expansion-final2}
m_s(r)
=
\sum_{k=0}^{\infty}
\left(\frac{(s)_k}{k!}\right)^2
r^{2k}.
\end{equation}
In particular, as $p\to\infty$,
\begin{equation}
\label{eq:local-random-zeta-moment-expansion}
m_s(p^{-\sigma})
=
1+s^2p^{-2\sigma}
+O_s(p^{-4\sigma}).
\end{equation}
Here, recall that $(r)_k$ denotes the Pochhammer symbol,
\[
(s)_0:=1,
\qquad
(s)_k:=s(s+1)\cdots(s+k-1),\quad k\ge1.
\]
Since $2\sigma>1$, it follows that
\[
\sum_p |m_s(p^{-\sigma})-1|<\infty,
\]
and therefore
\begin{equation}
\label{eq:local-moment-product-finite}
0<
\prod_p m_s(p^{-\sigma})
<\infty.
\end{equation}
By the independence of the random variables $\{X(p)\}_p$,
\begin{equation}
\label{eq:finite-random-euler-moment}
\mathbb E\left[
|\zeta_N(\sigma,X)|^{2s}
\right]
=
\prod_{p\leq N}m_s(p^{-\sigma}).
\end{equation}
Moreover, $\zeta_N(\sigma,X)\to\zeta(\sigma,X)$ almost surely as
$N\to\infty$. To pass to the limit in \eqref{eq:finite-random-euler-moment},
fix $q>1$. By \eqref{eq:local-moment-product-finite},
\[
\sup_{N\geq2}
\mathbb E\left[
|\zeta_N(\sigma,X)|^{2sq}
\right]
=
\sup_{N\geq2}
\prod_{p\leq N}m_{sq}(p^{-\sigma})
<\infty.
\]
Hence the family
\begin{equation}\label{UNIINtegrabilityofzetaN}
\left\{
|\zeta_N(\sigma,X)|^{2s}
\right\}_{N\geq2}
\end{equation}
is uniformly integrable. Since
\[
|\zeta_N(\sigma,X)|^{2s}
\longrightarrow
|\zeta(\sigma,X)|^{2s}
\qquad\text{almost surely},
\]
Vitali's convergence theorem yields
\begin{equation}
\label{eq:infinite-random-euler-moment}
\mathbb E\left[
|\zeta(\sigma,X)|^{2s}
\right]
=
\prod_p
\mathbb E\left[
\left|1-\frac{X(p)}{p^\sigma}\right|^{-2s}
\right]
\in(0,\infty).
\end{equation}
By \eqref{eq:infinite-random-euler-moment}, the denominator below is finite and strictly positive. We may therefore introduce the following well-defined tilted probability measure and the corresponding tilted expectation.

\begin{definition}
\label{def:tilted-expectation}
Let $\sigma>1/2$ and $s>0$. We define the tilted probability measure
$\widetilde{\mathbb P}_{\sigma,s}$ on $(\Omega,\mathcal F)$ by
\begin{equation}
\label{eq:tilted-measure-definition}
\frac{\d\widetilde{\mathbb P}_{\sigma,s}}{\d\mathbb P}
:=
\frac{
|\zeta(\sigma,X)|^{2s}
}{
\mathbb E\left[|\zeta(\sigma,X)|^{2s}\right]
}.
\end{equation}
Equivalently, for every $A\in\mathcal F$,
\[
\widetilde{\mathbb P}_{\sigma,s}(A)
=
\frac{
\mathbb E\left[
|\zeta(\sigma,X)|^{2s}\mathbf 1_A
\right]
}{
\mathbb E\left[
|\zeta(\sigma,X)|^{2s}
\right]
}.
\]
For every complex-valued random variable $F$ satisfying
\[
\mathbb E\left[
|\zeta(\sigma,X)|^{2s}|F|
\right]
<\infty,
\]
we define the corresponding tilted expectation by
\begin{equation}
\label{eq:tilted-expectation-definition}
\widetilde{\mathbb E}_{\sigma,s}[F]
:=
\int_\Omega F\,\d\widetilde{\mathbb P}_{\sigma,s}
=
\frac{
\mathbb E\left[
|\zeta(\sigma,X)|^{2s}F
\right]
}{
\mathbb E\left[
|\zeta(\sigma,X)|^{2s}
\right]
}.
\end{equation}
\end{definition}
In fact, the Radon-Nikodym derivative in
\eqref{eq:tilted-measure-definition} is nonnegative and satisfies
\[
\int_\Omega
\frac{
|\zeta(\sigma,X)|^{2s}
}{
\mathbb E[|\zeta(\sigma,X)|^{2s}]
}
\,\d\mathbb P
=1.
\]
Hence $\widetilde{\mathbb P}_{\sigma,s}$ is a probability measure.
Moreover, by \eqref{eq:random-zeta-positive-finite}, its
Radon--Nikodym derivative is strictly positive almost surely, so
$\widetilde{\mathbb P}_{\sigma,s}$ and $\mathbb P$ are equivalent
probability measures. Finally, \eqref{eq:tilted-expectation-definition}
is precisely the condition that
\[
F\in L^1(\Omega,\widetilde{\mathbb P}_{\sigma,s}).
\]
We finally record that the logarithmic derivative
$L'(\sigma,X)$ is integrable with respect to the tilted measure.
Recall from \eqref{eq:derivative-random-log-series} that
\begin{equation}
\label{eq:L-prime-prime-series}
L'(\sigma,X)
=
-\sum_p
\frac{(\log p)X(p)}
     {p^\sigma-X(p)}
\qquad\text{almost surely},
\end{equation}
where the series is understood as the limit of its partial sums over
the primes. Set
\[
Y_p
:=
-\frac{(\log p)X(p)}
       {p^\sigma-X(p)}
=
-(\log p)
\sum_{m=1}^{\infty}
\frac{X(p)^m}{p^{m\sigma}}.
\]
Since
\[
\mathbb E[X(p)^m]=0,
\qquad m\geq1,
\]
we have $\mathbb E[Y_p]=0$. Furthermore, the orthogonality of the
characters on $\mathbb T$ gives
\begin{align}
\mathbb E[|Y_p|^2]
&=
(\log p)^2
\sum_{m=1}^{\infty}p^{-2m\sigma}=
\frac{(\log p)^2}{p^{2\sigma}-1}.
\label{eq:Yp-second-moment}
\end{align}
Since the random variables $\{Y_p\}_p$ are independent and centered,
and
\[
\sum_p
\frac{(\log p)^2}{p^{2\sigma}-1}
<\infty,
\]
the partial sums in \eqref{eq:L-prime-prime-series} converge in
$L^2(\Omega,\mathbb P)$. Consequently,
\begin{equation}
\label{eq:L-prime-L2}
L'(\sigma,X)\in L^2(\Omega,\mathbb P),
\end{equation}
and in fact
\[
\mathbb E\left[
|L'(\sigma,X)|^2
\right]
=
\sum_p
\frac{(\log p)^2}{p^{2\sigma}-1}
<\infty.
\]
Combining \eqref{eq:infinite-random-euler-moment}, with $r=2s$,
and \eqref{eq:L-prime-L2}, the Cauchy--Schwarz inequality gives
\begin{align}
&\mathbb E\left[
|\zeta(\sigma,X)|^{2s}
|L'(\sigma,X)|
\right]
\leq
\left(
\mathbb E\left[
|\zeta(\sigma,X)|^{4s}
\right]
\right)^{1/2}
\left(
\mathbb E\left[
|L'(\sigma,X)|^2
\right]
\right)^{1/2}
<\infty.
\label{eq:L-prime-tilted-integrable}
\end{align}
Hence
\[
L'(\sigma,X)
\in
L^1(\Omega,\widetilde{\mathbb P}_{\sigma,s}),
\]
and therefore
\begin{equation}
\label{eq:tilted-expectation-L-prime}
\widetilde{\mathbb E}_{\sigma,s}
\left[L'(\sigma,X)\right]
=
\frac{
\mathbb E\left[
|\zeta(\sigma,X)|^{2s}
L'(\sigma,X)
\right]
}{
\mathbb E\left[
|\zeta(\sigma,X)|^{2s}
\right]
}
\end{equation}
is a well-defined finite complex number.

We remark that all the statements in this section hold unconditionally for every
$\sigma>1/2$ and $s>0$.

\subsection{Independence of prime-indexed random variables under the tilted measure}

Under the original probability measure
$\mathbb P$, the random variables
$\{X(p)\}_p$ are independent and uniformly
distributed on the unit circle.
We now investigate whether this independence
is preserved under the tilted probability
measure $\widetilde{\mathbb P}_{\sigma,s}$
defined in \eqref{eq:tilted-measure-definition}.
Lemma~\ref{lem:independence-under-tilted-measure}
below shows that the prime-indexed random
variables remain independent under the
tilted measure, although their marginal
distributions are modified by the
corresponding Euler factors.

\begin{lem}
\label{lem:independence-under-tilted-measure}
Let $\sigma>1/2$ and $s>0$. Let $m_{s}(p^{-\sigma})$ be defined in (\ref{defofmp2}). Under the tilted probability measure
$\widetilde{\mathbb P}_{\sigma,s}$ defined in
\eqref{def:tilted-expectation}, the random variables
$\{X(p)\}_p$ remain independent. More precisely, for each prime $p$,
the distribution of $X(p)$ under $\widetilde{\mathbb P}_{\sigma,s}$
is the probability measure
\begin{equation}
\label{eq:local-tilted-measure}
\d\widetilde{\mu}_{p,\sigma,s}(z)
:=
\frac{
|1-p^{-\sigma}z|^{-2s}
}{
m_s(p^{-\sigma})
}
\,\d m_{\mathbb T}(z),
\qquad z\in\mathbb T.
\end{equation}
Consequently,
\begin{equation}
\label{eq:tilted-product-measure}
\widetilde{\mathbb P}_{\sigma,s}
=
\bigotimes_p
\widetilde{\mu}_{p,\sigma,s}.
\end{equation}
\end{lem}

\begin{proof}
Since $p^{-\sigma}<1$, 
\[
0<m_s(p^{-\sigma})<\infty,
\]
and therefore \eqref{eq:local-tilted-measure} indeed defines a
probability measure on $\mathbb T$. For $N\geq2$, let
\[
\zeta_N(\sigma,X)
:=
\prod_{p\leq N}
\left(
1-\frac{X(p)}{p^\sigma}
\right)^{-1}
\]
and define
\[
R_N
:=
\frac{
|\zeta_N(\sigma,X)|^{2s}
}{
\mathbb E[|\zeta_N(\sigma,X)|^{2s}]
}.
\]
By the independence of the random variables $\{X(p)\}_p$ under
$\mathbb P$,
\[
\mathbb E[|\zeta_N(\sigma,X)|^{2s}]
=
\prod_{p\leq N}m_s(p^{-\sigma}),
\]
and hence
\begin{equation}
\label{eq:finite-tilted-density-product}
R_N
=
\prod_{p\leq N}
\frac{
|1-p^{-\sigma}X(p)|^{-2s}
}{
m_s(p^{-\sigma})
}.
\end{equation}
We first show that $R_N$ converges in $L^1(\Omega,\mathbb P)$ to the
Radon--Nikodym derivative defining
$\widetilde{\mathbb P}_{\sigma,s}$. Set
\[
Y_N:=|\zeta_N(\sigma,X)|^{2s},
\qquad
Y:=|\zeta(\sigma,X)|^{2s}.
\]
By  (\ref{UNIINtegrabilityofzetaN}),
\[
Y_N\longrightarrow Y
\qquad
\text{in }L^1(\Omega,\mathbb P),
\]
and therefore
\[
\mathbb E[Y_N]\longrightarrow \mathbb E[Y].
\]
Moreover, by \eqref{eq:infinite-random-euler-moment},
\[
0<\mathbb E[Y]<\infty.
\]
Thus, for all sufficiently large $N$, $M_N$ is bounded away from zero,
and
\begin{align*}
\mathbb E\left[
\left|
\frac{Y_N}{\mathbb E[Y_N]}
-
\frac{Y}{\mathbb E[Y]}
\right|
\right]
&\leq
\frac{1}{\mathbb E[Y_{N}]}
\mathbb E[|Y_N-Y|]
+
\left|
\frac{1}{\mathbb E[Y_{N}]}
-
\frac{1}{\mathbb E[Y]}
\right|
\mathbb E[Y]\longrightarrow0.
\end{align*}
Consequently,
\begin{equation}
\label{eq:RN-density-L1-convergence}
R_N
\longrightarrow
\frac{
|\zeta(\sigma,X)|^{2s}
}{
\mathbb E[|\zeta(\sigma,X)|^{2s}]
}
\qquad
\text{in }L^1(\Omega,\mathbb P).
\end{equation}
Now let $p_1,\ldots,p_k$ be distinct primes and let
$f_1,\ldots,f_k$ be bounded Borel measurable functions on $\mathbb T$.
Choose $N$ sufficiently large so that $p_j\leq N$ for all
$1\leq j\leq k$. Using \eqref{eq:finite-tilted-density-product} and
the independence of the variables $\{X(p)\}_p$ under $\mathbb P$, we
obtain
\begin{align*}
\mathbb E\Bigg[
R_N
\prod_{j=1}^k f_j(X(p_j))
\Bigg]=
\prod_{j=1}^k
\frac{
\mathbb E\left[
|1-p_j^{-\sigma}X(p_j)|^{-2s}
f_j(X(p_j))
\right]
}{
m_{s}(p_{j}^{-\sigma})
}
\prod_{\substack{p\leq N\\
p\notin\{p_1,\ldots,p_k\}}}
\frac{
\mathbb E\left[
|1-p^{-\sigma}X(p)|^{-2s}
\right]
}{
m_s(p^{-\sigma})
}.
\end{align*}
The second product on the right-hand side is equal to $1$. Hence
\[
\mathbb E\left[
R_N
\prod_{j=1}^k f_j(X(p_j))
\right]
=
\prod_{j=1}^k
\frac{
\mathbb E\left[
|1-p_j^{-\sigma}X(p_j)|^{-2s}
f_j(X(p_j))
\right]
}{
m_{s}(p_{j}^{-\sigma})
}.
\]
Since
\[
\prod_{j=1}^k f_j(X(p_j))
\]
is bounded, the $L^1$ convergence in
\eqref{eq:RN-density-L1-convergence} implies that
\[
\mathbb E\Bigg[
R_N
\prod_{j=1}^k f_j(X(p_j))
\Bigg]
\longrightarrow
\mathbb E\Bigg[
\frac{
|\zeta(\sigma,X)|^{2s}
}{
\mathbb E[|\zeta(\sigma,X)|^{2s}]}
\prod_{j=1}^k f_j(X(p_j))
\Bigg].
\]
By the definition of $\widetilde{\mathbb P}_{\sigma,s}$, we therefore
obtain
\begin{align}
\widetilde{\mathbb E}_{\sigma,s}
\Bigg[
\prod_{j=1}^k f_j(X(p_j))
\Bigg]
&=
\prod_{j=1}^k
\frac{
\mathbb E\Big[
|1-p_j^{-\sigma}X(p_j)|^{-2s}
f_j(X(p_j))
\Big]
}{
m_{s}(p_{j}^{-\sigma})
}
\nonumber\\
&=
\prod_{j=1}^k
\int_{\mathbb T}
f_j(z)\,
\d\widetilde{\mu}_{p_j,\sigma,s}(z).
\label{eq:infinite-dimensional-tilted-factorization}
\end{align}
Taking
\[
f_j=\mathbf 1_{A_j},
\qquad
A_j\in\mathcal B(\mathbb T),
\]
in \eqref{eq:infinite-dimensional-tilted-factorization}, we obtain
\begin{equation}
\label{eq:cylinder-factorization}
\widetilde{\mathbb P}_{\sigma,s}
\bigl(
X(p_1)\in A_1,\ldots,X(p_k)\in A_k
\bigr)
=
\prod_{j=1}^k
\widetilde{\mu}_{p_j,\sigma,s}(A_j).
\end{equation}
Taking $k=1$ shows that the marginal distribution of $X(p)$ under
$\widetilde{\mathbb P}_{\sigma,s}$ is
$\widetilde{\mu}_{p,\sigma,s}$. Since
\eqref{eq:cylinder-factorization} holds for every finite collection of
distinct primes and every collection of Borel sets
$A_1,\ldots,A_k\subset\mathbb T$, it follows directly from the
definition of independence that the random variables $\{X(p)\}_p$
are independent under $\widetilde{\mathbb P}_{\sigma,s}$. Finally, let
\[
\mu^\ast
:=
\bigotimes_p\widetilde{\mu}_{p,\sigma,s}.
\]
For every cylinder set of the form
\[
C
=
\left\{
\omega\in\Omega:
X(p_j)(\omega)\in A_j,\quad 1\leq j\leq k
\right\},
\]
with $A_j\in\mathcal B(\mathbb T)$, we have by
\eqref{eq:cylinder-factorization}
\[
\widetilde{\mathbb P}_{\sigma,s}(C)
=
\prod_{j=1}^k
\widetilde{\mu}_{p_j,\sigma,s}(A_j)
=
\mu^\ast(C).
\]
Since such cylinder sets form a $\pi$-system generating
\[
\mathcal F
=
\bigotimes_p\mathcal B(\mathbb T),
\]
the uniqueness theorem for probability measures yields
\[
\widetilde{\mathbb P}_{\sigma,s}
=
\mu^\ast
=
\bigotimes_p
\widetilde{\mu}_{p,\sigma,s}.
\]
\end{proof}
\subsection{Gaussian limit for the logarithmic derivative under the tilted measure}\label{33}

The independence established in
Lemma~\ref{lem:independence-under-tilted-measure}
allows us to analyze the normalized
logarithmic derivative
$(2\sigma-1)L'(\sigma,X)$
as a sum of independent random variables
and establish its Gaussian limit
as $\sigma\downarrow1/2$. We now prove Proposition~\ref{prop:tilted-gaussian-limit}.

\begin{proof}[Proof of Proposition~\ref{prop:tilted-gaussian-limit}]
It suffices to prove the assertion along an arbitrary sequence
$\sigma_n\downarrow1/2$. Set
\[
\delta_n:=2\sigma_n-1,
\]
so that $\delta_n\downarrow0$, and define
\[
W_{p,n}
:=
\frac{X(p)}{p^{\sigma_n}-X(p)}
=
\frac{p^{-\sigma_{n}}X(p)}
{1-p^{-\sigma_{n}}X(p)}.
\]
By Lemma~\ref{lem:independence-under-tilted-measure}, under
$\widetilde{\mathbb P}_{\sigma_n,s}$ the variables $\{X(p)\}_p$
are independent, and the law of $X(p)$ is
\[
\d\widetilde\mu_{p,\sigma_n,s}(z)
=
\frac{
|1-p^{-\sigma_{n}}z|^{-2s}
}{
m_s(p^{-\sigma_{n}})
}
\,\d m_{\mathbb T}(z),
\]
where $m_{s}(r)$ is defined in (\ref{defofmp2}).
We first establish the local moment estimates. By (\ref{eq:ms-expansion-final2}),
\begin{equation}
\label{eq:ms-expansion-final}
m_s(p^{-\sigma_{n}})
=
1+s^2p^{-2\sigma_{n}}+
\frac{s^2(s+1)^2}{4}p^{-4\sigma_{n}}
+O_s(p^{-6\sigma_{n}}).
\end{equation}
Since
\[
0\leq p^{-\sigma_{n}}\leq 2^{-1/2}<1,
\]
the above estimates may be taken uniformly in $p$ and $n$. Moreover, for $0\le r<1$,
\[
\frac{\d}{\d r}
\int_{\mathbb T}|1-rz|^{-2s}\,\d m_{\mathbb T}(z)
=
s\int_{\mathbb T}
|1-rz|^{-2s}
\left(
\frac{z}{1-rz}
+
\frac{\bar z}{1-r\bar z}
\right)
\,\d m_{\mathbb T}(z).
\]
The local tilted measure is invariant under complex conjugation, and
hence
\[
\widetilde{\mathbb E}_{\sigma_n,s}[W_{p,n}]
\in\mathbb R.
\]
It follows that
\[
\widetilde{\mathbb E}_{\sigma_n,s}[W_{p,n}]
=
\frac{p^{-\sigma_{n}}}{2s}
\frac{1}{m_s(p^{-\sigma_{n}})}\frac{\d}{\d r}
\int_{\mathbb T}|1-rz|^{-2s}\,\d m_{\mathbb T}(z)\Bigg|_{r=p^{-\sigma_{n}}}.
\]
By \eqref{eq:ms-expansion-final2}, we obtain
\begin{equation}
\label{eq:local-mean-final}
\widetilde{\mathbb E}_{\sigma_n,s}[W_{p,n}]
=
s p^{-2\sigma_{n}}+O_s(p^{-4\sigma_{n}}).
\end{equation}
Next, note that
\[
|W_{p,n}|^2
=
\frac{p^{-2\sigma_{n}}}
{|1-p^{-\sigma_{n}}X(p)|^2},
\]
so that
\[
\widetilde{\mathbb E}_{\sigma_n,s}
[|W_{p,n}|^2]
=
p^{-2\sigma_{n}}
\frac{m_{s+1}(p^{-\sigma_{n}})}
{m_s(p^{-\sigma_{n}})}.
\]
Hence
\begin{equation}
\label{eq:local-absolute-second-final}
\widetilde{\mathbb E}_{\sigma_n,s}
[|W_{p,n}|^2]
=p^{-2\sigma_{n}}+O_s(p^{-4\sigma_{n}}).
\end{equation}
For the second complex moment, expanding the factors into power
series gives
\[
\widetilde{\mathbb E}_{\sigma_n,s}[W_{p,n}^2]
=
\frac{p^{-4\sigma_{n}}}{m_s(p^{-\sigma_{n}})}
\sum_{k=0}^\infty
\frac{(s+2)_k(s)_{k+2}}
{k!(k+2)!}
r_{p,n}^{2k},
\]
and consequently
\begin{equation}
\label{eq:local-complex-second-final}
\widetilde{\mathbb E}_{\sigma_n,s}[W_{p,n}^2]
=
\frac{s(s+1)}{2}p^{-4\sigma_{n}}
+O_s(p^{-6\sigma_{n}}).
\end{equation}
By \eqref{eq:local-mean-final}--\eqref{eq:local-complex-second-final},
\begin{equation}
\label{eq:local-var-final}
\widetilde{\mathbb E}_{\sigma_n,s}
\left[
\left|W_{p,n}-\widetilde{\mathbb E}_{\sigma_n,s}[W_{p,n}]\right|^2
\right]
=
p^{-2\sigma_{n}}+O_s(p^{-4\sigma_{n}}),
\end{equation}
and
\begin{equation}
\label{eq:local-pseudovar-final}
\widetilde{\mathbb E}_{\sigma_n,s}
\left[
(W_{p,n}-\widetilde{\mathbb E}_{\sigma_n,s}[W_{p,n}])^2
\right]
=
O_s(p^{-4\sigma_{n}}).
\end{equation}
We next verify that the infinite prime sum is well defined. For every
fixed $n$, \eqref{eq:local-mean-final} gives
\[
\sum_p
(\log p)\left|\widetilde{\mathbb E}_{\sigma_n,s}[W_{p,n}]\right|
<\infty,
\]
since $2\sigma_n>1$. and define
\[
Z_{p,n}
:=
(\log p)(W_{p,n}-\widetilde{\mathbb E}_{\sigma_n,s}[W_{p,n}]).
\]
Likewise,
\eqref{eq:local-var-final} gives
\[
\sum_p
\widetilde{\mathbb E}_{\sigma_n,s}
[|Z_{p,n}|^2]
<\infty.
\]
The variables $\{Z_{p,n}\}_p$ are independent and centered under
$\widetilde{\mathbb P}_{\sigma_n,s}$, and therefore their partial
sums form a Cauchy sequence in
$L^2(\Omega,\widetilde{\mathbb P}_{\sigma_n,s})$. We may thus define
\begin{equation}
\label{eq:S-n-definition-final}
S_n
:=
\sum_p(\log p)W_{p,n}
=
\sum_p(\log p)\widetilde{\mathbb E}_{\sigma_n,s}[W_{p,n}]
+
\sum_pZ_{p,n},
\end{equation}
where the second sum is understood in
$L^2(\Omega,\widetilde{\mathbb P}_{\sigma_n,s})$. This definition agrees with the almost sure prime-series
representation of the logarithmic derivative. Indeed, by
\eqref{eq:derivative-random-log-series},
\[
-\sum_p
\frac{(\log p)X(p)}
{p^{\sigma_n}-X(p)}
=
L'(\sigma_n,X)
\qquad
\mathbb P\text{-almost surely}.
\]
Since $\widetilde{\mathbb P}_{\sigma_n,s}$ is equivalent to $\mathbb P$,
every $\mathbb P$-null set is also
$\widetilde{\mathbb P}_{\sigma_n,s}$-null. Hence the same identity holds
$\widetilde{\mathbb P}_{\sigma_n,s}$-almost surely.
On the other hand,
the partial sums in \eqref{eq:S-n-definition-final} converge in
$L^2(\widetilde{\mathbb P}_{\sigma_n,s})$, hence in probability, to
$S_n$. Uniqueness of limits in probability therefore gives
\begin{equation}
\label{eq:S-identification-final}
S_n=-L'(\sigma_n,X)
\qquad
\widetilde{\mathbb P}_{\sigma_n,s}\text{-almost surely}.
\end{equation} We now determine the asymptotic mean. From
\eqref{eq:local-mean-final},
\begin{align}
\widetilde{\mathbb E}_{\sigma_n,s}[S_n]
&=
s\sum_p
\frac{\log p}{p^{2\sigma_n}}
+
O_s\left(
\sum_p
\frac{\log p}{p^{4\sigma_n}}
\right)=
s\sum_p
\frac{\log p}{p^{1+\delta_n}}
+O_s(1),
\label{eq:mean-sum-final}
\end{align}
because $4\sigma_n\geq2$. It is known that
\begin{equation}
\label{eq:prime-first-final}
\sum_p
\frac{\log p}{p^{1+\delta}}
=
\frac1\delta+O(1).
\end{equation}
Then, it follows from \eqref{eq:mean-sum-final} that
\begin{equation}
\label{eq:mean-limit-final}
\delta_n
\widetilde{\mathbb E}_{\sigma_n,s}[S_n]
\longrightarrow s.
\end{equation}
We next determine the covariance of the centered sum. By
\eqref{eq:local-var-final},
\begin{align}
\sum_p
\widetilde{\mathbb E}_{\sigma_n,s}[|Z_{p,n}|^2]
=
\sum_p
\frac{(\log p)^2}{p^{2\sigma_n}}
+
O_s\left(
\sum_p
\frac{(\log p)^2}{p^{4\sigma_n}}
\right)
=
\sum_p
\frac{(\log p)^2}{p^{1+\delta_n}}
+O_s(1).
\label{eq:variance-sum-final}
\end{align}
By \eqref{eq:local-pseudovar-final},
\begin{equation*}
\label{eq:pseudovar-sum-final}
\sum_p
\widetilde{\mathbb E}_{\sigma_n,s}[Z_{p,n}^2]
=
O_s(1).
\end{equation*}
It is known that 
\begin{equation*}
\label{eq:prime-second-final}
\sum_p
\frac{(\log p)^2}{p^{1+\delta}}
=
\frac1{\delta^2}+O(1).
\end{equation*}
Consequently,
\begin{equation}
\label{eq:variance-limit-final}
\delta_n^2\left(\sum_p
\widetilde{\mathbb E}_{\sigma_n,s}[|Z_{p,n}|^2]\right)\longrightarrow1, \quad\text{as } n\to\infty,
\end{equation}
while
\begin{equation}
\label{eq:pseudovar-limit-final}
\delta_n^2\left(\sum_p
\widetilde{\mathbb E}_{\sigma_n,s}[Z_{p,n}^2]\right)\longrightarrow0, \quad\text{as } n\to\infty.
\end{equation}
We now reduce the centered infinite sum to a finite triangular array.
For each $n$, since
\[
\sum_p
\widetilde{\mathbb E}_{\sigma_n,s}[|Z_{p,n}|^2]
<\infty,
\]
we may choose $P_n\geq\delta_n^{-1}$ sufficiently large so that
\begin{equation*}
\label{eq:tail-choice-final}
\sum_{p>P_n}
\widetilde{\mathbb E}_{\sigma_n,s}
[|Z_{p,n}|^2]
\leq1.
\end{equation*}
In particular, $P_n\to\infty$, and
\begin{equation}
\label{eq:tail-var-final}
\delta_n^2
\sum_{p>P_n}
\widetilde{\mathbb E}_{\sigma_n,s}
[|Z_{p,n}|^2]
\longrightarrow0,
\qquad n\to\infty.
\end{equation}
Furthermore,
\begin{equation}
\label{eq:tail-pseudovar-final}
\delta_n^2
\left|
\sum_{p>P_n}
\widetilde{\mathbb E}_{\sigma_n,s}[Z_{p,n}^2]
\right|
\longrightarrow0.
\end{equation}
Consider the finite triangular array
\[
Y_{p,n}
:=
\delta_n
\begin{pmatrix}
\Re Z_{p,n}\\
\textnormal{Im}Z_{p,n}
\end{pmatrix},
\qquad
p\leq P_n.
\]
Its members are independent and centered. By
\eqref{eq:variance-limit-final},
\eqref{eq:pseudovar-limit-final},
\eqref{eq:tail-var-final}, and
\eqref{eq:tail-pseudovar-final}, the covariance matrix of
\[
\sum_{p\leq P_n}Y_{p,n}
\]
converges to $\frac12 \textbf{I}_2$, as $n\to\infty$. Indeed, for any centered complex random
variable $Z$,
\[
\Var(\Re Z)
=
\frac12
\left(
\mathbb E|Z|^2+\Re\mathbb E[Z^2]
\right), \quad \Var(\textnormal{Im}Z)
=
\frac12
\left(
\mathbb E|Z|^2-\Re\mathbb E[Z^2]
\right),
\]
and
\[
\textnormal{Cov}(\Re Z,\textnormal{Im}Z)
=
\frac12\textnormal{Im}\mathbb E[Z^2].
\]
It remains to verify the Lindeberg condition. Since
\[
|W_{p,n}|
\leq
\frac1{p^{\sigma_n}-1}
\leq
\frac1{\sqrt p-1},
\]
we have
\[
(\log p)|W_{p,n}|
\leq
\frac{\log p}{\sqrt p-1}.
\]
The function on the right is bounded on the set of primes. Thus there
exists an absolute constant $C>0$ such that
\[
(\log p)|W_{p,n}|\leq C
\]
for every $p$ and $n$. Therefore
\[
|Z_{p,n}|
\leq
(\log p)|W_{p,n}|
+
\widetilde{\mathbb E}_{\sigma_n,s}
[(\log p)|W_{p,n}|]
\leq2C.
\]
Consequently,
\[
\sup_{p\leq P_n}|Y_{p,n}|
\leq2C\delta_n
\longrightarrow0, \qquad n\to\infty.
\]
Hence, for every $\varepsilon>0$, for all sufficiently large $n$,
$
|Y_{p,n}|<\varepsilon$
for every $p\leq P_n$, and therefore the Lindeberg condition is satisfied. By the multivariate Lindeberg--Feller central limit theorem,
\begin{equation}
\label{eq:finite-clt-final}
\delta_n
\sum_{p\leq P_n}
\begin{pmatrix}
\Re Z_{p,n}\\
\text{Im} Z_{p,n}
\end{pmatrix}
\xrightarrow{\mathrm d}
N\left(
\begin{pmatrix}0\\0\end{pmatrix},
\frac12 \textbf{I}_2
\right).
\end{equation}
Finally, we shall remove the truncation. By independence and centering,
\begin{align*}
\widetilde{\mathbb E}_{\sigma_n,s}
\left[
\left|
\delta_n
\sum_{p>P_n}Z_{p,n}
\right|^2
\right]
&=
\delta_n^2
\sum_{p>P_n}
\widetilde{\mathbb E}_{\sigma_n,s}
[|Z_{p,n}|^2]\longrightarrow 0, \qquad n\to\infty
\end{align*}
by \eqref{eq:tail-var-final}. Hence
\[
\delta_n\sum_{p>P_n}Z_{p,n}
\longrightarrow0, \qquad n\to\infty
\]
in $L^2(\widetilde{\mathbb P}_{\sigma_n,s})$ and therefore in
probability. By \eqref{eq:finite-clt-final} and Slutsky's theorem,
\begin{equation}
\label{eq:centered-clt-final}
\delta_n\sum_pZ_{p,n}
\xrightarrow{\mathrm d}G, \qquad n\to\infty
\end{equation}
where
\[
G=G_1+\mathrm{i}G_2,
\qquad
G_1,G_2\stackrel{\mathrm{iid}}{\sim}N(0,1/2).
\]
Therefore, by \eqref{eq:S-n-definition-final},
\[
\delta_nS_n
=
\delta_n
\widetilde{\mathbb E}_{\sigma_n,s}[S_n]
+
\delta_n\sum_pZ_{p,n}.
\]
Combining \eqref{eq:S-identification-final}, \eqref{eq:mean-limit-final} and
\eqref{eq:centered-clt-final}, and applying Slutsky's theorem once
more, yields
\[
-(2\sigma_n-1)L'(\sigma_{n},X)
\xrightarrow{\mathrm d}
s+G, \qquad n\to\infty.
\]
Since the sequence $\sigma_n\downarrow1/2$ was arbitrary, the stated
convergence holds as $\sigma\downarrow1/2$.
\end{proof}

\subsection{Uniform integrability and asymptotics of real moments}

Proposition~\ref{prop:tilted-gaussian-limit}
establishes the convergence in distribution
of the normalized logarithmic derivative
$T_\sigma :=(2\sigma-1)L'(\sigma,X)$ under the tilted probability
measure as $\sigma\downarrow1/2$.
However, convergence in distribution alone
does not imply convergence of the corresponding
moments. To obtain the asymptotic formula
for the real moments of
$\zeta'(\sigma,X)$, we establish uniform
integrability of the relevant powers of
$T_\sigma$ under the tilted measure.
Combining this with the Gaussian limit
then yields the following result.
\begin{prop}
\label{prop:random-zeta-derivative-moment-asymptotics}
Let $s>0$, and let $L'(\sigma,X)$ be defined in
\eqref{eq:derivative-random-log-series}. Set
\[
T_\sigma=(2\sigma-1)L'(\sigma,X).
\]
There exist $q>2s$ and $\sigma_0>1/2$ such that
\begin{equation}
\label{eq:tilted-q-moment-bound-prop}
\sup_{1/2<\sigma\leq\sigma_0}
\widetilde{\mathbb E}_{\sigma,s}
\left[|T_\sigma|^q\right]<\infty,
\end{equation}
where $\widetilde{\mathbb E}_{\sigma,s}$ is defined in
\eqref{eq:tilted-expectation-definition}.
Consequently, as $\sigma\downarrow1/2$,
\begin{equation}
\label{eq:random-zeta-derivative-moment-asymptotic}
\mathbb E\left[|\zeta'(\sigma,X)|^{2s}\right]
\sim
a_s\e^{-s^2}\Gamma(s+1)
{}_1F_1(s+1;1;s^2)
(2\sigma-1)^{-(s^2+2s)},
\end{equation}
where $a_s$ is defined in
\eqref{definitionofthearithmeticfactor}.
\end{prop}

\begin{proof}
Set
\[
\delta:=2\sigma-1.
\]
By (\ref{eq:random-zeta-prime-log-derivative}),
\[
\zeta'(\sigma,X)
=
\zeta(\sigma,X)L'(\sigma,X)
\]
almost surely, the definition of the tilted expectation gives
\begin{equation}
\label{eq:real-moment-tilted-decomposition}
\begin{aligned}
\delta^{s^2+2s}
\mathbb E\left[|\zeta'(\sigma,X)|^{2s}\right]=
\left(
\delta^{s^2}
\mathbb E\left[|\zeta(\sigma,X)|^{2s}\right]
\right)
\widetilde{\mathbb E}_{\sigma,s}
\left[|\delta L'(\sigma,X)|^{2s}\right],
\end{aligned}
\end{equation}
where the identity is initially understood in the extended
nonnegative sense. We first show that
\begin{equation}
\label{eq:real-zeta-moment-limit}
\delta^{s^2}
\mathbb E\left[|\zeta(\sigma,X)|^{2s}\right]
\longrightarrow a_s
\qquad
(\sigma\downarrow1/2).
\end{equation}
By \eqref{eq:infinite-random-euler-moment},
\eqref{eq:ms-expansion-final2}, and the definition of
$m_s(p^{-\sigma})$ in \eqref{defofmp2}, we have
\[
\mathbb E\left[|\zeta(\sigma,X)|^{2s}\right]
=
\prod_p m_s(p^{-\sigma}),
\]
where
\[
m_s(p^{-\sigma})
=
\sum_{k=0}^{\infty}
\left(\frac{(s)_k}{k!}\right)^2
p^{-2k\sigma}.
\]
The local expansion
\[
m_s(p^{-\sigma})
=
1+s^2p^{-2\sigma}
+O_s(p^{-4\sigma})
\]
implies
\[
(1-p^{-2\sigma})^{s^2}m_s(p^{-\sigma})
=
1+O_s(p^{-4\sigma}),
\]
uniformly for $\sigma\geq1/2$. Since these factors are
strictly positive, their logarithms satisfy
\[
\log\left[
(1-p^{-2\sigma})^{s^2}m_s(p^{-\sigma})
\right]
=
O_s(p^{-4\sigma}).
\]
As
\[
\sum_p p^{-4\sigma}
\leq
\sum_p p^{-2}<\infty,
\]
the product
\[
\mathcal A_\sigma(s)
:=
\prod_p
(1-p^{-2\sigma})^{s^2}m_s(p^{-\sigma})
\]
converges, and the convergence is uniform for
$\sigma\geq1/2$. Hence
\[
\mathcal A_\sigma(s)
\longrightarrow
\prod_p
(1-p^{-1})^{s^2}
\sum_{k=0}^{\infty}
\left(\frac{(s)_k}{k!}\right)^2p^{-k}
=a_s.
\]
Moreover,
\[
\mathbb E\left[|\zeta(\sigma,X)|^{2s}\right]
=
\zeta(2\sigma)^{s^2}\mathcal A_\sigma(s).
\]
Since
\[
\delta\zeta(1+\delta)\longrightarrow1,
\]
we obtain \eqref{eq:real-zeta-moment-limit}. We next establish the uniform moment bound
\eqref{eq:tilted-q-moment-bound-prop}.
Choose
\[
q>\max\{2s,2\}.
\]
By Lemma~\ref{lem:independence-under-tilted-measure},
the variables $\{X(p)\}_p$ remain independent under
$\widetilde{\mathbb P}_{\sigma,s}$. Put
\[
W_{p,\sigma}
:=
\frac{X(p)}{p^\sigma-X(p)}
\]
and
\[
Z_{p,\sigma}
:=
(\log p)
\left(
W_{p,\sigma}
-
\widetilde{\mathbb E}_{\sigma,s}
[W_{p,\sigma}]
\right).
\]
The centered variables $Z_{p,\sigma}$ are independent,
and the prime-series representation of the logarithmic
derivative gives
\begin{equation}
\label{eq:real-centered-decomposition}
-L'(\sigma,X)
=
-\widetilde{\mathbb E}_{\sigma,s}
[L'(\sigma,X)]
+
\sum_p Z_{p,\sigma},
\end{equation}
where the centered series converges in
$L^2(\widetilde{\mathbb P}_{\sigma,s})$. 
By (\ref{eq:mean-limit-final}) and (\ref{eq:variance-sum-final}),
\begin{equation}
\label{eq:real-scaled-mean-bound}
-\delta\,
\widetilde{\mathbb E}_{\sigma,s}
[L'(\sigma,X)]
\longrightarrow s
\end{equation}
and
\begin{equation}
\label{eq:real-centered-variance-bound}
\sum_p
\widetilde{\mathbb E}_{\sigma,s}
[|Z_{p,\sigma}|^2]
=
O_s(\delta^{-2}).
\end{equation}
Furthermore, since $\sigma>1/2$,
\[
(\log p)|W_{p,\sigma}|
\leq
\frac{\log p}{\sqrt p-1}.
\]
The expression on the right is bounded uniformly
over all primes. Consequently,
\begin{equation}
\label{eq:real-centered-uniform-bound}
|Z_{p,\sigma}|\leq C_s
\end{equation}
for every prime $p$ and every $\sigma>1/2$. Applying Rosenthal's inequality to the real and
imaginary parts of the centered sum, we obtain
\begin{align}
\widetilde{\mathbb E}_{\sigma,s}
\left[
\left|\sum_{p\leq P}Z_{p,\sigma}\right|^q
\right]\ll_q
\left(
\sum_{p\leq P}
\widetilde{\mathbb E}_{\sigma,s}
[|Z_{p,\sigma}|^2]
\right)^{q/2}
+
\sum_{p\leq P}
\widetilde{\mathbb E}_{\sigma,s}
[|Z_{p,\sigma}|^q].
\label{eq:real-Rosenthal-estimate}
\end{align}
By \eqref{eq:real-centered-uniform-bound},
\[
|Z_{p,\sigma}|^q
\leq
C_s^{q-2}|Z_{p,\sigma}|^2.
\]
Thus \eqref{eq:real-centered-variance-bound} gives
\[
\sum_p
\widetilde{\mathbb E}_{\sigma,s}
[|Z_{p,\sigma}|^q]
=
O_{s,q}(\delta^{-2}).
\]
Since $q>2$, it follows from
\eqref{eq:real-Rosenthal-estimate} that
\begin{equation}
\label{eq:real-finite-q-moment}
\widetilde{\mathbb E}_{\sigma,s}
\left[
\left|\sum_{p\leq P}Z_{p,\sigma}\right|^q
\right]
\ll_{s,q}\delta^{-q},
\end{equation}
uniformly in $P$, for $\sigma$ sufficiently close
to $1/2$. For each fixed $\sigma>1/2$, the partial sums
converge in $L^2(\widetilde{\mathbb P}_{\sigma,s})$
and hence in probability. Passing to an almost
surely convergent subsequence and applying
Fatou's lemma to \eqref{eq:real-finite-q-moment},
we obtain
\begin{equation}
\label{eq:real-infinite-q-moment}
\widetilde{\mathbb E}_{\sigma,s}
\left[
\left|\sum_p Z_{p,\sigma}\right|^q
\right]
\ll_{s,q}\delta^{-q}.
\end{equation}
Combining \eqref{eq:real-centered-decomposition},
\eqref{eq:real-scaled-mean-bound}, and
\eqref{eq:real-infinite-q-moment}, we conclude
that there exists $\sigma_0>1/2$ such that
\[
\sup_{1/2<\sigma\leq\sigma_0}
\widetilde{\mathbb E}_{\sigma,s}
\left[
|\delta L'(\sigma,X)|^q
\right]
<\infty.
\]
This proves \eqref{eq:tilted-q-moment-bound-prop}.
Since $q>2s$, the preceding bound implies
that the family
\[
\left\{
|\delta L'(\sigma,X)|^{2s}:
\frac12<\sigma\leq\sigma_0
\right\}
\]
is uniformly integrable under the corresponding
tilted probability measures. On the other hand,
Proposition~\ref{prop:tilted-gaussian-limit}
shows that
\[
-\delta L'(\sigma,X)
\xrightarrow{\mathrm d}
s+G
\qquad
(\sigma\downarrow1/2)
\]
under $\widetilde{\mathbb P}_{\sigma,s}$.
By the continuous mapping theorem and uniform
integrability, we therefore obtain
\begin{equation}
\label{eq:real-tilted-moment-limit}
\widetilde{\mathbb E}_{\sigma,s}
\left[
|\delta L'(\sigma,X)|^{2s}
\right]
\longrightarrow
\mathbb E\left[|s+G|^{2s}\right].
\end{equation}
Combining \eqref{eq:real-moment-tilted-decomposition},
\eqref{eq:real-zeta-moment-limit}, and
\eqref{eq:real-tilted-moment-limit}, we obtain
\begin{equation}
\label{eq:real-moment-limit-gaussian}
\delta^{s^2+2s}
\mathbb E\left[|\zeta'(\sigma,X)|^{2s}\right]
\longrightarrow
a_s\mathbb E\left[|s+G|^{2s}\right].
\end{equation}
It remains to evaluate the Gaussian moment.
Since $G$ has density
\[
\frac1\pi\e^{-|z|^2}\,\d^2z
\]
on $\mathbb C$, translation and polar coordinates
give
\[
\mathbb E\left[|s+G|^{2s}\right]
=
2\e^{-s^2}
\int_0^\infty
r^{2s+1}\e^{-r^2}I_0(2sr)\,\d r,
\]
where $I_0$ is the modified Bessel function of
the first kind of order zero. Using its series
expansion
\[
I_0(2sr)
=
\sum_{k=0}^{\infty}
\frac{s^{2k}r^{2k}}{(k!)^2}
\]
and applying Tonelli's theorem, we obtain
\begin{align*}
\mathbb E\left[|s+G|^{2s}\right]
&=
2\e^{-s^2}
\sum_{k=0}^{\infty}
\frac{s^{2k}}{(k!)^2}
\int_0^\infty
r^{2s+2k+1}\e^{-r^2}\,\d r\\
&=
\e^{-s^2}
\sum_{k=0}^{\infty}
\frac{s^{2k}\Gamma(s+k+1)}{(k!)^2}\\
&=
\e^{-s^2}\Gamma(s+1)
{}_1F_1(s+1;1;s^2).
\end{align*}
Substituting this identity into
\eqref{eq:real-moment-limit-gaussian} proves
\eqref{eq:random-zeta-derivative-moment-asymptotic}
and completes the proof.
\end{proof}

\section{Complex moments via Fourier--Riesz analysis}
\label{sec:complex-moments-random-zeta-derivative}

In this section we extend the asymptotic formula for the moments of
the derivative of the random zeta function to complex moment
parameters. Throughout this section, let
$
\sigma>\frac12$.
Recall that
\[
L'(\sigma,X)
=
-\sum_p
\frac{(\log p)X(p)}{p^\sigma-X(p)},
\]
so that
\[
\zeta'(\sigma,X)
=
\zeta(\sigma,X)L'(\sigma,X)
\]
almost surely. For $x>0$ and $s\in\mathbb C$, we use the convention
\[
x^s:=\exp(s\log x),
\]
where $\log x$ denotes the real logarithm. If
$\operatorname{Re}(s)>0$, we set $0^s:=0$. Since $2\sigma-1>0$, we
also define
\[
(2\sigma-1)^s:=\exp\big(s\log(2\sigma-1)\big).
\]
In particular,
\[
(2\sigma-1)^{s^2+2s}
=
\exp\left((s^2+2s)\log(2\sigma-1)\right).
\]
Recall the arithmetic factor
\begin{equation*}
a_s
:=
\prod_p
(1-p^{-1})^{s^2}
\sum_{m=0}^{\infty}
\left(\frac{(s)_m}{m!}\right)^2p^{-m}.
\end{equation*}

Our main result in this section is the following
asymptotic formula for complex moments of the
derivative of the random zeta function.

\begin{thm}[Restatement of Theorem~\ref{thm:complex-random-derivative-moments}]
For every $s\in\mathbb C$ with
$\operatorname{Re}(s)>0$, one has
\begin{equation}
\lim_{\sigma\downarrow1/2}
(2\sigma-1)^{s^2+2s}
\mathbb E\left[
|\zeta'(\sigma,X)|^{2s}
\right]
=
a_s \mathrm e^{-s^2}\Gamma(s+1)
{}_1F_1(s+1,1;s^2),
\end{equation}
where $a_s$ is the arithmetic factor defined
in \eqref{definitionofthearithmeticfactor},
and ${}_1F_1(a,b;z)$ denotes the confluent
hypergeometric function of the first kind.
Moreover, the convergence is locally uniform
in the half-plane $\operatorname{Re}(s)>0$.
\end{thm}

Combining Theorem~\ref{thm:complex-random-derivative-moments}
with Theorem~\ref{mainresult}, which identifies
the mean values of complex moments of
$\zeta'(\sigma+\mathrm{i}t)$ with the corresponding
moments of the random model for each fixed
$\sigma>1/2$, we obtain Theorem~\ref{onemainresult}.

\subsection{Proof of Corollary~\ref{onemainresult}
assuming Theorem~\ref{thm:complex-random-derivative-moments}}\label{41}
\begin{proof}[Deduction of Corollary~\ref{onemainresult}
from Theorem~\ref{thm:complex-random-derivative-moments}]
Assume the Lindel\"of hypothesis.
By Theorem~\ref{mainresult}, for every fixed
$\sigma>1/2$ and every $s\in\mathbb C$ with
$\operatorname{Re}(s)>0$, we have
\begin{equation}
\label{eq:complex-moment-main2}
\lim_{T\to\infty}
\frac{1}{T}
\int_1^T
|\zeta'(\sigma+\mathrm{i}t)|^{2s}\,\d t
=
\mathbb E\left[
|\zeta'(\sigma,X)|^{2s}
\right].
\end{equation}
Consequently, Theorem~
\ref{thm:complex-random-derivative-moments}
gives
\begin{align*}
&\lim_{\sigma\downarrow1/2}
(2\sigma-1)^{s^2+2s}
\left(
\lim_{T\to\infty}
\frac{1}{T}
\int_1^T
|\zeta'(\sigma+\mathrm{i}t)|^{2s}\,\d t
\right)\\
&\qquad=
a_s\e^{-s^2}\Gamma(s+1)
{}_1F_1(s+1;1;s^2).
\end{align*}
Moreover, since \eqref{eq:complex-moment-main2}
holds unconditionally for
$0<\operatorname{Re}(s)<2$,
the same asymptotic formula holds
unconditionally in this range. For $s=2$, the result follows directly from
\cite[Theorem 1.7]{simm2026moments}. 
This completes the proof.
\end{proof}

It therefore remains to prove
Theorem~\ref{thm:complex-random-derivative-moments},
whose proof is given in Subsection~\ref{thm41}.
We begin by establishing several preliminary estimates
that will be used in the proof.

\subsection{Preliminary estimates}

We first show the convergence and integrability properties of the
random Euler product which will be used below.

\begin{lem}
\label{lem:basic-random-Euler-product}
Let $\sigma>1/2$. Define
\bea\label{defineofzetaP}
\zeta_P(\sigma,X)
=
\prod_{p\leq P}
(1-p^{-\sigma}X(p))^{-1}
\eea
and
\begin{equation}\label{defofLP}
L'_P(\sigma,X)
=
-\sum_{p\leq P}
\frac{(\log p)X(p)}{p^\sigma-X(p)}.
\end{equation}
Then
\[
\zeta_P(\sigma,X)\longrightarrow\zeta(\sigma,X)
\]
and
\[
L'_P(\sigma,X)\longrightarrow L'(\sigma,X)
\]
almost surely as $P\to\infty$. Moreover, for every $q>0$,
\begin{equation}
\label{eq:fixed-sigma-all-moments}
\mathbb E\left[|\zeta(\sigma,X)|^q\right]<\infty,
\qquad
\mathbb E\left[|L'(\sigma,X)|^q\right]<\infty,
\qquad
\mathbb E\left[|\zeta'(\sigma,X)|^q\right]<\infty.
\end{equation}
\end{lem}

\begin{proof}
Write
\[
\log\zeta_P(\sigma,X)
=
\sum_{p\leq P}
\sum_{k\geq1}
\frac{X(p)^k}{kp^{k\sigma}}.
\]
The contribution from $k\geq2$ converges absolutely, since
\[
\sum_p\sum_{k\geq2}\frac1{kp^{k\sigma}}<\infty.
\]
The random variables $X(p)p^{-\sigma}$ are independent and centered,
and
\[
\sum_p
\mathbb E\left[\left|
\frac{X(p)}{p^\sigma}
\right|^2\right]
=
\sum_pp^{-2\sigma}<\infty.
\]
Hence Kolmogorov's convergence theorem gives almost sure convergence
of
\[
\sum_p\frac{X(p)}{p^\sigma}.
\]
Thus $\log\zeta_P(\sigma,X)$ converges almost surely, and therefore so
does $\zeta_P(\sigma,X)$.
For the logarithmic derivative, put
\[
Z_{p,\sigma}
:=
\frac{(\log p)X(p)}{p^\sigma-X(p)}
=
(\log p)
\sum_{k\geq1}\frac{X(p)^k}{p^{k\sigma}}.
\]
Then
\[
\mathbb E\left[Z_{p,\sigma}\right]=0
\]
and
\[
\mathbb E\left[|Z_{p,\sigma}|^2\right]
\ll_\sigma
\frac{(\log p)^2}{p^{2\sigma}}.
\]
Hence
\[
\sum_p\mathbb E\left[|Z_{p,\sigma}|^2\right]<\infty,
\]
and Kolmogorov's convergence theorem again gives almost sure
convergence of $\sum_p Z_{p,\sigma}$.
For the moment estimates, independence gives, for the partial Euler
products,
\[
\mathbb E[|\zeta_P(\sigma,X)|^q]
=
\prod_{p\leq P}
\int_{\mathbb T}
|1-p^{-\sigma}z|^{-q}\,\d m_{\mathbb T}(z).
\]
For any $q>0$, the local factors satisfy
\[
\int_{\mathbb T}
|1-p^{-\sigma}z|^{-q}\,dm_{\mathbb T}(z)
=
1+O_{q,\sigma}(p^{-2\sigma}),
\]
and hence
\[
\sup_P
\mathbb E\left[|\zeta_P(\sigma,X)|^q\right]<\infty.
\]
Fatou's lemma therefore gives
\[
\mathbb E\left[|\zeta(\sigma,X)|^q\right]<\infty.
\]
For $q\geq2$,
\[
\sum_p\mathbb E\left[|Z_{p,\sigma}|^q\right]
\ll_{q,\sigma}
\sum_p\frac{(\log p)^q}{p^{q\sigma}}
<\infty.
\]
Rosenthal's inequality, first applied to finite partial sums and then
combined with Fatou's lemma, gives
\[
\mathbb E\left[|L'(\sigma,X)|^q\right]<\infty.
\]
The case $0<q<2$ follows from the $L^2$ bound. Hence, we have for any $q>0$, $
\mathbb E\left[|L'(\sigma,X)|^q\right]<\infty$.
Finally, Using 
\[
|\zeta'(\sigma,X)|
=
|\zeta(\sigma,X)|\,|L'(\sigma,X)|,
\]
and H\"older's inequality, we have $\mathbb{E}\left[|\zeta'(\sigma,X)|^{q}\right]<\infty$ for any $q>0$.
\end{proof}
We now introduce a local moment function that
extends the quantity $m_s(r)$ defined in
\eqref{defofmp2} to complex values of $s$.
\begin{lem}
\label{lem:ms-complex}
For $0\leq r<1$ and $s\in\mathbb C$, define
\begin{equation}\label{defofms}
m_s(r)
:=
\int_{\mathbb T}
|1-rz|^{-2s}\,\d m_{\mathbb T}(z).
\end{equation}
Then, 
for $0\leq r<1$,
\begin{equation}
\label{eq:ms-power-series-final}
m_s(r)
=
\sum_{m=0}^{\infty}
\left(\frac{(s)_m}{m!}\right)^2r^{2m}.
\end{equation}
Moreover, for every compact set $K\subset\mathbb C$,
\begin{equation}
\label{eq:ms-expansion-final}
m_s(r)
=
1+s^2r^2+O_K(r^4)
\end{equation}
uniformly for $s\in K$ and $0\leq r\leq2^{-1/2}$.
\end{lem}

\begin{proof}
Since $r<1$, the number $1-rz$ lies in the open right half-plane
for $z\in\mathbb T$. Hence the principal logarithm may be used and
\[
|1-rz|^{-2s}
=
(1-rz)^{-s}(1-r\overline z)^{-s}.
\]
Using
\[
(1-rz)^{-s}
=
\sum_{m\geq0}\frac{(s)_m}{m!}r^mz^m
\]
and Haar orthogonality gives \eqref{eq:ms-power-series-final}.
The expansion \eqref{eq:ms-expansion-final} follows immediately.
\end{proof}
We shall use the following standard estimates
for sums over primes, which follow from the
classical properties of the Riemann zeta
function near its pole at $1$.

\begin{lem}
\label{lem:prime-sums-final}
As \(\delta\downarrow0\), the following asymptotic estimates hold,
where all sums are taken over primes:
\begin{align}
\sum_p \frac{1}{p^{1+\delta}}
&=
\log\frac{1}{\delta}+O(1),
\label{eq:prime-sum-zero-final}
\\
\sum_p \frac{\log p}{p^{1+\delta}}
&=
\frac{1}{\delta}+O(1),
\label{eq:prime-sum-one-final}
\\
\sum_p \frac{(\log p)^2}{p^{1+\delta}}
&=
\frac{1}{\delta^2}+O(1).
\label{eq:prime-sum-two-final}
\end{align}
These asymptotic estimates remain valid if any fixed finite set of
primes is omitted from the sums.
\end{lem}
\begin{proof}
The first estimate follows from
\[
\log\zeta(1+\delta)
=
\sum_p\sum_{k\geq1}
\frac1{kp^{k(1+\delta)}}
=
\log\frac1\delta+O(1),
\]
since the contribution from $k\geq2$ is uniformly bounded. Similarly,
\[
-\frac{\zeta'}{\zeta}(1+\delta)
=
\sum_p\sum_{k\geq1}
\frac{\log p}{p^{k(1+\delta)}}
=
\frac1\delta+O(1),
\]
which gives \eqref{eq:prime-sum-one-final}. Differentiating with
respect to $\delta$ yields
\[
\sum_p\sum_{k\geq1}
\frac{k(\log p)^2}{p^{k(1+\delta)}}
=
\frac1{\delta^2}+O(1),
\]
and the $k\geq2$ terms are again uniformly bounded.
\end{proof}

We also require the following uniform tail estimate.

\begin{lem}
\label{lem:uniform-prime-tail-final}
There exist $P_0\geq2$ and $c>0$ such that, uniformly for
$P\geq P_0$ and $0<\delta\leq1/2$,
\begin{equation}
\label{eq:uniform-prime-tail-final}
\sum_{p\geq P}
\frac{(\log p)^2}{p^{1+\delta}}
\geq
cP^{-\delta}
\left(
\frac{\log P}{\delta}
+
\frac1{\delta^2}
\right).
\end{equation}
\end{lem}

\begin{proof}
Let
\[
\vartheta(x):=\sum_{p\leq x}\log p.
\]
By the prime number theorem (see \cite[Chapter~4, Theorem~4.4]{Apostol1976}), after increasing $P_0$ if necessary,
there exists $c_0>0$ such that
\[
\vartheta(2x)-\vartheta(x)\geq c_0x,
\qquad x\geq P_0.
\]
Choosing a $P$ with $P\geq P_{0}$. 
Put $P_k=2^kP$. Then
\begin{align*}
\sum_{P_k<p\leq2P_k}
\frac{(\log p)^2}{p^{1+\delta}}
&\geq
(\log P_k)(2P_k)^{-1-\delta}
\bigl(\vartheta(2P_k)-\vartheta(P_k)\bigr)\\
&\gg
(\log P_k)P_k^{-\delta}.
\end{align*}
Here and below, for nonnegative quantities $f$ and $g$, we use the notation \(f\ll g\) (equivalently, \(g\gg f\)) to mean
that \(f\leq cg\) for some positive constant \(c\), and
\(f\asymp g\) to mean that both \(f\ll g\) and \(g\ll f\).
Summing over $k\geq0$ gives
\[
\sum_{p\geq P}
\frac{(\log p)^2}{p^{1+\delta}}
\gg
P^{-\delta}
\sum_{k=0}^{\infty}
(\log P+k\log2)\e^{-k\delta\log2}.
\]
Since, uniformly for $0<\delta\leq1/2$,
\[
\sum_{k\geq0}\e^{-k\delta\log2}\asymp\delta^{-1},
\qquad
\sum_{k\geq0}k\e^{-k\delta\log2}\asymp\delta^{-2},
\]
the assertion follows.  
\end{proof}
Finally, we record the real Euler-product bound required later.
\begin{lem}
\label{lem:real-zeta-moment-uniform}
Let $0<a_0\leq a\leq a_1$.
Uniformly for $a\in[a_0,a_1]$ and $\sigma$ sufficiently close to
$1/2$,
\begin{equation}
\label{eq:real-Euler-product-uniform}
(2\sigma-1)^{a^2}
\prod_p m_a(p^{-\sigma})
\ll_{a_0,a_1}1.
\end{equation}
Here, $m_{a}(r)$ is defined by (\ref{defofms}).
Consequently,
\begin{equation}
\label{eq:real-zeta-moment-uniform}
(2\sigma-1)^{a^2}
\mathbb E\left[|\zeta(\sigma,X)|^{2a}\right]
\ll_{a_0,a_1}1.
\end{equation}
\end{lem}

\begin{proof}
By \eqref{eq:ms-expansion-final},
\[
(1-p^{-2\sigma})^{a^2}
m_a(p^{-\sigma})
=
1+O_{a_0,a_1}(p^{-4\sigma}).
\]
Therefore
\[
A_\sigma(a)
:=
\prod_p
(1-p^{-2\sigma})^{a^2}
m_a(p^{-\sigma})
\]
is uniformly bounded for $\sigma$ sufficiently close to $1/2$.
Since
\[
\prod_p(1-p^{-2\sigma})^{-a^2}
=
\zeta(2\sigma)^{a^2},
\]
we obtain
\[
(2\sigma-1)^{a^2}\prod_pm_a(p^{-\sigma})
=
\bigl((2\sigma-1)\zeta(2\sigma)\bigr)^{a^2}
A_\sigma(a)
\ll1.
\]
This proves \eqref{eq:real-Euler-product-uniform}.
For the partial Euler products,
\[
\mathbb E\left[|\zeta_P(\sigma,X)|^{2a}\right]
=
\prod_{p\leq P}m_a(p^{-\sigma}).
\]
Letting $P\to\infty$ and using Fatou's lemma gives
\[
\mathbb E\left[|\zeta(\sigma,X)|^{2a}\right]
\leq
\prod_pm_a(p^{-\sigma}),
\]
which proves \eqref{eq:real-zeta-moment-uniform}.
\end{proof}

\subsection{The weighted joint Fourier transform and Proof of Proposition \ref{prop:compact-joint-transform-rigorous}}\label{43}

We begin by establishing an Euler product representation
for the weighted joint Fourier transform.

\begin{lem}
\label{lem:H-factorization-rigorous}
Let $\sigma>1/2$, $s\in\mathbb C$ with
$\operatorname{Re}(s)>0$, and
$(\lambda,u)\in\mathbb R^2$.
Recall the weighted joint Fourier transform
$H_{\sigma,s}(\lambda,u)$ defined in
\eqref{eq:Hsigma-def-final}.
Then $H_{\sigma,s}(\lambda,u)$ admits
the following Euler-product factorization, 
\begin{equation}
\label{eq:H-factorization-rigorous}
H_{\sigma,s}(\lambda,u)
=
(2\sigma-1)^{s^2}
\prod_{p~\mathrm{prime}}
\Phi_{p,\sigma}(s;\lambda,u),
\end{equation}
where, for each prime \(p\),
\begin{align}
&\Phi_{p,\sigma}(s;\lambda,u)\nonumber\\
:=&
\int_{\mathbb T}
|1-p^{-\sigma}z|^{-2s}
\exp\left\{
-\i(2\sigma-1)\log p
\left[
\lambda\operatorname{Re}
\frac{p^{-\sigma}z}{1-p^{-\sigma}z}
+
u\operatorname{Im}
\frac{p^{-\sigma}z}{1-p^{-\sigma}z}
\right]
\right\}
\d m_{\mathbb T}(z).
\label{eq:Phi-def-final}
\end{align}
\end{lem}

\begin{proof}
Let $\zeta_{P}$  and \(L_P'\) be defined in \eqref{defineofzetaP} and \eqref{defofLP}, respectively.
By independence,
\begin{align*}
&
\mathbb E\left[
|\zeta_P(\sigma,X)|^{2s}
\e^{\i\lambda(2\sigma-1)\operatorname{Re}L'_P(\sigma,X)
+\i u(2\sigma-1)\operatorname{Im}L'_P(\sigma,X)}
\right]
\\
&\qquad=
\prod_{p\leq P}
\Phi_{p,\sigma}(s;\lambda,u).
\end{align*}
By Lemma~\ref{lem:basic-random-Euler-product}, the integrands converge
almost surely to the corresponding infinite-product integrand. Since $\lambda,u$ are real, the modulus of the finite-product
integrand is
\[
|\zeta_P(\sigma,X)|^{2\textnormal{Re}(s)}.
\]
Choose $q>1$. The proof of
Lemma~\ref{lem:basic-random-Euler-product} gives
\[
\sup_P
\mathbb E\left[|\zeta_P(\sigma,X)|^{2\textnormal{Re}(s)q}\right]<\infty.
\]
Hence the integrands are uniformly integrable. Vitali's convergence
theorem therefore allows us to pass to the limit in the expectation,
which proves \eqref{eq:H-factorization-rigorous}.
\end{proof}

We next obtain a uniform expansion of the local factors.

\begin{lem}
\label{lem:Phi-local-expansion-rigorous}
Let \(K\) be a compact subset of the open right half-plane
$
\{s\in\mathbb C:\operatorname{Re}(s)>0\}
$
and let $B_0\subset\mathbb C^2$ be bounded. 
For each prime \(p\), set
$
t_p=(2\sigma-1)\log p$.
Uniformly for
$s\in K$, $(\lambda,u)\in B_0$, and $\sigma$ sufficiently close to
$1/2$,
\begin{align}
\Phi_{p,\sigma}(s;\lambda,u)
={}&
m_s(p^{-\sigma})
-\i s\lambda t_pp^{-2\sigma}
-\frac{\lambda^2+u^2}{4}t_p^2p^{-2\sigma}
+
R_{p,\sigma},
\label{eq:Phi-local-expansion-rigorous}
\end{align}
where
\begin{equation}
\label{eq:Phi-local-remainder-rigorous}
|R_{p,\sigma}|
\ll_{K,B_0}
t_pp^{-4\sigma}+t_p^2p^{-3\sigma}+t_p^3p^{-3\sigma}.
\end{equation}
Here, $\Phi_{p,\sigma}(s;\lambda,u)$ is defined in (\ref{eq:Phi-def-final}).
Consequently,
\begin{equation}
\label{eq:Phi-remainder-sum-rigorous}
\sum_p
\sup_{\substack{s\in K\\(\lambda,u)\in B_0}}
|R_{p,\sigma}|
=
O_{K,B_0}((2\sigma-1)).
\end{equation}
\end{lem}

\begin{proof}
Write
\[
t=t_p,\qquad
W_r(z)=\frac{rz}{1-rz},
\]
and
\[
Q_{\lambda,u}(z)
=
\lambda\operatorname{Re}W_r(z)
+
u\operatorname{Im}W_r(z).
\]
So, we have 
\[
t_{p}Q_{\lambda,u}(z)=O_{B_0}((2\sigma-1))
\]
uniformly in primes $p$ and $z$ with $|z|=1$. We note that
\[
W_r(\overline z)=\overline{W_r(z)},
\]
we have
\[
\operatorname{Im}W_r(\overline z)
=
-\operatorname{Im}W_r(z).
\]
By the change of variables $z\mapsto\bar z$ and the invariance of Haar
measure on $\mathbb T$ under complex conjugation,
\[
\int_{\mathbb T}
|1-rz|^{-2s}
\operatorname{Im}W_r(z)\,\d m_{\mathbb T}(z)=0.
\]
Differentiating $m_s(r)$ gives
\[
m_s'(r)
=
\frac{2s}{r}
\int_{\mathbb T}
|1-rz|^{-2s}
\operatorname{Re}W_r(z)\,\d m_{\mathbb T}(z).
\]
Since
$
K\subset\{\operatorname{Re}(s)>0\},
$
we have
$
\inf_{s\in K}|s|>0,
$
so division by $s$ is uniform on $K$. From
\eqref{eq:ms-expansion-final},
\[
m_s'(r)=2s^2r+O_K(r^3),
\]
and hence
\begin{equation}
\label{eq:Q-first-rigorous}
\int_{\mathbb T}
|1-rz|^{-2s}
Q_{\lambda,u}(z)\,\d m_{\mathbb T}(z)
=
s\lambda r^2+O_{K,B_0}(r^4).
\end{equation}
Moreover,
\[
W_r(z)=rz+O(r^2),
\qquad
|1-rz|^{-2s}=1+O_K(r),
\]
so Haar orthogonality gives
\begin{equation}
\label{eq:Q-second-rigorous}
\int_{\mathbb T}
|1-rz|^{-2s}
Q_{\lambda,u}(z)^2\,\d m_{\mathbb T}(z)
=
\frac{\lambda^2+u^2}{2}r^2
+
O_{K,B_0}(r^3).
\end{equation}
Taylor's formula yields
\[
\e^{-\i tQ}
=
1-\i t Q-\frac{t^2Q^2}{2}
+
O_{B_0}(t^3r^3).
\]
Substitution of
\eqref{eq:Q-first-rigorous} and
\eqref{eq:Q-second-rigorous} gives
\eqref{eq:Phi-local-expansion-rigorous}.
Finally,
\begin{align*}
\sum_p t_pp^{-4\sigma}
&\ll
(2\sigma-1)\sum_{n\geq2}\frac{\log n}{n^2}
\ll(2\sigma-1),\\
\sum_p t_p^2p^{-3\sigma}
&\ll
(2\sigma-1)^2\sum_{n\geq2}\frac{(\log n)^2}{n^{3/2}}
\ll(2\sigma-1)^2,\\
\sum_p t_p^3p^{-3\sigma}
&\ll
(2\sigma-1)^3\sum_{n\geq2}\frac{(\log n)^3}{n^{3/2}}
\ll (2\sigma-1)^3.
\end{align*}
This proves \eqref{eq:Phi-remainder-sum-rigorous}.
\end{proof}

Combining the local expansion in
Lemma~\ref{lem:Phi-local-expansion-rigorous}
with the prime-sum estimates in
Lemma~\ref{lem:prime-sums-final},
we now prove
Proposition~\ref{prop:compact-joint-transform-rigorous},
which establishes the asymptotic behaviour
of the weighted joint Fourier transform
as $\sigma\downarrow1/2$.

\begin{proof}[Proof of Proposition \ref{prop:compact-joint-transform-rigorous}]
Fix a compact subset
\[
K\subset\{\operatorname{Re}(s)>0\}
\]
and a bounded closed polydisc \(B_0\subset\mathbb C^2\).
By \eqref{eq:ms-expansion-final}, there exists \(P_0=P_0(K)\)
such that, for all \(p>P_0\), \(s\in K\), and \(\sigma\geq1/2\),
\[
|m_s(p^{-\sigma})|\geq\frac34.
\]
Here $m_{s}(r)$ is defined by (\ref{defofms}).
Thus we divide only by the large-prime factors and write
\begin{align}
\prod_p\Phi_{p,\sigma}
={}&
\left(
\prod_{p\leq P_0}\Phi_{p,\sigma}
\right)
\left(
\prod_{p>P_0}m_s(p^{-\sigma})
\right)
\prod_{p>P_0}
\frac{\Phi_{p,\sigma}}
{m_s(p^{-\sigma})},
\label{eq:three-products-rigorous}
\end{align}
where $\Phi_{p,\sigma}(s;\lambda,u)$ is defined by (\ref{eq:Phi-def-final}). For the finitely many primes \(p\leq P_0\), as
\(\sigma\downarrow1/2\),
\bea\label{singleconvergence}
\Phi_{p,\sigma}(s;\lambda,u)
\longrightarrow
m_s(p^{-1/2})
\eea
uniformly for \(s\in K\) and \((\lambda,u)\in B_0\). For the unperturbed large-prime factors,
\eqref{eq:ms-expansion-final} gives
\[
(1-p^{-2\sigma})^{s^2}m_s(p^{-\sigma})
=
1+O_K(p^{-4\sigma})
\]
uniformly for \(s\in K\) and \(\sigma\geq1/2\).
Since
$
p^{-4\sigma}\leq p^{-2},
$
the corresponding infinite product converges locally uniformly, and as $\sigma\downarrow \frac{1}{2}$,
\[
\prod_{p>P_0}
(1-p^{-2\sigma})^{s^2}m_s(p^{-\sigma})
\longrightarrow
\prod_{p>P_0}
(1-p^{-1})^{s^2}m_s(p^{-1/2})
\]
uniformly for \(s\in K\). Moreover,
\begin{align}
(2\sigma-1)^{s^2}
\prod_{p>P_0}m_s(p^{-\sigma})
={}&
\bigl((2\sigma-1)\zeta(2\sigma)\bigr)^{s^2}
\prod_{p\leq P_0}(1-p^{-2\sigma})
^{s^2}
\prod_{p>P_0}
(1-p^{-2\sigma})^{s^2}m_s(p^{-\sigma}).
\label{eq:tail-unperturbed-rigorous}
\end{align}
It follows that, as $\sigma\downarrow \frac{1}{2}$,
\begin{align}
(2\sigma-1)^{s^2}
\prod_{p>P_0}m_s(p^{-\sigma})
\longrightarrow
\prod_{p\leq P_0}(1-p^{-1})
^{s^2}
\prod_{p>P_0}
(1-p^{-1})^{s^2}m_s(p^{-1/2})
\label{eq:tail-unperturbed-limit-rigorous}
\end{align}
uniformly for \(s\in K\). Observe that multiplying the limit in \eqref{eq:tail-unperturbed-limit-rigorous} by the finite-prime factor
\[
\prod_{p\leq P_0}m_s(p^{-1/2})
\]
and using \eqref{eq:ms-power-series-final} gives exactly the arithmetic factor \(a_s\) defined in \eqref{definitionofthearithmeticfactor}.
Combining this with \eqref{eq:tail-unperturbed-limit-rigorous} and \eqref{singleconvergence},  we have
\bea\label{obtainas}
(2\sigma-1)^{s^2}\left(
\prod_{p\leq P_0}\Phi_{p,\sigma}
\right)
\left(
\prod_{p>P_0}m_s(p^{-\sigma})
\right)
\longrightarrow a_{s}, \qquad \sigma\downarrow \frac{1}{2}
\eea
uniformly on \(K\times B_0\). For \(p>P_0\), set
\bea\label{defofD}
D_{p,\sigma}=D_{p,\sigma}(s;\lambda,\mu)
:=
\frac{\Phi_{p,\sigma}}
{m_s(p^{-\sigma})}-1.
\eea
Lemma~\ref{lem:Phi-local-expansion-rigorous}, together with
\[
m_s(p^{-\sigma})^{-1}
=
1+O_K(p^{-2\sigma}),
\]
gives
\begin{align}
D_{p,\sigma}
={}&
-\i s\lambda(2\sigma-1)
\frac{\log p}{p^{1+\delta}}
-
\frac{\lambda^2+u^2}{4}
(2\sigma-1)^2
\frac{(\log p)^2}{p^{1+\delta}}
+
E_{p,\sigma},
\label{eq:D-expansion-rigorous}
\end{align}
where
\begin{equation}
\label{eq:E-sum-rigorous}
\sum_{p>P_0}
\sup_{\substack{s\in K\\(\lambda,u)\in B_0}}
|E_{p,\sigma}|
=o(1),
\qquad
\sigma\downarrow1/2.
\end{equation}
The same estimates also give
\bea\label{addestimat1}
\sup_{\substack{
p>P_0,\ s\in K\\
(\lambda,u)\in B_0}}
|D_{p,\sigma}|
=o(1), \qquad \sigma\downarrow1/2,
\eea
and
\bea\label{addestimat2}
\sum_{p>P_0}
\sup_{\substack{s\in K\\(\lambda,u)\in B_0}}
|D_{p,\sigma}|
=
O_{K,B_0}(1)
\eea
by Lemma~\ref{lem:prime-sums-final}. 
Hence the branch of \(\log(1+D_{p,\sigma})\) near \(1\) is
well defined uniformly on \(K\times B_0\). Furthermore,
\[
\sum_{p>P_0}
\sup_{\substack{s\in K\\(\lambda,u)\in B_0}}
|D_{p,\sigma}|^2
\leq
\left(
\sup_{\substack{
p>P_0,\ s\in K\\
(\lambda,u)\in B_0}}
|D_{p,\sigma}|
\right)
\sum_{p>P_0}
\sup_{\substack{s\in K\\(\lambda,u)\in B_0}}
|D_{p,\sigma}|
=o(1).
\]
Therefore
\bea\label{estimationonlog}
\sum_{p>P_0}
\log(1+D_{p,\sigma})
=
\sum_{p>P_0}D_{p,\sigma}+o(1)
\eea
uniformly on \(K\times B_0\). Using \eqref{defofD}, \eqref{estimationonlog}, \eqref{eq:D-expansion-rigorous},
and Lemma~\ref{lem:prime-sums-final}, we obtain we obtain
\bea\label{productoflargeprimes}
\sum_{p>P_0}
\log
\frac{\Phi_{p,\sigma}}
{m_s(p^{-\sigma})}
\longrightarrow
-\i s\lambda-\frac{\lambda^2+u^2}{4}, \qquad \sigma\downarrow \frac{1}{2},
\eea
uniformly on \(K\times B_0\). Combining this with
\eqref{eq:three-products-rigorous} and
\eqref{obtainas} yields
\[
(2\sigma-1)^{s^2}\prod_p\Phi_{p,\sigma}(s;\lambda,u)
\longrightarrow
a_s
\e^{-\i s\lambda-(\lambda^2+u^2)/4}
\]
uniformly on \(K\times B_0\).

For complex \((\lambda,u)\in\widetilde B_0\), where
\(\widetilde B_0\subset\mathbb C^2\) is a bounded open polydisc, the
local factors \(\Phi_{p,\sigma}(s;\lambda,u)\) are entire in
\((\lambda,u)\). Moreover, the estimates
\eqref{eq:D-expansion-rigorous}--\eqref{addestimat2} and
\eqref{estimationonlog} remain valid locally uniformly on
\(K\times\widetilde B_0\). Hence the same Euler-product argument
yields a holomorphic continuation to \(\widetilde B_0\), with the
same convergence locally uniformly on
\(K\times\widetilde B_0\). On the real plane, this continuation
agrees with \(H_{\sigma,s}\) by
Lemma~\ref{lem:H-factorization-rigorous}. This proves
\eqref{eq:compact-joint-transform-rigorous}, together with the
asserted holomorphic continuation and locally uniform convergence.
\end{proof}

Having established the locally uniform convergence
of the weighted joint Fourier transform in
Proposition~\ref{prop:compact-joint-transform-rigorous},
we next derive uniform Fourier estimates that will
be used to justify the convergence of the
Gaussian-regularized moments and the subsequent
limiting arguments involving the regularized
Fourier--Riesz representation. We begin with a global decay estimate on the real
Fourier plane.

\subsection{Uniform Fourier estimates}
\subsubsection{Global decay on the real Fourier axis and proof of Proposition \ref{prop:global-Fourier-decay-rigorous}}\label{441}

Write
\[
s=a+\i b.
\]
For each prime $p$, define the probability measure
\[
\d \mu_{p,a}(z)
=
\frac{|1-p^{-\sigma}z|^{-2a}}
{m_a(p^{-\sigma})}
\,\d m_{\mathbb T}(z),
\]
where \(m_s(r)\) is defined in \eqref{defofms}.
Let $\Phi_{p,\sigma}(s;\lambda,u)$ be defined in (\ref{eq:Phi-def-final}).
Then
\[
\Phi_{p,\sigma}(s;\lambda,u)
=
m_a(p^{-\sigma})\chi_{p,\sigma},
\]
where
\begin{equation}\label{defofchi}
\chi_{p,\sigma}
=
\int_{\mathbb T}
\e^{\i\Psi_{p,\sigma}(z)}
\,\d \mu_{p,a}(z)
\end{equation}
and
\begin{equation}\label{defofPsi}
\Psi_{p,\sigma}(z)
=
-2b\log|1-p^{-\sigma}z|
-
(2\sigma-1)\log p\left(
\lambda\operatorname{Re}\frac{p^{-\sigma}z}{1-p^{-\sigma}z}
+
u\operatorname{Im}\frac{p^{-\sigma}z}{1-p^{-\sigma}z}
\right).
\end{equation}
For real $\lambda,u$, the phase is real, and therefore
$
|\chi_{p,\sigma}|\leq1.
$
\begin{lem}
\label{lem:local-damping-rigorous}
Let \(s=a+ib\in\mathbb C\) with
\(\operatorname{Re}s=a>0\), and suppose that
\[
0<a_0\leq a\leq a_1,
\qquad
|b|\leq B.
\]
Let \(\sigma>1/2\) and let \(p\) be a prime. Let
\(\chi_{p,\sigma}\) be defined as in \eqref{defofchi}, and set
\[
t_p:=(2\sigma-1)\log p.
\]
Then there exist constants
\[
r_0=r_0(a_0,a_1,B)>0,
\qquad
\eta=\eta(a_0,a_1,B)>0,
\]
such that, whenever
\[
0<p^{-\sigma}\leq r_0
\]
and
\[
p^{-\sigma}
\sqrt{
(2b-t_p\lambda)^2+t_p^2\mu^2
}
\leq\eta,
\]
one has
\begin{align}
\log|\chi_{p,\sigma}(\lambda,\mu)|
={}&
-\frac{p^{-2\sigma}}{4}
\left[
(2b-t_p\lambda)^2+t_p^2\mu^2
\right]
+
E_{p,\sigma}(\lambda,\mu),
\label{eq:local-damping-rigorous}
\end{align}
where
\begin{align}
|E_{p,\sigma}(\lambda,\mu)|
\ll_{a_0,a_1,B}
{}&
B^2p^{-3\sigma}
+t_p^2(\lambda^2+\mu^2)p^{-3\sigma}
+
B^4p^{-4\sigma}
+t_p^4(\lambda^2+\mu^2)^2p^{-4\sigma}.
\label{eq:local-damping-error-rigorous}
\end{align}
The implied constant is uniform in
\(\sigma,p,\lambda,\mu\), and \(s=a+\i b\) subject to the above
conditions.
\end{lem}

\begin{proof}
Note that 
uniformly for \(z\in\mathbb T\), as \(p^{-\sigma}\to0\),
\[
-\log|1-p^{-\sigma}z|
=
p^{-\sigma}\operatorname{Re}z
+
O(p^{-2\sigma})
\]
and
\[
\frac{p^{-\sigma}z}{1-p^{-\sigma}z}
=
p^{-\sigma}z
+
O(p^{-2\sigma}).
\]
Hence, by the definition of \(\Psi_{p,\sigma}\) in (\ref{defofPsi}),
\begin{align}\label{expansionforPsi}
\Psi_{p,\sigma}(z)
={}&
p^{-\sigma}
\left[
(2b-t_p\lambda)\operatorname{Re}z
-
t_p\mu\operatorname{Im}z
\right]+
O\left(
p^{-2\sigma}
\left(
B+t_p\sqrt{\lambda^2+\mu^2}
\right)
\right),
\end{align}
uniformly for \(z\in\mathbb T\). Moreover, by the triangle inequality,
\[
t_p\sqrt{\lambda^2+\mu^2}
\leq
\sqrt{
(2b-t_p\lambda)^2+t_p^2\mu^2
}
+
2|b|,
\]
and therefore
\[
p^{-\sigma}t_p\sqrt{\lambda^2+\mu^2}
\leq
\eta+2Bp^{-\sigma}.
\]
It follows that
\[
\sup_{z\in\mathbb T}
|\Psi_{p,\sigma}(z)|
\ll_{a_0,a_1,B}
\eta+Bp^{-\sigma}.
\]
Thus, after choosing \(\eta>0\) and \(r_0>0\) sufficiently small,
\begin{equation}\label{farfrom0}
|\chi_{p,\sigma}(\lambda,\mu)-1|<\frac12.
\end{equation}
In particular,
\(\chi_{p,\sigma}(\lambda,\mu)\) stays uniformly away from zero. Next, from
\[
\frac{\d\mu_{p,a}}{\d m_{\mathbb T}}(z)
=
\frac{|1-p^{-\sigma}z|^{-2a}}
{m_a(p^{-\sigma})},
\]
together with \(a_0\leq a\leq a_1\), we have uniformly for
\(z\in\mathbb T\),
\begin{equation}\label{relationbetweenmuandlesbergue}
\frac{\d \mu_{p,a}}{\d m_{\mathbb T}}(z)
=
1+2ap^{-\sigma}\operatorname{Re}z
+O_{a_0,a_1}(p^{-2\sigma}).
\end{equation}
Using the expansion of \(\Psi_{p,\sigma}\) in
\eqref{expansionforPsi} and (\ref{relationbetweenmuandlesbergue}), a direct calculation gives
\begin{equation}
\label{eq:phase-variance-rigorous}
\begin{aligned}
\operatorname{Var}_{\mu_{p,a}}(\Psi_{p,\sigma}):=
&\int_{\mathbb T}\Psi_{p,\sigma}(z)^2\,\d\mu_{p,a}(z)
-
\left(
\int_{\mathbb T}\Psi_{p,\sigma}(z)\,\d\mu_{p,a}(z)
\right)^2\nonumber\\
=&
\frac{p^{-2\sigma}}{2}
\left[
(2b-t_p\lambda)^2+t_p^2\mu^2
\right]+
O_{a_0,a_1,B}\left(
B^2p^{-3\sigma}
+
t_p^2(\lambda^2+\mu^2)p^{-3\sigma}
\right).
\end{aligned}
\end{equation}
Let \(\Psi_{p,\sigma}'\) be an independent copy of
\(\Psi_{p,\sigma}\) under \(\mu_{p,a}\). Then
\[
|\chi_{p,\sigma}(\lambda,\mu)|^2
=
\mathbb E_{\mu_{p,a}\otimes\mu_{p,a}}
\left[
\mathrm e^{\i(\Psi_{p,\sigma}-\Psi_{p,\sigma}')}
\right]
=
\mathbb E_{\mu_{p,a}\otimes\mu_{p,a}}
\left[
\cos(\Psi_{p,\sigma}-\Psi_{p,\sigma}')
\right].
\]
Using
\[
\cos x
=
1-\frac{x^2}{2}+O(x^4)
\]
and
\[
\frac12
\mathbb E_{\mu_{p,a}\otimes\mu_{p,a}}
\left[
(\Psi_{p,\sigma}-\Psi_{p,\sigma}')^2
\right]
=
\operatorname{Var}_{\mu_{p,a}}(\Psi_{p,\sigma}),
\]
we obtain
\[
|\chi_{p,\sigma}(\lambda,\mu)|^2
=
1-
\operatorname{Var}_{\mu_{p,a}}(\Psi_{p,\sigma})
+
O\left(
\mathbb E_{\mu_{p,a}\otimes\mu_{p,a}}
\left[|\Psi_{p,\sigma}-\Psi_{p,\sigma}'|^4\right]
\right).
\]
The bounds on \(\Psi_{p,\sigma}\) in
\eqref{expansionforPsi} imply
\[
\mathbb E_{\mu_{p,a}\otimes\mu_{p,a}}
\left[|\Psi_{p,\sigma}-\Psi_{p,\sigma}'|^4\right]
\ll_{a_0,a_1,B}
B^4p^{-4\sigma}
+
t_p^4(\lambda^2+\mu^2)^2p^{-4\sigma}.
\]
Consequently,
\[
\begin{aligned}
|\chi_{p,\sigma}(\lambda,\mu)|^2
={}&
1-
\frac{p^{-2\sigma}}{2}
\left[
(2b-t_p\lambda)^2+t_p^2\mu^2
\right]\\
&+
O_{a_0,a_1,B}\left(
B^2p^{-3\sigma}
+t_p^2(\lambda^2+\mu^2)p^{-3\sigma}
+B^4p^{-4\sigma}
+t_p^4(\lambda^2+\mu^2)^2p^{-4\sigma}
\right).
\end{aligned}
\]
By (\ref{farfrom0}), \(|\chi_{p,\sigma}(\lambda,\mu)-1|<1/2\), the quantity
\(|\chi_{p,\sigma}(\lambda,\mu)|^2\) is uniformly bounded away from
zero. Hence we may take logarithms and use
\[
\log(1+x)=x+O(x^2)
\]
for \(x\) sufficiently small. Hence, we obtain
\eqref{eq:local-damping-rigorous}--\eqref{eq:local-damping-error-rigorous}.
\end{proof}

We now establish a uniform global decay estimate
for the weighted joint Fourier transform
$H_{\sigma,s}$ on $\mathbb R^2$ in Proposition \ref{prop:global-Fourier-decay-rigorous}.
\begin{proof}[Proof of Proposition \ref{prop:global-Fourier-decay-rigorous}]
Write 
\[
s=a+\i b,\qquad
\delta=2\sigma-1,\qquad
\rho=(\lambda^2+u^2)^{1/2}.
\]
Since \(K\) is a compact subset of the open right half-plane, there
exist constants \(0<a_0\leq a_1\) and \(B>0\), depending only on
\(K\), such that
\[
a_0\leq a\leq a_1,
\qquad
|b|\leq B,
\qquad s\in K.
\]
Fix \(N>0\), and choose \(A=A_{K,N}>0\) sufficiently large.
The precise lower bound required for \(A\) will be specified at the
end of the proof. Fix \(\sigma_0\in(1/2,1)\). Let \(r_0,\eta_0>0\) be as in
Lemma~\ref{lem:local-damping-rigorous}, and let \(C_E=C_E(K)\)
denote the constant occurring in its error estimate. Choose
\(0<\eta\leq\eta_0\) sufficiently small so that
\bea\label{chooseofeta}
\frac{C_E\eta^2}{16}\leq\frac1{32}.
\eea
Let \(c_*>0\) and \(P_*\geq2\) be such that
Lemma~\ref{lem:prime-sums-final} gives
\[
\sum_{p\geq P}
\frac{(\log p)^2}{p^{1+\delta}}
\geq
c_*
\frac{P^{-\delta}}{\delta^2}
\bigl(1+\delta\log P\bigr),
\qquad
P\geq P_*,
\]
uniformly for \(0<\delta\leq 2\sigma_{0}-1\).
Set
\bea\label{defofc1}
c_1:=\frac{c_*\eta^2}{16}>0.
\eea
We now choose \(Y_0>3\) sufficiently large so that
\bea\label{constra1}
\e^{-Y_0/2}\leq r_0,
\qquad
2B\e^{-Y_0/2}\leq\frac{\eta}{4},
\qquad
C_E\e^{-Y_0/2}\leq\frac1{32},
\eea
and
\bea\label{constra2}
\e^{Y_0}\geq P_*.
\eea
We also require
\bea\label{constra3}
c_1\frac{\e^y}{y^2}\geq4Ay,
\qquad y\geq Y_0.
\eea
Having fixed \(\eta\) and \(Y_0\), choose
\begin{equation}\label{defofc0}
c_0:=\frac{\eta \e^{\sigma_0Y_0}}{2Y_0}.
\end{equation}
Finally, choose \(P_0\) sufficiently large so that, for every
\(p>P_0\),
\[
p^{-1/2}\leq r_0
\]
and
\[
2Bp^{-1/2}
+
c_0\frac{\log p}{\sqrt p}
<\eta.
\]
We now distinguish the two cases
\[
\delta\rho\leq c_0
\qquad\text{and}\qquad
\delta\rho>c_0.
\]
\medskip
\noindent
\emph{Case 1: \(\delta\rho\leq c_0\).} 
By  Lemma \ref{lem:H-factorization-rigorous}, and  the independence
of \(\{X(p)\}_p\), we have
\begin{align*}
H_{\sigma,s}(\lambda,u)
&=
\delta^{s^2}
\prod_p
\int_{\mathbb T}
|1-p^{-\sigma}z|^{-2a}
\mathrm e^{\i\Psi_{p,\sigma}(z)}
\,\d m_{\mathbb T}(z)\\
&=
\delta^{s^2}
\prod_p
m_a(p^{-\sigma})\chi_{p,\sigma}(\lambda,u),
\end{align*}
where \(m_a(r)\), \(\Psi_{p,\sigma}(z)\), and \(\chi_{p,\sigma}\)
are defined in \eqref{defofms}, \eqref{defofPsi}, and
\eqref{defofchi}, respectively.
Since
$
|\delta^{s^2}|
=
\delta^{a^2-b^2},
$
combining Lemma~\ref{lem:real-zeta-moment-uniform} with the above
prime factorization of \(H_{\sigma,s}\) gives, 
\begin{equation}
\label{eq:H-modulus-preliminary}
|H_{\sigma,s}(\lambda,u)|
\ll_K
\delta^{-b^2}
\prod_p|\chi_{p,\sigma}(\lambda,u)|.
\end{equation}
If
$
\chi_{p,\sigma}(\lambda,u)=0
$
for some \(p\leq P_0\), then
\(H_{\sigma,s}(\lambda,u)=0\), and there is nothing to prove.
Otherwise, since
$
|\chi_{p,\sigma}(\lambda,u)|\leq1
$ by (\ref{defofchi}),
the finitely many primes \(p\leq P_0\) may be discarded when deriving
an upper bound. Note that for $p>P_{0}$,
\[
p^{-\sigma}
\sqrt{
(2b-\delta\lambda\log p)^2
+
\delta^2u^2(\log p)^2
}
\leq
2Bp^{-\sigma}
+
c_0(\log p)p^{-\sigma},
\]
which is smaller than the constant \(\eta\) in
Lemma~\ref{lem:local-damping-rigorous}. So applying Lemma~\ref{lem:local-damping-rigorous},
\begin{align*}
\sum_p\log|\chi_{p,\sigma}(\lambda,u)|
&\leq
-\frac14
\sum_{p>P_0}
\frac{
(2b-\delta\lambda\log p)^2
+
\delta^2u^2(\log p)^2
}{
p^{1+\delta}
}
+
\sum_{p>P_0}|E_{p,\sigma}|.
\end{align*}
Under the assumption \(\delta\rho\leq c_0\),
\eqref{eq:local-damping-error-rigorous} gives
\[
\sum_{p>P_0}|E_{p,\sigma}|
=
O_K(1).
\]
Expanding the quadratic term and applying
Lemma~\ref{lem:uniform-prime-tail-final}, we obtain
\begin{align}
\frac14
\sum_{p>P_0}
\frac{
(2b-\delta\lambda\log p)^2
+
\delta^2u^2(\log p)^2
}{
p^{1+\delta}
}
=
b^2\log\frac1\delta
-b\lambda
+\frac{\rho^2}{4}
+
O_K(1).
\label{eq:case1-quadratic}
\end{align}
Since \(|b|\leq B\),
\[
b\lambda-\frac{\rho^2}{4}
\leq
-\frac{\rho^2}{8}+2B^2.
\]
Consequently,
\[
\prod_p|\chi_{p,\sigma}(\lambda,u)|
\ll_K
\delta^{b^2}\e^{-\rho^2/8}.
\]
Substituting this into \eqref{eq:H-modulus-preliminary} gives
\[
|H_{\sigma,s}(\lambda,u)|
\ll_K
\e^{-\rho^2/8}.
\]
In particular, for every \(N>0\),
\[
|H_{\sigma,s}(\lambda,u)|
\ll_{K,N}
(1+\rho^2)^{-N}.
\]

\medskip
\noindent
\emph{Case 2: \(\delta\rho>c_0\).}
By the definition of \(c_0\) in (\ref{defofc0}),
\[
\frac{\eta}{4c_0}
=
\frac12Y_0\e^{-\sigma_0Y_0}.
\]
Hence, since \(\delta\rho>c_0\) and \(\sigma\leq\sigma_0\),
\[
\frac{\eta}{4\delta\rho}
<
\frac{\eta}{4c_0}
=
\frac12Y_0e^{-\sigma_0Y_0}
<
Y_0e^{-\sigma Y_0}.
\]
Since \(Y_0>3\), the function $
x\longmapsto x\e^{-\sigma x}$
is strictly decreasing on \([Y_0,\infty)\) and tends to zero as
\(x\to\infty\). Therefore there exists a unique \(y>Y_0\) such that
\begin{equation}
\label{eq:y-definition-rigorous}
\delta\rho\,y \e^{-\sigma y}
=
\frac{\eta}{4}.
\end{equation}
Set
\bea\label{defofP}
P:=\e^y.
\eea
We first verify that Lemma~\ref{lem:local-damping-rigorous} applies
for every prime \(p\geq P\). Since \(p\geq P=\e^y\), we have
\(\log p\geq y> Y_0\). By \eqref{eq:y-definition-rigorous},
\begin{equation}
\label{eq:case2-smallness-rho}
\delta\rho(\log p)p^{-\sigma}
\leq
\frac{\eta}{4}.
\end{equation}
Moreover, by our choice of \(Y_0\) in (\ref{constra1}),
\bea\label{comparision1}
p^{-\sigma}
\leq
\e^{-\sigma y}
<
\e^{-Y_0/2}
\leq r_0
\eea
and
\bea\label{comparision2}
2Bp^{-\sigma}
\leq
2B\e^{-Y_0/2}
\leq
\frac{\eta}{4}.
\eea
Therefore
\[
\begin{aligned}
&p^{-\sigma}
\sqrt{(2b-\delta\lambda\log p)^2
+\delta^2u^2(\log p)^2}
\\
&\qquad\leq
2Bp^{-\sigma}
+
\delta\rho(\log p)p^{-\sigma}
\leq
\frac{\eta}{2}
<\eta.
\end{aligned}
\]
Thus the local damping estimate in Lemma \ref{lem:local-damping-rigorous} applies for every \(p\geq P\). If \(\chi_{p,\sigma}(\lambda,u)=0\) for some prime \(p\), then
\(H_{\sigma,s}(\lambda,u)=0\), and there is nothing to prove.
Otherwise, since \(|\chi_{p,\sigma}|\leq1\),
\[
\sum_p\log|\chi_{p,\sigma}|
\leq
\sum_{p\geq P}\log|\chi_{p,\sigma}|.
\]
Lemma~\ref{lem:local-damping-rigorous} gives
\[
\log|\chi_{p,\sigma}(\lambda,u)|
=
-\frac{1}{4p^{1+\delta}}
\left(
(2b-\delta\lambda\log p)^2
+
\delta^2u^2(\log p)^2
\right)
+
E_{p,\sigma},
\]
where
\[
|E_{p,\sigma}|
\leq
C_E\left(
B^2p^{-3\sigma}
+
\delta^2\rho^2(\log p)^2p^{-3\sigma}
+
B^4p^{-4\sigma}
+
\delta^4\rho^4(\log p)^4p^{-4\sigma}
\right).
\]
Put
\[
Q_\sigma(\rho)
:=
\rho^2\delta^2
\sum_{p\geq P}
\frac{(\log p)^2}{p^{1+\delta}}.
\]
Note that \(1+\delta=2\sigma\). By (\ref{comparision1}),
\[
\begin{aligned}
\sum_{p\geq P}
\delta^2\rho^2(\log p)^2p^{-3\sigma}
&\leq
P^{-\sigma}Q_\sigma(\rho)<
\e^{-Y_0/2}Q_\sigma(\rho).
\end{aligned}
\]
Furthermore, by \eqref{eq:case2-smallness-rho},
\[
\begin{aligned}
\sum_{p\geq P}
\delta^4\rho^4(\log p)^4p^{-4\sigma}
&\leq 
\frac{\eta^2}{16}Q_\sigma(\rho).
\end{aligned}
\]
The remaining terms satisfy
\[
\sum_{p\geq P}
\left(
B^2p^{-3\sigma}+B^4p^{-4\sigma}
\right)
\ll_K1.
\]
By the choices of \(Y_0\) and \(\eta\) in \eqref{constra1} and \eqref{chooseofeta}, respectively,
\[
C_E\e^{-Y_0/2}\leq\frac1{32},
\qquad
\frac{C_E\eta^2}{16}\leq\frac1{32}.
\]
It follows that
\begin{equation}
\label{eq:error-absorption-rigorous}
\sum_{p\geq P}|E_{p,\sigma}|
\leq
\frac1{16}Q_\sigma(\rho)+C_K.
\end{equation}
Next,
\[
(2b-\delta\lambda\log p)^2
+
\delta^2u^2(\log p)^2
\geq
\frac12\delta^2\rho^2(\log p)^2-4B^2.
\]
Therefore
\[
\begin{aligned}
\sum_p\log|\chi_{p,\sigma}|
&\leq
-\frac18Q_\sigma(\rho)
+
B^2\sum_{p\geq P}\frac1{p^{1+\delta}}
+
\sum_{p\geq P}|E_{p,\sigma}|
\\
&\leq
-\frac1{16}Q_\sigma(\rho)
+
B^2\sum_{p\geq P}\frac1{p^{1+\delta}}
+
C_K.
\end{aligned}
\]
By (\ref{eq:prime-sum-zero-final}),
we obtain
\begin{equation}
\label{eq:case2-logchi}
\sum_p\log|\chi_{p,\sigma}|
\leq
-\frac1{16}Q_\sigma(\rho)
+
C_K\log\frac1\delta
+
C_K.
\end{equation}
Here, throughout the proof, \(C_K>0\) denotes a constant depending only on
\(K\), whose value may change from line to line.
By \eqref{defofP}, \eqref{constra2}, and the fact that \(y\geq Y_0\), we have
\[
P=\e^y\geq \e^{Y_0}\geq P_*.
\]
Hence Lemma~\ref{lem:uniform-prime-tail-final} gives
\[
\sum_{p\geq P}
\frac{(\log p)^2}{p^{1+\delta}}
\geq
c_*
\frac{P^{-\delta}}{\delta^2}
(1+\delta\log P).
\]
Combining the above estimate with \eqref{eq:y-definition-rigorous},  we obtain
\begin{equation}
\label{eq:Q-case2-lower}
Q_\sigma(\rho)
\geq
c_1
\frac{\e^y(1+\delta y)}{\delta^2y^2}
\end{equation}
with $c_{1}$ as in (\ref{defofc1}). By our initial choice of \(Y_0\) as in (\ref{constra3}),
\[
c_1\frac{\e^y}{y^2}\geq4Ay,
\qquad y\geq Y_0>1.
\]
Hence
\[
Q_\sigma(\rho)
\geq
\frac{4Ay}{\delta^2}.
\]
Since \(0<\delta<1\), we have
\[
\delta^{-2}\geq1
\qquad\text{and}\qquad
\delta^{-2}\geq\log\frac1\delta.
\]
It follows that
\[
Q_\sigma(\rho)
\geq
2A\left(
y+\log\frac1\delta
\right).
\]
Moreover, from \eqref{eq:y-definition-rigorous},
\[
\log(2+\rho)
\leq
y+\log\frac1\delta+C_K.
\]
Consequently,
\begin{equation}
\label{eq:Q-coercivity-final}
Q_\sigma(\rho)
\geq
A\left(
\log\frac1\delta+\log(2+\rho)
\right)
-C_{K,A}.
\end{equation}
Substituting \eqref{eq:Q-coercivity-final} into
\eqref{eq:case2-logchi}, we obtain
\[
\sum_p\log|\chi_{p,\sigma}|
\leq
-\frac{A}{16}\log(2+\rho)
-
\left(\frac{A}{16}-C_K\right)
\log\frac1\delta
+
\widetilde{C}_{K,A}.
\]
Therefore
\[
\prod_p|\chi_{p,\sigma}|
\ll_{K,A}
\delta^{A/16-C_K}(2+\rho)^{-A/16}.
\]
Combining this with \eqref{eq:H-modulus-preliminary} and using
\(|b|\leq B\), we find
\[
|H_{\sigma,s}(\lambda,u)|
\ll_{K,A}
\delta^{A/16-C_K-B^2}
(2+\rho)^{-A/16}.
\]
We now specify the lower bound required for \(A\). Choose
\(A=A_{K,N}\) sufficiently large so that
\[
\frac{A}{16}\geq C_K+B^2
\qquad\text{and}\qquad
\frac{A}{16}\geq2N.
\]
Since \(0<\delta<1\),
\[
|H_{\sigma,s}(\lambda,u)|
\ll_{K,N}
(2+\rho)^{-2N}
\ll_{K,N}
(1+\rho^2)^{-N}.
\]
This completes Case~2, and hence the proof of the proposition.
\end{proof}

In addition to the global decay estimate established
in Proposition~\ref{prop:global-Fourier-decay-rigorous},
we require local bounds for the derivatives of
$H_{\sigma,s}$ to control the Taylor remainders
arising in the subsequent regularization of
Riesz distributions.

\subsubsection{Local Fourier derivative bounds}
\begin{lem}
\label{lem:local-Fourier-derivatives-rigorous}
Let \(K\) be a compact subset of the open right half-plane
$
\{s\in\mathbb C:\operatorname{Re}(s)>0\}
$
and let $R>0$. For every pair of nonnegative integers
$\alpha,\beta$, there exist $\sigma_0>1/2$ and
$C_{K,R,\alpha,\beta}>0$ such that
\begin{equation}
\label{eq:local-Fourier-derivatives-bound-rigorous}
\sup_{\substack{
1/2<\sigma\leq\sigma_0\\
s\in K\\
|\lambda|\leq R,\ |u|\leq R}}
\left|
\partial_\lambda^\alpha
\partial_u^\beta
H_{\sigma,s}(\lambda,u)
\right|
\leq
C_{K,R,\alpha,\beta}.
\end{equation}
Here \(H_{\sigma,s}(\lambda,u)\) is defined by
\eqref{eq:Hsigma-def-final}, and \(\partial_\lambda\) and
\(\partial_u\) denote partial differentiation with respect to
\(\lambda\) and \(u\), respectively. Moreover, as $\sigma\downarrow 1/2$,
\begin{align}
\partial_\lambda^\alpha
\partial_u^\beta
H_{\sigma,s}(\lambda,u)
\longrightarrow
\partial_\lambda^\alpha
\partial_u^\beta
\left[
a_s \e^{-\i s\lambda-(\lambda^2+u^2)/4}
\right]
\label{eq:local-Fourier-derivatives-limit-rigorous}
\end{align}
locally uniformly for bounded real $(\lambda,u)$ and locally uniformly on
$\operatorname{Re}(s)>0$. Here $a_{s}$ is defined in (\ref{definitionofthearithmeticfactor}). 
\end{lem}

\begin{proof}
Choose
\[
R<R_2<R_1.
\]
By Proposition~\ref{prop:compact-joint-transform-rigorous}, for
\(\sigma>1/2\) sufficiently close to \(1/2\), the function
\(H_{\sigma,s}\) admits a holomorphic continuation
$
\mathcal H_{\sigma,s}(\lambda,u)
$
to the polydisc
$
|\lambda|, |u|<R_1$,
and
\[
\mathcal H_{\sigma,s}(\lambda,u)
\longrightarrow
a_s\e^{-\i s\lambda-(\lambda^2+u^2)/4}
\]
uniformly for \(s\in K\) and
\[
|\lambda|\leq R_2,\qquad |u|\leq R_2.
\]
Hence we may choose \(\sigma_0>1/2\) such that the family
\(\mathcal H_{\sigma,s}\) is uniformly bounded on this closed
polydisc for
\[
1/2<\sigma\leq\sigma_0,\qquad s\in K.
\]
Set
\[
r:=\frac{R_2-R}{2}>0.
\]
For \(|\lambda|,|u|\leq R\), the circles
\[
|z-\lambda|=r,\qquad |w-u|=r
\]
are contained in the polydisc
\[
|z|<R_2,\qquad |w|<R_2.
\]
The two-variable Cauchy integral formula therefore gives
\[
\left|
\partial_\lambda^\alpha\partial_u^\beta
\mathcal H_{\sigma,s}(\lambda,u)
\right|
\leq
\frac{\alpha!\beta!}{r^{\alpha+\beta}}
\sup_{\substack{|z|\leq R_2\\|w|\leq R_2}}
|\mathcal H_{\sigma,s}(\zeta,\omega)|.
\]
Since \(\mathcal H_{\sigma,s}=H_{\sigma,s}\) for real
\((\lambda,u)\), this proves
\eqref{eq:local-Fourier-derivatives-bound-rigorous}. Applying the same Cauchy estimate to
\[
\mathcal H_{\sigma,s}(\lambda,u)
-
a_s\e^{-\i s\lambda-(\lambda^2+u^2)/4}
\]
and using its uniform convergence to zero on
\(K\times\{|\lambda|,|u|\leq R_2\}\), we obtain
\[
\sup_{\substack{
s\in K\\
|\lambda|\leq R,\ |u|\leq R}}
\left|
\partial_\lambda^\alpha\partial_u^\beta
H_{\sigma,s}(\lambda,u)
-
\partial_\lambda^\alpha\partial_u^\beta
\left[
a_s\e^{-\i s\lambda-(\lambda^2+u^2)/4}
\right]
\right|
\longrightarrow0
\]
as \(\sigma\downarrow1/2\). This proves
\eqref{eq:local-Fourier-derivatives-limit-rigorous}.
\end{proof}

In view of the global Fourier decay and local
derivative bounds established above for $H_{\sigma,s}$,
we now turn to the Fourier representation of
Gaussian-regularized complex moments.
We first derive a Fourier inversion formula for
$|x|^{2s}\e^{-t|x|^2}$, which will be used to
express the Gaussian-regularized moments in
terms of $H_{\sigma,s}$.

\subsection{Fourier analysis for the Gaussian-regularized moment}
\begin{lem}
\label{lem:fourier-inversion-power-gaussian}
Let $s\in\mathbb C$ with
$
\operatorname{Re}(s)>0.
$
Then, for every $(x,y)\in\mathbb R^2$ and $t\in \mathbb{R}_{> 0}$,
\begin{align}
(x^2+y^2)^s \e^{-t(x^2+y^2)}
=
\frac{\Gamma(s+1)}{4\pi t^{s+1}}
\int_{\mathbb R^2}
{}_1F_1\left(
s+1;1;-\frac{\lambda^2+u^2}{4t}
\right)
\e^{\i(\lambda x+uy)}
\,\d\lambda\,\d u.
\label{eq:fourier-inversion-power-gaussian}
\end{align}
The integral on the right-hand side is absolutely convergent. Let $T_{\sigma}$ be as in (\ref{defofTsigma}), that is 
\begin{equation*}
T_{\sigma}=(2\sigma-1)L'(\sigma,X).
\end{equation*}
Consequently, taking
\[
x=\operatorname{Re}(T_{\sigma}),
\qquad
y=\operatorname{Im}\left(T_{\sigma}\right),
\]
we obtain, pathwise,
\begin{align}
|T_\sigma|^{2s}\e^{-t|T_\sigma|^2}
=
\frac{\Gamma(s+1)}{4\pi t^{s+1}}
\int_{\mathbb R^2}
{}_1F_1\left(
s+1;1;-\frac{\lambda^2+u^2}{4t}
\right)
\e^{\i\lambda\operatorname{Re}T_\sigma
+\i u\operatorname{Im}T_\sigma}
\,\d\lambda\,\d u.
\label{eq:fourier-inversion-Tsigma}
\end{align}
Here powers of positive real numbers are defined using the real
logarithm.
\end{lem}

\begin{proof}
Let
\[
f_{s,t}(x,y)
:=(x^2+y^2)^s\e^{-t(x^2+y^2)}.
\]
Note that $f_{s,t}$ is radial. Using polar coordinate, its Fourier transform
is given by
\[
\widehat f_{s,t}(\lambda,u)
=
2\pi\int_0^\infty
r^{2s+1}\e^{-tr^2}J_0(\rho r)\,\d r,
\qquad
\rho=\sqrt{\lambda^2+u^2},
\]
where $J_0$ is the Bessel function of the first kind of order
zero.
By the standard Hankel transform formula
(see \cite[Eq.~(13.4.3)]{NIST}),
\bea\label{eq:fst-Fourier}
\widehat {f_{s,t}}(\lambda,u)
=
\frac{\pi\Gamma(s+1)}{t^{s+1}}
{}_1F_1\left(
s+1;1;-\frac{\rho^2}{4t}
\right).
\eea
It remains to verify that
\[
\widehat {f_{s,t}}\in L^1(\mathbb R^2),
\]
so that the ordinary Fourier inversion theorem is applicable. For bounded $\rho$, the elementary estimate
\begin{equation}\label{basicboundforFourier}
|\widehat {f_{s,t}}(\lambda,u)|
\leq
\|f_{s,t}\|_{L^1}
\end{equation}
is sufficient. For large $\rho$, we use the standard Kummer asymptotic on the
negative real axis. By the large-argument expansion for Kummer's function
\cite[Eq.~(13.7.2)]{NIST},
with $a=s+1$, $b=1$, and $z=-x$, we obtain, as
$x\to+\infty$,
\begin{equation}
\label{eq:Kummer-negative-axis}
{}_1F_1(s+1;1;-x)
=
\frac{x^{-s-1}}{\Gamma(-s)}
\left(1+O_s(x^{-1})\right)
+
O_s\left(
\e^{-x/2}(1+x)^{C_s}
\right),
\end{equation}
with the algebraic term absent when $s\in\mathbb Z_{\geq0}$,
since then $1/\Gamma(-s)=0$. In the latter case, Kummer's
transformation gives
\[
{}_1F_1(s+1;1;-x)
=
\e^{-x}{}_1F_1(-s;1;x),
\]
and the second factor is a polynomial. It follows from \eqref{eq:fst-Fourier} that, for $\rho\geq1$,
\begin{equation}
\label{eq:fst-decay}
|\widehat {f_{s,t}}(\lambda,u)|
\leq
C_{s,t}
\left[
\rho^{-2\operatorname{Re}(s)-2}
+
\e ^{-c_t\rho^2}(1+\rho)^{C_s}
\right].
\end{equation}
Therefore, using polar coordinates,
\begin{align*}
\int_{\rho\geq1}
|\widehat {f_{s,t}}(\lambda,u)|
\,\d\lambda\,\d u
&\ll_{s,t}
\int_1^\infty
\rho^{-2\operatorname{Re}(s)-2}
\rho\,\d\rho+
\int_1^\infty
\e^{-c_t\rho^2}(1+\rho)^{C_s}
\rho\,\d\rho.
\end{align*}
The second integral is finite, and the first is finite since $\textnormal{Re}(s)>0$. Hence
\begin{equation}\label{conclusion}
\widehat {f_{s,t}}\in L^1(\mathbb R^2).  
\end{equation}
We have therefore proved that
$
f_{s,t},\widehat {f_{s,t}}\in L^1(\mathbb R^2),
$
and $f_{s,t}$ is continuous.  Hence, the Fourier inversion theorem gives
\eqref{eq:fourier-inversion-power-gaussian}.
Finally, substituting
\[
x=\operatorname{Re}T_\sigma,
\qquad
y=\operatorname{Im}T_\sigma
\]
yields \eqref{eq:fourier-inversion-Tsigma}
pathwise.
\end{proof}

Applying the Fourier inversion formula established
in Lemma~\ref{lem:fourier-inversion-power-gaussian}
and taking expectations, we obtain the following
Fourier representation of the Gaussian-regularized
complex moments in terms of $H_{\sigma,s}$.

\begin{lem}[Fourier representation of the Gaussian-regularized moment]
\label{lem:smoothed-moment-fourier}
Let  $s\in \mathbb{C}$ with $\textnormal{Re}(s)>0$.  Let $T_{\sigma}$ be as in (\ref{defofTsigma}). 
Then for any $\sigma$ with $\sigma>1/2$ and any $t$ with $t>0$,
\begin{align}
&
(2\sigma-1)^{s^2}
\mathbb E\left[
|\zeta(\sigma,X)|^{2s}
|T_\sigma|^{2s}\e^{-t|T_\sigma|^2}
\right]
\nonumber\\
&\qquad=
\frac{\Gamma(s+1)}{4\pi t^{s+1}}
\int_{\mathbb R^2}
{}_1F_1\left(
s+1;1;-\frac{\lambda^2+u^2}{4t}
\right)
H_{\sigma,s}(\lambda,u)
\,\d\lambda\,\d u.
\label{eq:smoothed-moment-fourier}
\end{align}
Moreover, the integral on the right-hand side is absolutely
convergent.
\end{lem}

\begin{proof}
We first multiply \eqref{eq:fourier-inversion-Tsigma} by
$
(2\sigma-1)^{s^2}|\zeta(\sigma,X)|^{2s}
$
and take expectations. We now verify that Fubini's theorem is
applicable. This is because that
\begin{align}\label{eq:Fubini-smoothed-moment}
&\int_{\mathbb R^2}
\mathbb E\left[
\left|
(2\sigma-1)^{s^2}
|\zeta(\sigma,X)|^{2s}
{}_1F_1\left(
s+1;1;-\frac{\lambda^2+u^2}{4t}
\right)
\e^{\i\lambda\operatorname{Re}T_\sigma
+\i u\operatorname{Im}T_\sigma}
\right|
\right]
\,\d\lambda\,\d u \nonumber
\\
&\qquad=
\left|(2\sigma-1)^{s^2}\right|
\mathbb E\left[
|\zeta(\sigma,X)|^{2\textnormal{Re}(s)}
\right]
\int_{\mathbb R^2}
\left|{}_1F_1\left(
s+1;1;-\frac{\lambda^2+u^2}{4t}
\right)\right|
\,\d\lambda\,\d u<\infty.
\end{align} In fact,
the first factor is finite by (\ref{eq:infinite-random-euler-moment}), that is 
\[
\mathbb E\left[|\zeta(\sigma,X)|^{2\textnormal{Re}(s)}\right]<\infty,
\]
and the second one is finite by \eqref{conclusion}. 
Hence Fubini's theorem allows us to interchange expectation and the
Fourier integral. We obtain
\begin{align*}
&
(2\sigma-1)^{s^2}
\mathbb E\left[
|\zeta(\sigma,X)|^{2s}
|T_\sigma|^{2s}\e^{-t|T_\sigma|^2}
\right]
\\
&\qquad=
\frac{\Gamma(s+1)}{4\pi t^{s+1}}
\int_{\mathbb R^2}
{}_1F_1\left(
s+1;1;-\frac{\lambda^2+u^2}{4t}
\right)
(2\sigma-1)^{s^2}
\mathbb E\left[
|\zeta(\sigma,X)|^{2s}
\e^{\i\lambda\operatorname{Re}T_\sigma
+\i u\operatorname{Im}T_\sigma}
\right]
\,\d\lambda\,\d u.
\end{align*}
By the definition of $H_{\sigma,s}$ in (\ref{eq:Hsigma-def-final}), this is exactly
\eqref{eq:smoothed-moment-fourier}. Finally, the same estimate
\eqref{eq:Fubini-smoothed-moment} shows that the Fourier integral in
\eqref{eq:smoothed-moment-fourier} is absolutely convergent.
\end{proof}

Using the Fourier representation established in
the preceding lemma, together with the locally
uniform convergence in
Proposition~\ref{prop:compact-joint-transform-rigorous}
and the global decay estimate in
Proposition~\ref{prop:global-Fourier-decay-rigorous},
we obtain the following limit of the
Gaussian-regularized complex moments as
$\sigma\downarrow1/2$.

\begin{lem}[Limit of the Gaussian-regularized moments]
\label{lem:fixed-t-smoothed-moment-limit}
Let
$
t>0
$
be fixed.  Let $a_{s}$ be the arithmetic factor as in (\ref{definitionofthearithmeticfactor}).
Then, as $\sigma\downarrow1/2$,
\begin{align}
&(2\sigma-1)^{s^2}
\mathbb E\left[
|\zeta(\sigma,X)|^{2s}
|T_\sigma|^{2s}
\e^{-t|T_\sigma|^2}
\right]\nonumber\\
&\longrightarrow
\mathcal M_t(s)
:=
\frac{a_s\Gamma(s+1)}{4\pi t^{s+1}}
\int_{\mathbb R^2}
{}_1F_1\left(
s+1;1;-\frac{\lambda^2+u^2}{4t}
\right)
\e^{-\i s\lambda-(\lambda^2+u^2)/4}
\,d\lambda\,du,
\label{eq:fixed-t-smoothed-limit}
\end{align}
locally uniformly for
\[
\operatorname{Re}(s)>0.
\]
Here $T_{\sigma}$ is as in (\ref{defofTsigma}).
\end{lem}

\begin{proof}
Let
\[
K\subset
\{s\in\mathbb C:\operatorname{Re}(s)>0\}.
\]
We shall prove that the convergence in
\eqref{eq:fixed-t-smoothed-limit} is uniform for $s\in K$. By Lemma \ref{lem:smoothed-moment-fourier}, it suffices to prove that as $\sigma\downarrow 1/2$,
\begin{equation}
\label{eq:integral-convergence-fixed-t}
\sup_{s\in K}
\left|
\int_{\mathbb R^2}
{}_1F_1\left(
s+1;1;-\frac{\lambda^2+u^2}{4t}
\right)
\left(
H_{\sigma,s}(\lambda,u)-a_s \e^{-\i s\lambda-(\lambda^2+u^2)/4}
\right)
\,d\lambda\,du
\right|
\longrightarrow0.
\end{equation}
We first obtain a bound for ${}_1F_1\left(
s+1;1;-\frac{\lambda^2+u^2}{4t}\right)$ uniformly bounded for 
Fourier variables. By equation (\ref{eq:fst-Fourier}) and (\ref{basicboundforFourier}),
\begin{align}\label{controlF}
\left|{}_1F_1\left(
s+1;1;-\frac{\lambda^2+u^2}{4t}
\right)\right|
&\leq
\frac{
\|f\|_{L^{1}(\mathbb{R}^2)}
}{
\pi|\Gamma(s+1)|t^{-a-1}
},
\end{align}
where
\begin{equation*}
f(x,y)
=(x^2+y^2)^s \e^{-t(x^2+y^2)}.
\end{equation*}
Using polar coordinates,
\begin{align*}
\|f\|_{L^1(\mathbb R^2)}
=\pi t^{-a-1}\Gamma(a+1).
\end{align*}
Then, by (\ref{controlF}),
\begin{equation*}
\left|{}_1F_1\left(
s+1;1;-\frac{\lambda^2+u^2}{4t}
\right)\right|\leq \frac{\Gamma(\operatorname{Re}(s)+1)}
{|\Gamma(s+1)|}.
\end{equation*}
Since $K$ is compact and the Gamma function has no zeros,
\begin{equation}
\label{eq:K-kernel-uniform-K}
\sup_{\substack{
s\in K\\
(\lambda,u)\in\mathbb R^2}}
\left|{}_1F_1\left(
s+1;1;-\frac{\lambda^2+u^2}{4t}
\right)\right|
\leq C_K.
\end{equation}
Then for any $R>0$, by \eqref{eq:K-kernel-uniform-K}, on the disk
\[
B_R
=
\{(\lambda,u)\in\mathbb R^2:
\lambda^2+u^2\leq R^2\},
\]
we have
\begin{align}\label{inthedisk}
&
\sup_{s\in K}
\left|
\int_{B_R}
{}_1F_1\left(
s+1;1;-\frac{\lambda^2+u^2}{4t}
\right)
\left(
H_{\sigma,s}(\lambda,u)-a_s \e^{-\i s\lambda-(\lambda^2+u^2)/4}
\right)
\,d\lambda\,du
\right|\nonumber
\\
&\qquad\leq
C_K |B_R|
\sup_{\substack{
s\in K\\
(\lambda,u)\in B_R}}
\left|\left(
H_{\sigma,s}(\lambda,u)-a_s \e^{-\i s\lambda-(\lambda^2+u^2)/4}
\right)\right|\rightarrow 0,
\end{align}
the last expression tends to zero as $\sigma\downarrow1/2$ Proposition \ref{prop:compact-joint-transform-rigorous}.

On the other hand, applying the global Fourier decay estimate with the fixed
exponent $N=2$, there exist $\sigma_0>1/2$ and $C_K>0$ such that
\begin{equation}
\label{eq:Hsigma-fixed-decay-fixed-t}
|H_{\sigma,s}(\lambda,u)|
\leq
C_K(1+\lambda^2+u^2)^{-2}
\end{equation}
for
\[
\frac12<\sigma\leq\sigma_0,
\qquad
s\in K,
\qquad
(\lambda,u)\in\mathbb R^2.
\]
It is not to check that the limiting function satisfies a bound of the same form, that is, there is $C_K'>0$ such that
\begin{equation}
\label{eq:Hs-fixed-decay-fixed-t}
|a_s \e^{-\i s\lambda-(\lambda^2+u^2)/4}|
\leq
C_K'
(1+\lambda^2+u^2)^{-2},
\qquad
s\in K.
\end{equation}
Then by
\eqref{eq:K-kernel-uniform-K},
\eqref{eq:Hsigma-fixed-decay-fixed-t}, and
\eqref{eq:Hs-fixed-decay-fixed-t}, for any $\epsilon>0$, we may choose $R>0$ so large
that
\begin{align}
\sup_{\substack{
1/2<\sigma\leq\sigma_0\\
s\in K}}
\int_{|(\lambda,u)|>R}
&
|K_{s,t}(\lambda,u)|
\left(
|H_{\sigma,s}(\lambda,u)|
+
|H_s(\lambda,u)|
\right)
\,\d\lambda\,\d u
<\varepsilon.
\label{eq:fixed-t-tail-small}
\end{align}
Combining (\ref{inthedisk}) with
\eqref{eq:fixed-t-tail-small}, we obtain
\eqref{eq:integral-convergence-fixed-t}.
\end{proof}

We now evaluate the limiting Gaussian-regularized
moment obtained in the preceding lemma.
By expressing the corresponding Gaussian integral
in polar coordinates, we obtain an explicit formula
in terms of the confluent hypergeometric function.

\begin{lem}[Evaluation of the limiting Gaussian-regularized moment]
\label{lem:evaluation-fixed-t-limit}
Let
$
s\in\mathbb C$ with 
$\operatorname{Re}(s)>0.
$ Let $t>0$. Let $\mathcal M_t(s)$ be defined in (\ref{eq:fixed-t-smoothed-limit}). Then
\begin{equation}
\label{eq:Mt-explicit}
\mathcal M_t(s)
=
a_s\e^{-s^2}\Gamma(s+1)
(1+t)^{-s-1}
{}_1F_1\left(
s+1;1;\frac{s^2}{1+t}
\right).
\end{equation}
Here $a_{s}$ is the arithmetic factor as in (\ref{definitionofthearithmeticfactor}). Moreover, both sides are holomorphic in
$\operatorname{Re}(s)>0$, and the identity holds locally uniformly
in this half-plane.
\end{lem}

\begin{proof}
Using the fact that for any $u\in \mathbb{R}$,
\[
\int_{\mathbb R}
\e^{\i u y}\e^{-y^2}\,\d y
=
\sqrt{\pi}\,\e^{-u^2/4},
\]
it is not hard to check that

\begin{equation}
\label{eq:gs-Fourier-for-Mt}
\frac{1}{\pi}\int_{\mathbb{R}^2}
\e^{\i(\lambda x+uy)}\e^{-(x+s)^2-y^2}\,\d x\,\d y
=\e^{-\i s\lambda-(\lambda^2+u^2)/4}.
\end{equation}
By Lemma~\ref{lem:fourier-inversion-power-gaussian} that the Fourier transform of 
$
(x^2+y^2)^s \e^{-t(x^2+y^2)}
$
is
\begin{equation}
\label{eq:fst-fourier-for-Mt}
\pi\Gamma(s+1)t^{-s-1}
{}_1F_1\left(
s+1;1;-\frac{\lambda^2+u^2}{4t}
\right).
\end{equation}
In the following, we denote
\[
f_{s,t}(x,y):=(x^2+y^2)^s \e^{-t(x^2+y^2)},\quad g_{s}(x,y)=\frac{1}{\pi}\e^{-(x+s)^2-y^2},
\]
and use the convention
\[
\widehat f(\lambda,u)
=
\int_{\mathbb R^2}
\e^{-\i(\lambda x+uy)}
f(x,y)\,\d x\,\d y.
\]
Hence, by \eqref{eq:gs-Fourier-for-Mt} and
\eqref{eq:fst-fourier-for-Mt}, together with the definition of
$\mathcal M_t(s)$ in \eqref{eq:fixed-t-smoothed-limit}, 
\begin{align*}
\mathcal M_t(s)
&=
\frac{a_{s}}{(2\pi)^2}
\int_{\mathbb R^2}
\widehat f_{s,t}(\lambda,u)
\int_{\mathbb{R}^2}
\left(\e^{\i(\lambda x+uy)}g_{s}(x,y)\,dx\,dy\right)
\,\d \lambda\,\d u.
\end{align*}
Since
\[
f_{s,t},\widehat f_{s,t}\in L^1(\mathbb R^2)
\]
for $\operatorname{Re}(s)>0$, while $g_s\in L^1(\mathbb R^2)$,
Fubini's theorem and Fourier inversion yield
\begin{align}
\mathcal M_t(s)
=
\int_{\mathbb R^2}
f_{s,t}(x,y)g_s(x,y)\,\d x\,\d y
=
\frac{a_s}{\pi}
\int_{\mathbb R^2}
(x^2+y^2)^s
\e^{-t(x^2+y^2)}
\e^{-(x+s)^2-y^2}
\,dx\,dy.
\label{eq:Mt-physical-space}
\end{align}
The integral in \eqref{eq:Mt-physical-space} is absolutely
convergent. Indeed, 
\[
\left|(x^2+y^2)^s\right|
=
(x^2+y^2)^{\textnormal{Re}(s)}
\]
and
\[
\left|
\e^{-(x+s)^2-y^2}
\right|
=
\e^{(\textnormal{Im}(s))^2}\e^{-(x+\textnormal{Re}(s))^2-y^2}.
\]
Hence the absolute value of the integrand is bounded by a polynomial
times a Gaussian. We may therefore pass to polar coordinates,
\[
x=r\cos\theta,
\qquad
y=r\sin\theta,
\qquad
dx\,dy=r\,dr\,\d \theta.
\]
Since
\[
(x+s)^2+y^2
=
r^2+2sr\cos\theta+s^2,
\]
equation \eqref{eq:Mt-physical-space} becomes
\begin{align}
\mathcal M_t(s)
&=
\frac{a_s\e^{-s^2}}{\pi}
\int_0^\infty
r^{2s+1}\e^{-(1+t)r^2}
\left(
\int_0^{2\pi}
\e^{-2sr\cos\theta}\,\d \theta
\right)\d r.
\label{eq:Mt-polar}
\end{align}
Recall that the modified Bessel function of the first kind of
order zero satisfies (see, e.g., \cite[Eq.~(10.32.1)]{NIST})
\[
I_0(z)
=
\frac1{2\pi}
\int_0^{2\pi}
\e^{z\cos\theta}\,d\theta
=
\sum_{k=0}^{\infty}
\frac1{(k!)^2}
\left(\frac z2\right)^{2k},
\qquad z\in\mathbb C.
\]
In particular, $I_0$ is an even function. Hence
\begin{equation}
\label{eq:Mt-I0}
\mathcal M_t(s)
=
2a_s\e^{-s^2}
\int_0^\infty
r^{2s+1}
\e^{-(1+t)r^2}
I_0(2sr)\,\d r.
\end{equation}
Note that the series
\begin{align*}
\sum_{k=0}^{\infty}
\frac{|s|^{2k}}{(k!)^2}
\int_0^\infty
r^{2\textnormal{Re}(s)+2k+1}
\e^{-(1+t)r^2}\,\d r=
\frac12
\sum_{k=0}^{\infty}
\frac{|s|^{2k}}{(k!)^2}
(1+t)^{-\textnormal{Re}(s)-k-1}
\Gamma(\textnormal{Re}(s)+k+1).
\end{align*}
is convergent.This justifies the termwise integration obtained by substituting
the power series expansion of $I_0(2sr)$ into
\eqref{eq:Mt-I0}. Consequently,
\begin{align*}
\mathcal M_t(s)
&=
2a_s\e^{-s^2}
\sum_{k=0}^{\infty}
\frac{s^{2k}}{(k!)^2}
\int_0^\infty
r^{2s+2k+1}
\e^{-(1+t)r^2}\,\d r\\
&=
a_s\e^{-s^2}
\Gamma(s+1)(1+t)^{-s-1}
\sum_{k=0}^{\infty}
\frac{(s+1)_k}{(k!)^2}
\left(
\frac{s^2}{1+t}
\right)^k.
\end{align*}
Finally, by the series expansion of the confluent hypergeometric
function, see, e.g., \cite[Eq.~(13.2.2)]{NIST},
\[
{}_1F_1(a;1;z)
=
\sum_{k=0}^{\infty}
\frac{(a)_k}{(k!)^2}z^k,
\]
we conclude that \eqref{eq:Mt-explicit}. 
Since all the preceding interchanges and identities are locally
uniform in $s$, \eqref{eq:Mt-explicit} holds locally uniformly in
the half-plane
$
\operatorname{Re}(s)>0.
$
\end{proof}

It is important to note that the limit $t\downarrow0$ cannot be
interchanged directly with the limit $\sigma\downarrow1/2$.
Indeed, for each fixed $\sigma>1/2$, dominated convergence gives
\[
\mathcal M_{\sigma,t}(s):=(2\sigma-1)^{s^2}
\mathbb E\left[
|\zeta(\sigma,X)|^{2s}
|T_\sigma|^{2s}
\e^{-t|T_\sigma|^2}
\right]
\longrightarrow
\mathcal M_{\sigma,0}(s),
\qquad t\downarrow0,
\]
where
\[
T_\sigma=(2\sigma-1) L'(\sigma,X),\quad \mathcal M_{\sigma,0}(s)
=
(2\sigma-1)^{s^2}
\mathbb E\left[
|\zeta(\sigma,X)|^{2s}|T_\sigma|^{2s}
\right].
\]
However, to interchange the two limits one would need a bound of the
form
\[
\sup_{\sigma\downarrow1/2}
\left|
\mathcal M_{\sigma,0}(s)-\mathcal M_{\sigma,t}(s)
\right|
\longrightarrow0,
\qquad t\downarrow0.
\]
Note that
\begin{align*}
&
\left|
\mathcal M_{\sigma,0}(s)-\mathcal M_{\sigma,t}(s)
\right|\leq
(2\sigma-1)^{(\textnormal{Re}(s))^2-(\textnormal{Im}(s))^2}
\mathbb E\left[
|\zeta(\sigma,X)|^{2\textnormal{Re}(s)}
|T_\sigma|^{2\textnormal{Re}(s)}
\left(1-\e^{-t|T_\sigma|^2}\right)
\right].
\end{align*}
Taking absolute values replaces the complex weight
$|\zeta|^{2s}|T_\sigma|^{2s}$ by the positive weight
$|\zeta|^{2a}|T_\sigma|^{2a}$ and therefore destroys the
cancellation present in the original complex expectation. One obtains only
\[
\left|
\mathcal M_{\sigma,0}(s)-\mathcal M_{\sigma,t}(s)
\right|
\leq
t\,(2\sigma-1)^{(\textnormal{Re}(s))^2-(\textnormal{Im}(s))^2}
\mathbb E\left[
|\zeta(\sigma,X)|^{2\textnormal{Re}(s)}
|T_\sigma|^{2\textnormal{Re}(s)+2}
\right].
\]
Even if the corresponding real-weighted moment
satisfies the natural bound
\[
(2\sigma-1)^{a^2}
\mathbb E\left[
|\zeta(\sigma,X)|^{2a}|T_\sigma|^{2a+2}
\right]
\ll 1,
\]
the preceding estimate yields only
\[
\left|
\mathcal M_{\sigma,0}(s)-\mathcal M_{\sigma,t}(s)
\right|
\ll t(2\sigma-1)^{-b^2},
\]
which is not uniform as $\sigma\downarrow1/2$
when $b\neq0$.
Thus a direct absolute-value argument does not provide the uniform
control required to interchange the limits.  This is the reason for
using the regularized Fourier--Riesz representation, which preserves
the analytic cancellation in the complex parameter $s$.

For complex $s$, rather than interchanging the limits
$t\downarrow0$ and $\sigma\downarrow1/2$,
we first fix $\sigma>1/2$ and remove the Gaussian
cutoff by letting $t\downarrow0$.
Using the theory of Riesz distributions and
their meromorphic continuation, we then obtain
a Fourier representation of the unregularized
complex moments. This representation allows us
to subsequently take the limit $\sigma\downarrow1/2$.
We begin by removing the Gaussian cutoff
in the space of tempered distributions.

\subsection{Riesz distributions and Fourier--Riesz identities}

We first observe that, for every $w\in\mathbb C$
with $\operatorname{Re}(w)>-1$,
\begin{equation}
\label{eq:fw-Sprime-limit}
|x|^{2w}\e^{-t|x|^2}
\longrightarrow
|x|^{2w}
\qquad\text{in }\mathcal S'(\mathbb R^2),
\qquad t\downarrow0.
\end{equation}
Indeed, for every $\phi\in\mathcal S(\mathbb R^2)$,
\[
\left|
|x|^{2w}
\bigl(\e^{-t|x|^2}-1\bigr)\phi(x)
\right|
\leq
2|x|^{2\operatorname{Re}(w)}|\phi(x)|.
\]
Since the right-hand side is integrable
over $\mathbb R^2$ when
$\operatorname{Re}(w)>-1$,
the convergence follows from the dominated
convergence theorem.
\subsubsection{Fourier transform of homogeneous distributions}
Since the Fourier transform is a continuous automorphism of
$\mathcal S'(\mathbb R^2)$,
\eqref{eq:fw-Sprime-limit} implies
\[
\widehat {|x|^{2w}\e^{-t|x|^2}}
\longrightarrow
\widehat{|x|^{2w}}
\qquad\text{in }\mathcal S'(\mathbb R^2).
\]
The theory of homogeneous distributions and
Riesz distributions is classical; see
\cite[Chapter~III, Section~3.2]{Hormander1983}
and \cite[Chapter~V, Section~1]{Stein1970}.
In particular, the meromorphic continuation
of radial powers $r^\lambda$ is treated in
\cite[Chapter~I, \S3.9, Eq.~(3) and the
discussion immediately following it]{GelfandShilov},
while their Fourier transforms are computed
in \cite[Chapter~II, \S3.3, Eq.~(2)]{GelfandShilov}.
Thus, with $n=2$ and $\lambda=2w$, we obtain the Fourier--Riesz identity
\eqref{1023evening}, that is 
\begin{equation}\label{other1023evening}
\widehat{|x|^{2w}}(\xi)
=
2^{2w+2}\pi
\frac{\Gamma(w+1)}{\Gamma(-w)}
\,\operatorname{Reg}|\xi|^{-2w-2},
\end{equation}
where $\operatorname{Reg}|\xi|^{-2w-2}$
denotes the meromorphic continuation of
the homogeneous distribution
$|\xi|^{-2w-2}$, which is locally
integrable for $\operatorname{Re}w<0$.
The strip
\[
-1<\operatorname{Re}w<0
\]
is the common range in which both
$|x|^{2w}$ and $|\xi|^{-2w-2}$
are represented by locally integrable functions.

To make the identity
\eqref{other1023evening} explicit for the subsequent
analysis of complex moments, we now construct
the meromorphic continuation of the Riesz
distribution through Taylor subtraction at the
origin. 
This yields a regularized functional
that provides an explicit representation
of the pairing between the Riesz distribution
and the weighted joint Fourier transform
$H_{\sigma,s}$.
\noindent
\subsubsection{Taylor-subtraction realization of Riesz distributions}

The following construction is the two-dimensional specialization
of the meromorphic continuation of the radial homogeneous
distribution $r^\lambda$ developed by Gelfand and Shilov;
see \cite[Chapter~I, \S3.9, pp.~71--74]{GelfandShilov}.
In particular, the spherical mean appearing there becomes the
circular mean introduced below, while the corresponding
Pizzetti formula gives the coefficients involving the iterated
Laplacian. 
For $H\in\mathcal S(\mathbb R^2)$,
recall that $\mathcal A H(r)$ denotes
the circular average of $H$, as defined
in \eqref{defofcirmean}. For $-1<\operatorname{Re}w<0$, whenever the pairing is
absolutely convergent, polar coordinates give
\[
\left\langle |\xi|^{-2w-2},H\right\rangle
=
2\pi
\int_0^\infty
\mathcal A H(r)r^{-2w-1}\,\d r.
\]
This is the two-dimensional specialization of the spherical-mean
representation of $r^\lambda$ in
\cite[Chapter~I, \S3.9, Eq.~(2)]{GelfandShilov},
with
\[
\lambda=-2w-2.
\]
In what follows, we use the notation
\[
\partial_\lambda:=\frac{\partial}{\partial\lambda},
\qquad
\partial_u:=\frac{\partial}{\partial u}.
\]
Thus,
\[
\cos\theta\,\partial_\lambda+\sin\theta\,\partial_u
\]
denotes the directional derivative in the direction
$(\cos\theta,\sin\theta)$. Let $M\geq0$. Taylor's theorem at the origin gives
\[
H(r\cos\theta,r\sin\theta)
=
\sum_{k=0}^{2M+1}
\frac{r^k}{k!}
\left(
\cos\theta\,\partial_\lambda
+
\sin\theta\,\partial_u
\right)^k H(0)
+
O(r^{2M+2}),
\qquad r\downarrow0,
\]
where the implied constant depends only on $M$ and on a
$C^{2M+2}$-bound for $H$ in a fixed neighbourhood of the origin,
and is uniform in $\theta\in[0,2\pi]$. After averaging over $\theta$, all the odd-order terms vanish.
Moreover,
\[
\frac1{2\pi}
\int_0^{2\pi}
\left(
\cos\theta\,\partial_\lambda
+
\sin\theta\,\partial_u
\right)^{2j}
\,\d \theta
=
\frac{(2j)!}{4^j(j!)^2}\Delta^j.
\]
Here recall that $\Delta=\partial_\lambda^2+\partial_u^2$ denotes the Laplacian in the variables $(\lambda,u)$. 
Consequently,
\begin{equation}
\label{eq:circular-average-Taylor}
\mathcal A H(r)
=
\sum_{j=0}^M
\frac{r^{2j}}{4^j(j!)^2}
\Delta^j H(0)
+
O(r^{2M+2}),
\qquad r\downarrow0.
\end{equation}
This is precisely the two-dimensional form of the spherical-mean
expansion in
\cite[Chapter~I, \S3.9, Eq.~(7)]{GelfandShilov}. 
The original radial Riesz integral
\[
\int_0^\infty \mathcal A H(r)r^{-2w-1}\,\d r
\]
is absolutely convergent near the origin only in the initial strip
$-1<\operatorname{Re}w<0$.  To continue it to larger values of
$\operatorname{Re}w$, we subtract sufficiently many terms of the
Taylor expansion of the circular average at the origin and add
back their explicitly integrated contributions.  This leads to the
following regularized Riesz functional. This is the two-dimensional specialization
of the Taylor-subtraction procedure for algebraic singularities
developed by Gelfand and Shilov
\cite[Chapter~I, Section~3.2, Eq.~(3)]{GelfandShilov}, combined with their
spherical-mean representation of the radial homogeneous
distribution $r^\lambda$ in
\cite[Chapter~I, \S3.9]{GelfandShilov} (here we set $\lambda=-2w-1$).
For
\[
-1<\operatorname{Re}w<M+1,
\qquad
w\notin\{0,1,\ldots,M\},
\]
let the regularized Riesz functional
$\mathcal R_{M,w}[H]$ be defined in
\eqref{eq:Riesz-functional-abstract}. In particular,  for
\[
-1<\operatorname{Re}w<0,
\]
no regularization is needed at the origin. Indeed,
\[
\int_0^1 r^{2j-2w-1}\,dr
=
\frac1{2(j-w)},
\]
and therefore
\begin{equation}
\label{eq:Riesz-functional-initial}
\mathcal R_{M,w}[H]
=
\int_0^\infty
\mathcal AH(r)r^{-2w-1}\,\d r.
\end{equation}
Thus, using polar coordinates,
\[
\begin{aligned}
\left\langle
|\xi|^{-2w-2},H
\right\rangle
=
\int_{\mathbb R^2}
|\xi|^{-2w-2}H(\xi)\,\d\xi=
2\pi
\int_0^\infty
\mathcal A H(r)r^{-2w-1}\,\d r=
2\pi\mathcal R_{M,w}[H],
\end{aligned}
\]
in the initial strip $-1<\operatorname{Re}w<0$.

The Taylor subtraction in
\eqref{eq:Riesz-functional-abstract} is the explicit
two-dimensional realization of the analytic-continuation
procedure for $r^{-2w-2}$ in
\cite[Chapter~I, \S3.9]{GelfandShilov}.
Indeed, by \eqref{eq:circular-average-Taylor},
\[
\mathcal A H(r)
-
\sum_{j=0}^M
\frac{r^{2j}}{4^j(j!)^2}\Delta^jH(0)
=
O(r^{2M+2}),
\]
so that
\[
\left[
\mathcal A H(r)
-
\sum_{j=0}^M
\frac{r^{2j}}{4^j(j!)^2}\Delta^jH(0)
\right]
r^{-2w-1}
=
O\!\left(r^{2M+1-2\operatorname{Re}w}\right).
\]
Hence the first integral in
\eqref{eq:Riesz-functional-abstract} is locally uniformly
absolutely convergent, and therefore holomorphic in $w$, for
\[
\operatorname{Re}w<M+1.
\]
The correction terms
\[
\frac{\Delta^jH(0)}
{2\,4^j(j!)^2(j-w)}
\]
give the meromorphic continuation of the subtracted Taylor
monomials and exhibit the possible simple poles at
\[
w=0,1,\ldots,M.
\]
This agrees with the pole structure of the meromorphic
continuation of $r^{-2w-2}$ described in
\cite[Chapter~I, \S3.9, p.~72]{GelfandShilov}. Provided that $H$ has sufficient decay at infinity, the second
integral in \eqref{eq:Riesz-functional-abstract} is also locally
uniformly convergent in the range under consideration.

By the Gelfand-Shilov meromorphic continuation of radial
homogeneous distributions \cite[Chapter~I, \S3.9]{GelfandShilov}, the map
\[
w\longmapsto
\left\langle
\operatorname{Reg}|\xi|^{-2w-2},H
\right\rangle
\]
is meromorphic in the strip under consideration. Here $H\in\mathcal S(\mathbb R^2)$. On the other
hand, the first two integrals defining
$\mathcal R_{M,w}[H]$ are holomorphic there by locally uniform
absolute convergence, while the finite correction sum is
meromorphic with possible simple poles at
$w=0,1,\ldots,M$.  Hence both sides are meromorphic functions
of $w$. Since they agree for $-1<\operatorname{Re}w<0$, uniqueness of
meromorphic continuation yields
\begin{equation}\label{anaexten}
\left\langle
\operatorname{Reg}|\xi|^{-2w-2},H
\right\rangle
=
2\pi\mathcal R_{M,w}[H].
\end{equation}

Recall the circular average $\mathcal A H(r)$
introduced for Schwartz functions in
\eqref{defofcirmean}.
This definition extends naturally to functions
$H\in C(\mathbb R^2)$
by the same formula:
\begin{equation}
\label{def:circular-mean2}
\mathcal A H(r)
:=
\frac{1}{2\pi}
\int_0^{2\pi}
H(r\cos\theta,r\sin\theta)
\,\d\theta,
\qquad r\geq0.
\end{equation}
We also recall the regularized Riesz
functional $\mathcal R_{M,w}[H]$
defined in \eqref{eq:Riesz-functional-abstract}. Note that it involves derivatives of $H$
at the origin. Hence, for
$M\in\mathbb Z_{\geq0}$, formula
\eqref{eq:Riesz-functional-abstract}
is well-defined whenever
\[
H\in C^{2M+2}(\mathbb R^2)
\cap L^1(\mathbb R^2),
\]
and
\[
-1<\operatorname{Re}(w)<M+1,
\qquad
w\notin\{0,1,\ldots,M\}.
\]
Indeed, Taylor's theorem gives
\[
\mathcal A H(r)
-
\sum_{j=0}^M
\frac{r^{2j}}{4^j(j!)^2}
\Delta^jH(0)
=
O(r^{2M+2}),
\qquad r\downarrow0,
\]
which ensures convergence of the regularized
integral at the origin whenever
$\operatorname{Re}(w)<M+1$.
The $L^1$-assumption, together with
$\operatorname{Re}(w)>-1$, guarantees
absolute convergence of the integral
at infinity.

The following proposition shows that the regularized
functional $\mathcal R_{M,w}[H]$ agrees with an
extension of the Riesz distributional pairing
to functions in
$C^{2M+2}(\mathbb R^2)\cap L^1(\mathbb R^2)$.
In particular, this extension allows us to
express the original complex moments,
without Gaussian regularization, in terms of
the Riesz functional $\mathcal R_{M,w}[H_{\sigma,s}]$.

\begin{prop}
\label{lem:Reg-Riesz-Hsigma}
Let \(M\in\mathbb Z_{\geq0}\). Suppose that
\[
H\in
C^\infty(\mathbb R^2)\cap L^1(\mathbb R^2).
\]
Define
\[
\Omega_M
:=
\left\{
w\in\mathbb C:
-1<\operatorname{Re}w<M+1
\right\}
\setminus\{0,1,\ldots,M\}.
\]
Let $\mathcal R_{M,w}[H]$ be defined in (\ref{eq:Riesz-functional-abstract}) for $w\in \Omega_M$.
Then the Gelfand--Shilov pairing with the regularized homogeneous
distribution admits an extension to \(H\) on
\(\Omega_M\), given by
\begin{equation}
\label{eq:Reg-Riesz-Hsigma-lemma}
\left\langle
\operatorname{Reg}|\xi|^{-2w-2},
H
\right\rangle_{\mathrm{ext}}
=
2\pi\mathcal R_{M,w}[H],
\qquad
w\in\Omega_M.
\end{equation}
The extended pairing is independent of the cutoff used in its
definition and is holomorphic in \(w\) on \(\Omega_M\).
\end{prop}

\begin{proof}
Choose a radial function
\(\chi\in C_c^\infty(\mathbb R^2)\) such that
\[
0\leq\chi\leq1,
\qquad
\chi(\xi)=1
\quad\text{for }|\xi|\leq1.
\]
For \(R\geq1\), set
\[
\chi_R(\xi):=\chi(\xi/R),
\qquad
H_R(\xi):=\chi_R(\xi)H(\xi).
\]
Then
\[
H_R\in C_c^\infty(\mathbb R^2)
\subset\mathcal S(\mathbb R^2),
\]
and hence the Gelfand--Shilov distribution
\(\operatorname{Reg}|\xi|^{-2w-2}\) acts on \(H_R\) in the usual
distributional sense. By (\ref{anaexten}),
\begin{equation}
\label{eq:cutoff-Reg-Riesz}
\left\langle
\operatorname{Reg}|\xi|^{-2w-2},
H_R
\right\rangle
=
2\pi\mathcal R_{M,w}[H_R],
\qquad
w\in\Omega_M
\end{equation}
and the later one is holomorphic in $w$ on $\Omega_{M}$. We next pass to the limit \(R\to\infty\). Since
\(\chi_R=1\) on \(B(0,1)\) for every \(R\geq1\),
\[
H_R(\xi)=H(\xi),
\qquad |\xi|\leq1.
\]
Consequently,
\[
\mathcal A H_R(r)
=
\mathcal A H(r),
\qquad 0\leq r\leq1,
\]
and
\[
\Delta^jH_R(0)
=
\Delta^jH(0),
\qquad 0\leq j\leq M.
\]
Thus the first and third terms in
\(\mathcal R_{M,w}[H_R]\) coincide exactly with the corresponding
terms in
\(\mathcal R_{M,w}[H]\). It therefore remains only to consider the tail integral. Let  \(K\) be a compact set contained in
\(\Omega_M\). 
For \(r\geq1\),  
\[
r^{-2\operatorname{Re}w-2}
\leq r^{-2\inf_{w\in K}\operatorname{Re}w-2}
\leq1.
\]
Hence
\begin{align}
\sup_{w\in K}
\left|
\int_1^\infty
\mathcal A(H_R-H)(r)
r^{-2w-1}\,\d r
\right|
&\leq
\int_1^\infty
\left|
\mathcal A(H_R-H)(r)
\right|r\,\d r\nonumber\\
&\leq
\frac{1}{2\pi}
\int_{|\xi|\geq1}
|H_R(\xi)-H(\xi)|\,\d \xi.
\label{eq:tail-cutoff-estimate}
\end{align}
Moreover,
\[
H_R(\xi)\longrightarrow H(\xi)
\qquad
(R\to\infty)
\]
for every \(\xi\), and, since \(0\leq\chi_R\leq1\),
\[
|H_R(\xi)-H(\xi)|
\leq
|H(\xi)|.
\]
Because
\(H\in L^1(\mathbb R^2)\), dominated convergence gives
\[
\|H_R-H\|_{L^1(\mathbb R^2)}
\longrightarrow0.
\]
It follows from \eqref{eq:tail-cutoff-estimate} that
\begin{equation}
\label{eq:R-functional-cutoff-limit}
\mathcal R_{M,w}[H_R]
\longrightarrow
\mathcal R_{M,w}[H]
\end{equation}
locally uniformly for \(w\in\Omega_M\). This motivates the following extension of the Gelfand--Shilov
pairing to \(H\):
\begin{equation}
\label{eq:extended-Reg-pairing}
\left\langle
\operatorname{Reg}|\xi|^{-2w-2},
H
\right\rangle_{\mathrm{ext}}
:=
\lim_{R\to\infty}
\left\langle
\operatorname{Reg}|\xi|^{-2w-2},
H_R
\right\rangle .
\end{equation}
The limit exists locally uniformly in \(w\in\Omega_M\).
Combining \eqref{eq:cutoff-Reg-Riesz} with
\eqref{eq:R-functional-cutoff-limit} yields
\[
\left\langle
\operatorname{Reg}|\xi|^{-2w-2},
H
\right\rangle_{\mathrm{ext}}
=
2\pi\mathcal R_{M,w}[H],
\]
which proves \eqref{eq:Reg-Riesz-Hsigma-lemma}.
The right-hand side depends only on \(H\), and not on the
particular cutoff function \(\chi\). Hence the extended pairing is
independent of the choice of cutoff. Finally, since for every \(R\)
the function
\[
w\longmapsto
\left\langle
\operatorname{Reg}|\xi|^{-2w-2},H_R
\right\rangle
\]
is holomorphic on \(\Omega_M\), and the convergence in
\eqref{eq:extended-Reg-pairing} is locally uniform there, the
Weierstrass theorem implies that the extended pairing is holomorphic
on \(\Omega_M\).
\end{proof}

The extended Riesz distributional pairing established
in Proposition~\ref{lem:Reg-Riesz-Hsigma}
will allow us to express the weighted mixed moments
of $\zeta(\sigma,X)$ and $L'(\sigma,X)$ in terms of
the regularized Riesz functional
$\mathcal R_{M,w}[H_{\sigma,s}]$.
To establish this representation, we first prove
the integrability of the relevant mixed moments
and then introduce the weighted measure
$\nu_{\sigma,s}$ defined in \eqref{defofthemeasure}. 
We subsequently derive the corresponding
Fourier--Riesz identity. We first establish the integrability of the
corresponding mixed moments.

\noindent
\subsubsection{The Fourier--Riesz identity for the weighted measure}

\begin{lem}
\label{lem:mixed-moment}
Let \(\sigma>1/2\), \(s\in\mathbb C\) with
\(\operatorname{Re}s>0\), and \(w\in\mathbb C\) with
$
\operatorname{Re}w>-1.
$ Let $T_{\sigma}$ be as in (\ref{defofTsigma}). 
Then
\begin{equation}
\label{eq:mixed-moment-finite}
\mathbb E\left[
|\zeta(\sigma,X)|^{2\operatorname{Re}s}
|T_\sigma|^{2\operatorname{Re}w}
\right]
<\infty.
\end{equation}
Consequently, for every \(t>0\),
\begin{equation}
\label{eq:mixed-moment-gaussian-finite}
\mathbb E\left[
|\zeta(\sigma,X)|^{2\operatorname{Re}s}
|T_\sigma|^{2\operatorname{Re}w}
\mathrm e^{-t|T_\sigma|^2}
\right]
<\infty.
\end{equation}
\end{lem}

\begin{proof}
We distinguish two cases. Suppose first that \(-1<\textnormal{Re}(w)<0\). Since
\(H_{\sigma,\textnormal{Re}(s)}\in L^1(\mathbb R^2)\) by (\ref{eq:global-Fourier-decay-rigorous}), Fourier inversion yields a
bounded continuous density
\[
h_{\sigma,\textnormal{Re}(s)}(x)
=
\frac{1}{(2\pi)^2}
\int_{\mathbb R^2}
\mathrm e^{\i \langle x,\xi\rangle}
H_{\sigma,\textnormal{Re}(s)}(\xi)\,\d\xi
\]
such that, for every nonnegative Borel function \(G\),
\[
(2\sigma-1)^{(\textnormal{Re}(s))^2}
\mathbb E\left[
|\zeta(\sigma,X)|^{2\textnormal{Re}(s)}G(T_\sigma)
\right]
=
\int_{\mathbb R^2}G(x)h_{\sigma,\textnormal{Re}(s)}(x)\,\d x,
\]
whenever either side is finite, with the identity for general
nonnegative \(G\) obtained by monotone convergence. Taking
\[
G(x)=|x|^{2\textnormal{Re}(w)}\mathbf 1_{\{|x|<1\}},
\]
we obtain
\[
\begin{aligned}
&(2\sigma-1)^{(\textnormal{Re}(s))^2}
\mathbb E\left[
|\zeta(\sigma,X)|^{2\textnormal{Re}(s)}
|T_\sigma|^{2\textnormal{Re}(w)}
\mathbf 1_{\{|T_\sigma|<1\}}
\right]\\
&\qquad=
\int_{|x|<1}|x|^{2\textnormal{Re}(w)}h_{\sigma,\textnormal{Re}(s)}(x)\,dx\\
&\qquad\leq
2\pi\|h_{\sigma,\textnormal{Re}(s)}\|_{\infty}
\int_0^1 r^{2\textnormal{Re}(w)+1}\,dr
<\infty,
\end{aligned}
\]
since \(\textnormal{Re}(w)>-1\). On the other hand, since \(\textnormal{Re}(w)<0\),
\[
|T_\sigma|^{2\textnormal{Re}(w)}\leq1
\qquad\text{on }\{|T_\sigma|\geq1\},
\]
and therefore by (\ref{eq:fixed-sigma-all-moments}),
\[
\begin{aligned}
\mathbb E\left[
|\zeta(\sigma,X)|^{2\textnormal{Re}(s)}
|T_\sigma|^{2\textnormal{Re}(w)}
\mathbf 1_{\{|T_\sigma|\geq1\}}
\right]
&\leq
\mathbb E\left[
|\zeta(\sigma,X)|^{2\textnormal{Re}(s)}
\right]
<\infty.
\end{aligned}
\]
Thus \eqref{eq:mixed-moment-finite} holds when
\(-1<\textnormal{Re}(w)<0\). Suppose now that \(\Re(w)\geq0\). If \(\textnormal{Re}(w)=0\), then
\[
\mathbb E\left[
|\zeta(\sigma,X)|^{2\textnormal{Re}(s)}|T_\sigma|^{2\textnormal{Re}(w)}
\right]
=
\mathbb E\left[
|\zeta(\sigma,X)|^{2\textnormal{Re}(s)}
\right]
<\infty.
\]
If \(\textnormal{Re}(w)>0\), using
\[
T_\sigma=(2\sigma-1)L'(\sigma,X)
\] and (\ref{eq:fixed-sigma-all-moments}), together with
the Cauchy--Schwarz inequality gives
\[
\begin{aligned}
&\mathbb E\left[
|\zeta(\sigma,X)|^{2\textnormal{Re}(s)}|T_\sigma|^{2\textnormal{Re}(w)}
\right]\\
&\leq
(2\sigma-1)^{2\textnormal{Re}(w)}
\left(
\mathbb E\left[
|\zeta(\sigma,X)|^{4\textnormal{Re}(s)}
\right]
\right)^{1/2}
\left(
\mathbb E\left[
|L'(\sigma,X)|^{4\textnormal{Re}(w)}
\right]
\right)^{1/2}<\infty,
\end{aligned}
\]
This proves
\eqref{eq:mixed-moment-finite} for all \(\textnormal{Re}(w)>-1\). Finally, since
\[
0<\mathrm e^{-t|T_\sigma|^2}\leq1,
\]
we have
\[
|\zeta(\sigma,X)|^{2\textnormal{Re}(s)}
|T_\sigma|^{2\textnormal{Re}(w)}
\mathrm e^{-t|T_\sigma|^2}
\leq
|\zeta(\sigma,X)|^{2\textnormal{Re}(s)}
|T_\sigma|^{2\textnormal{Re}(w)},
\]
and hence \eqref{eq:mixed-moment-gaussian-finite} follows immediately
from \eqref{eq:mixed-moment-finite}.
\end{proof}

The integrability established in
Lemma~\ref{lem:mixed-moment} allows us to
study the weighted mixed moments through
a complex measure on $\mathbb R^2$.
To this end, we introduce the following
weighted measure, whose Fourier transform
is given by $H_{\sigma,s}$.

Recall the weighted complex Borel measure
$\nu_{\sigma,s}$ introduced in
\eqref{defofthemeasure}.
We first verify that this measure is finite. Fix $\sigma>1/2$ and $s\in\mathbb C$
with $\operatorname{Re}(s)>0$.
By the moment estimate
\eqref{eq:fixed-sigma-all-moments},
\[
\mathbb E\left[
|\zeta(\sigma,X)|^{2\operatorname{Re}(s)}
\right]<\infty.
\]
Moreover, for every Borel set
$B\subset\mathbb R^2$,
\[
|\nu_{\sigma,s}(B)|
\leq
(2\sigma-1)^{\operatorname{Re}(s^2)}
\mathbb E\left[
|\zeta(\sigma,X)|^{2\operatorname{Re}(s)}
\mathbf 1_{\{T_\sigma\in B\}}
\right].
\]
Consequently, the total variation measure
satisfies
\[
|\nu_{\sigma,s}|(\mathbb R^2)
\leq
(2\sigma-1)^{\operatorname{Re}(s^2)}
\mathbb E\left[
|\zeta(\sigma,X)|^{2\operatorname{Re}(s)}
\right]
<\infty.
\]
Hence $\nu_{\sigma,s}$ is a finite
complex Borel measure.

The measure $\nu_{\sigma,s}$ thus encodes the
weighted distribution of $T_\sigma$. Recall that $T_{\sigma}$ denotes $(2\sigma-1)L'(\sigma,X)$.
The following lemma establishes the corresponding
integration formula, allowing weighted expectations
of Borel functions of $T_\sigma$ to be expressed
as integrals with respect to $\nu_{\sigma,s}$.

\begin{lem}
\label{lem:weighted-pushforward-nu}
Fix \(\sigma>1/2\) and \(s\in\mathbb C\) with
\(\operatorname{Re}s>0\). Let $T_{\sigma}$ be as in (\ref{defofTsigma}). Let $\nu_{\sigma,s}$ be the measure defined in (\ref{defofthemeasure}).
Let \(G:\mathbb R^2\to\mathbb C\) be a Borel measurable function
such that
\begin{equation}
\label{eq:G-weighted-integrability}
\mathbb E\left[
|\zeta(\sigma,X)|^{2\operatorname{Re}s}
|G(T_\sigma)|
\right]<\infty.
\end{equation}
Then \(G\in L^1(\mathbb R^2,d|\nu_{\sigma,s}|)\), and
\begin{equation}
\label{eq:weighted-pushforward-identity}
\int_{\mathbb R^2}
G(x)\,d\nu_{\sigma,s}(x)
=
(2\sigma-1)^{s^2}
\mathbb E\left[
|\zeta(\sigma,X)|^{2s}G(T_\sigma)
\right].
\end{equation}
\end{lem}

\begin{proof}
We first consider a Borel simple function
\[
G(x)=\sum_{j=1}^{N}c_j\mathbf 1_{B_j}(x),
\]
where the \(B_j\)'s are Borel sets. By the definition of
\(\nu_{\sigma,s}\),
\[
\begin{aligned}
\int_{\mathbb R^2}G(x)\,d\nu_{\sigma,s}(x)
&=
\sum_{j=1}^{N}c_j\nu_{\sigma,s}(B_j)\\
&=
(2\sigma-1)^{s^2}
\sum_{j=1}^{N}c_j
\mathbb E\left[
|\zeta(\sigma,X)|^{2s}
\mathbf 1_{\{T_\sigma\in B_j\}}
\right]\\
&=
(2\sigma-1)^{s^2}
\mathbb E\left[
|\zeta(\sigma,X)|^{2s}G(T_\sigma)
\right].
\end{aligned}
\]
Thus \eqref{eq:weighted-pushforward-identity} holds for Borel
simple functions. Moreover, for every nonnegative Borel simple function \(\Phi\),
the definition of the total variation measure gives
\[
\int_{\mathbb R^2}\Phi(x)\,\d|\nu_{\sigma,s}|(x)
\leq
(2\sigma-1)^{\operatorname{Re}(s^2)}
\mathbb E\left[
|\zeta(\sigma,X)|^{2\operatorname{Re}s}
\Phi(T_\sigma)
\right].
\]
By monotone convergence, the same inequality holds for every
nonnegative Borel function \(\Phi\). Taking
\(\Phi=|G|\) and using \eqref{eq:G-weighted-integrability}, we obtain
\[
\int_{\mathbb R^2}|G(x)|\,\d|\nu_{\sigma,s}|(x)
\leq
(2\sigma-1)^{\operatorname{Re}(s^2)}
\mathbb E\left[
|\zeta(\sigma,X)|^{2\operatorname{Re}s}
|G(T_\sigma)|
\right]
<\infty.
\]
Hence
\[
G\in L^1(\mathbb R^2,d|\nu_{\sigma,s}|).
\]
Now choose Borel simple functions \(G_n\) such that
\[
G_n(x)\longrightarrow G(x)
\]
pointwise and
\[
|G_n(x)|\leq |G(x)|,
\qquad x\in\mathbb R^2.
\]
For each \(n\), the simple-function case gives
\[
\int_{\mathbb R^2}G_n(x)\,\d\nu_{\sigma,s}(x)
=
(2\sigma-1)^{s^2}
\mathbb E\left[
|\zeta(\sigma,X)|^{2s}G_n(T_\sigma)
\right].
\]
Since
\[
|G_n(x)|
\leq |G(x)|
\in L^1(\d|\nu_{\sigma,s}|),
\]
dominated convergence on the left yields
\[
\int_{\mathbb R^2}G_n(x)\,\d\nu_{\sigma,s}(x)
\longrightarrow
\int_{\mathbb R^2}G(x)\,\d\nu_{\sigma,s}(x).
\]
On the other hand,
\[
\left|
|\zeta(\sigma,X)|^{2s}G_n(T_\sigma)
\right|
\leq
|\zeta(\sigma,X)|^{2\operatorname{Re}s}
|G(T_\sigma)|,
\]
which is integrable by \eqref{eq:G-weighted-integrability}.
Hence another application of the dominated convergence theorem gives
\[
\mathbb E\left[
|\zeta(\sigma,X)|^{2s}G_n(T_\sigma)
\right]
\longrightarrow
\mathbb E\left[
|\zeta(\sigma,X)|^{2s}G(T_\sigma)
\right].
\]
Passing to the limit proves
\eqref{eq:weighted-pushforward-identity}.
\end{proof}

We next establish the smoothness and integrability
of $H_{\sigma,s}$, which are needed to apply
the extended Riesz distributional pairing
from Proposition~\ref{lem:Reg-Riesz-Hsigma}.
The proof relies on the moment estimates
established earlier and the global Fourier decay
in Proposition~\ref{prop:global-Fourier-decay-rigorous}.

\begin{lem}
\label{lem:H-smooth-integrable}
Let \(K\) be a compact subset of the open right half-plane
\(\{s\in\mathbb C:\operatorname{Re}s>0\}\).
There exists \(\sigma_0>1/2\) such that, for every
\[
\frac12<\sigma\leq\sigma_0,
\qquad
s\in K,
\]
one has
\[
H_{\sigma,s}\in
C_b^\infty(\mathbb R^2)\cap L^1(\mathbb R^2)
\]
with $H_{\sigma,s}$ defined in (\ref{eq:Hsigma-def-final}). 
Moreover,
\[
\sup_{\substack{1/2<\sigma\leq\sigma_0\\ s\in K}}
\|H_{\sigma,s}\|_{L^1(\mathbb R^2)}
<\infty.
\]
Here \(C_b^\infty(\mathbb R^2)\) denotes the space of infinitely
differentiable functions on \(\mathbb R^2\) whose derivatives of all
orders are bounded.
\end{lem}
\begin{proof}
Recall that $T_{\sigma}$ is given in (\ref{defofTsigma}), that is 
\begin{equation*}
T_{\sigma}=(2\sigma-1)L'(\sigma,X).
\end{equation*}
Write
\[
T_\sigma=(T_{\sigma,1},T_{\sigma,2})
\]
and recall that in (\ref{eq:Hsigma-def-final}),
\[
H_{\sigma,s}(\xi)
=
(2\sigma-1)^{s^2}
\mathbb E\left[
|\zeta(\sigma,X)|^{2s}
\mathrm e^{-\i\langle\xi,T_\sigma\rangle}
\right].
\]
Let
$
\alpha=(\alpha_1,\alpha_2)\in\mathbb Z_{\geq0}^2
$
be a multi-index. By Lemma \ref{lem:basic-random-Euler-product} and the Cauchy-Schwarz inequality,
\[
\begin{aligned}
\mathbb E\left[
|\zeta(\sigma,X)|^{2\operatorname{Re}s}
|T_\sigma|^{|\alpha|}
\right]
&\leq
(2\sigma-1)^{|\alpha|}
\left(
\mathbb E\left[|\zeta(\sigma,X)|^{4\operatorname{Re}s}\right]
\right)^{1/2}
\left(
\mathbb E\left[|L'(\sigma,X)|^{2|\alpha|}\right]
\right)^{1/2}
<\infty.
\end{aligned}
\]
Hence differentiation under the expectation is justified and yields
\[
\partial^\alpha H_{\sigma,s}(\xi)
=
(-\i)^{|\alpha|}
(2\sigma-1)^{s^2}
\mathbb E\left[
|\zeta(\sigma,X)|^{2s}
T_{\sigma,1}^{\alpha_1}
T_{\sigma,2}^{\alpha_2}
\mathrm e^{-\i\langle\xi,T_\sigma\rangle}
\right].
\]
Dominated convergence shows that these derivatives are continuous, and
\[
\sup_{\xi\in\mathbb R^2}
|\partial^\alpha H_{\sigma,s}(\xi)|
\leq
(2\sigma-1)^{\operatorname{Re}(s^2)}
\mathbb E\left[
|\zeta(\sigma,X)|^{2\operatorname{Re}s}
|T_\sigma|^{|\alpha|}
\right]
<\infty.
\]
Since \(\alpha\) is arbitrary,
$
H_{\sigma,s}\in C_b^\infty(\mathbb R^2).
$ 
It remains to prove integrability. Applying
Proposition~\ref{prop:global-Fourier-decay-rigorous} with, for
instance, \(N=2\), and decreasing \(\sigma_0\) if necessary, gives
\[
|H_{\sigma,s}(\lambda,u)|
\leq
C_{K,2}
(1+\lambda^2+u^2)^{-2}
\]
uniformly for
\[
\frac12<\sigma\leq\sigma_0,
\qquad
s\in K.
\]
Therefore
\[
\begin{aligned}
\|H_{\sigma,s}\|_{L^1(\mathbb R^2)}
&\leq
C_{K,2}
\int_{\mathbb R^2}
(1+\lambda^2+u^2)^{-2}\,\d \lambda\,\d u\\
&=
2\pi C_{K,2}
\int_0^\infty
\frac{r}{(1+r^2)^2}\,dr\\
&=
\pi C_{K,2}<\infty.
\end{aligned}
\]
Hence
$
H_{\sigma,s}\in L^1(\mathbb R^2),
$
and the same estimate also gives
\[
\sup_{\substack{1/2<\sigma\leq\sigma_0\\ s\in K}}
\|H_{\sigma,s}\|_{L^1(\mathbb R^2)}
\leq
\pi C_{K,2}.
\]
This completes the proof.
\end{proof}

By Lemmas~\ref{lem:mixed-moment} and
\ref{lem:weighted-pushforward-nu}, the weighted
mixed moments of $\zeta(\sigma,X)$ and
$L'(\sigma,X)$ are well defined and can be
expressed as integrals with respect to
the weighted measure $\nu_{\sigma,s}$.
Moreover, Lemma~\ref{lem:H-smooth-integrable}
ensures that $H_{\sigma,s}$ satisfies
the conditions required to apply the extended
Riesz distributional pairing established in
Proposition~\ref{lem:Reg-Riesz-Hsigma}. 
Combining these results with the Fourier transform
identity for the homogeneous distribution
$|x|^{2s}$ established in \eqref{1023evening},
we first remove the Gaussian regularization
for each fixed $\sigma>1/2$ sufficiently close
to $1/2$ and derive a Fourier--Riesz representation
of the original complex moments as in Lemma~\ref{lem:moment-Riesz-Hsigma}.
We then use this representation to study
the limit as $\sigma\downarrow1/2$ as in Lemma \ref{lem:Riesz-functional-convergence},
without interchanging the two limits
$t\downarrow0$ and $\sigma\downarrow1/2$ .

\subsection{Removal of Gaussian regularization and convergence}

\subsubsection{Removal of Gaussian regularization and proof of Lemma~\ref{lem:moment-Riesz-Hsigma}}\label{471}

\begin{proof}[Proof of Lemma \ref{lem:moment-Riesz-Hsigma} ]
Let \(K\) be a compact subset of the open right half-plane
$
\{s\in\mathbb C:\operatorname{Re}s>0\},
$
and let \(\sigma_0>1/2\) be as in
Lemma~\ref{lem:H-smooth-integrable}. Fix
\[
1/2<\sigma\leq\sigma_0,
\qquad
s\in K\cap\Omega_M.
\]
Throughout the proof, \(s\) is fixed and \(w\) is regarded as an
auxiliary complex variable. By the definition of \(\nu_{\sigma,s}\) and the moment estimates
\eqref{eq:fixed-sigma-all-moments}, \(\nu\) is a finite complex Borel
measure. Applying Lemma~\ref{lem:weighted-pushforward-nu} with
$
G(x)=\mathrm e^{-\i\langle \xi,x\rangle},
$
we obtain, with the Fourier transform convention,
\[
H_{\sigma,s}(\xi)
=
\widehat{\nu_{\sigma,s}}(\xi)
=
\int_{\mathbb R^2}
\mathrm e^{-\i\langle \xi,x\rangle}\,\d\nu_{\sigma,s}(x).
\]
By Lemma~\ref{lem:H-smooth-integrable},
\[
H_{\sigma,s}\in C_b^\infty(\mathbb R^2)\cap L^1(\mathbb R^2).
\]
Therefore, Fourier inversion for finite measures implies that
\(\nu\) is absolutely continuous with respect to Lebesgue measure:
\begin{equation}\label{absoluteconvergence}
d\nu_{\sigma,s}(x)=h(x)\,\d x,
\end{equation}
where
\begin{equation}
\label{eq:h-density-proof}
h(x)
=
\frac{1}{(2\pi)^2}
\int_{\mathbb R^2}
\mathrm e^{\i \langle x,\xi\rangle}H_{\sigma,s
}(\xi)\,\d\xi.
\end{equation}
In particular,
\begin{equation}\label{propertyofh}
h\in C_0(\mathbb R^2)\cap L^\infty(\mathbb R^2),
\end{equation}
where \(C_0(\mathbb R^2)\) denotes the space of continuous functions
vanishing at infinity.  Hence,
by (\ref{propertyofh}),
\begin{equation}\label{1023formula2}
x\longmapsto |x|^{2\operatorname{Re}(w)}\in L^2(\d |\nu_{\sigma,s}|)
\end{equation}
when
$
-1<\operatorname{Re}(w)<0.
$
In fact, the integrability near the origin follows from the bounded density, while for
$|x|\geq1$ one has
$
|x|^{2\operatorname{Re}(w)}\leq1$.
Applying Lemmas \ref{lem:mixed-moment} and \ref{lem:weighted-pushforward-nu}, for any $w\in \mathbb{C}$ with $\textnormal{Re}(w)>-1$ and $t>0$,
\begin{equation}\label{1023formula}
(2\sigma-1)^{s^2}
\mathbb E\left[
|\zeta(\sigma,X)|^{2s}
|T_\sigma|^{2w}
\mathrm{e}^{-t|T_\sigma|^2}
\right]=
\int_{\mathbb R^2}
|x|^{2w}\e^{-t|x|^2}\,
\d \nu_{\sigma,s}(x)
\end{equation}
We first let $t\downarrow0$ in the strip
\[
-1<\operatorname{Re}(w)<0.
\]
By (\ref{1023formula2}) and the dominated convergence give
\begin{equation}
\label{eq:left-t-zero-removal}
\lim_{t\downarrow0}
\int_{\mathbb R^2}
|x|^{2w}\e^{-t|x|^2}\,
d\nu_{\sigma,s}(x)
=
\int_{\mathbb R^2}|x|^{2w}\,\d \nu_{\sigma,s}(x).
\end{equation}
We next consider the left-handside of \eqref{1023formula} as $t\downarrow 0$.  
By the Fourier representation proved in
Lemma~\ref{lem:smoothed-moment-fourier} gives
\begin{align}
(2\sigma-1)^{s^2}
\mathbb E\left[
|\zeta(\sigma,X)|^{2s}
|T_\sigma|^{2w}
\e^{-t|T_\sigma|^2}
\right]=
\frac{\Gamma(w+1)}{4\pi t^{w+1}}
\int_{\mathbb R^2}
{}_1F_1\left(
w+1;1;-\frac{|\xi|^2}{4t}
\right)
H_{\sigma,s}(\xi)\,\d\xi.
\label{eq:Gaussian-regularized-general-w}
\end{align}
For
$-1<\operatorname{Re}(w)<0
$,
the Kummer asymptotic proved in
(\ref{eq:Kummer-negative-axis}) gives, for every
$\xi\neq0$,
\begin{align}
&t^{-w-1}
{}_1F_1\left(
w+1;1;-\frac{|\xi|^2}{4t}
\right)
\longrightarrow
\frac{4^{w+1}}{\Gamma(-w)}
|\xi|^{-2w-2},
\qquad
t\downarrow0.
\label{eq:Kummer-t-limit-removal}
\end{align}
Moreover, on compact subsets of
\[
-1<\operatorname{Re}(w)<0,
\]
one has, uniformly for $0<t\leq1$,
\begin{equation}
\label{eq:Kummer-majorant-removal}
\left|
t^{-w-1}
{}_1F_1\left(
w+1;1;-\frac{|\xi|^2}{4t}
\right)
\right|
\leq
C_w|\xi|^{-2\operatorname{Re}(w)-2}.
\end{equation}
Since $\operatorname{Re}(w)<0$,
\[
|\xi|^{-2\operatorname{Re}(w)-2}
\]
is locally integrable at $\xi=0$ in two dimensions.  At infinity,
the decay condition (\ref{eq:Hsigma-fixed-decay-fixed-t})
implies that
\[
|\xi|^{-2\operatorname{Re}(w)-2}
|H_{\sigma,s}(\xi)|
\]
is integrable. Hence, by (\ref{eq:mixed-moment-finite}) and (\ref{eq:left-t-zero-removal}), together with the dominated convergence on the right-hand side
of \eqref{eq:Gaussian-regularized-general-w} yields
\begin{align}
&(2\sigma-1)^{s^2}
\mathbb E\left[
|\zeta(\sigma,X)|^{2s}
|T_\sigma|^{2w}\right]=
\int_{\mathbb R^2}
|x|^{2w}\,\d\nu_{\sigma,s}(x)
=
\frac{2^{2w}\Gamma(w+1)}
{\pi\Gamma(-w)}
\int_{\mathbb R^2}
|\xi|^{-2w-2}
H_{\sigma,s}(\xi)\,\d\xi.
\label{eq:Riesz-cartesian-removal}
\end{align}
In the strip $-1<\operatorname{Re}w<0$, the integral is absolutely convergent, and the
ordinary pairing agrees with the extended Gelfand--Shilov pairing:
\[
\int_{\mathbb R^2}
|\xi|^{-2w-2}H_{\sigma,s}(\xi)\,\d\xi
=
\left\langle
\operatorname{Reg}|\xi|^{-2w-2},
H_{\sigma,s}
\right\rangle_{\mathrm{ext}}.
\]
Applying Proposition~\ref{lem:Reg-Riesz-Hsigma}, we obtain
\[
\left\langle
\operatorname{Reg}|\xi|^{-2w-2},
H_{\sigma,s}
\right\rangle_{\mathrm{ext}}
=
2\pi\mathcal R_{M,w}[H_{\sigma,s}],
\]
and hence
\begin{equation}\label{eq:moment-Riesz-initial-strip}
\int_{\mathbb R^2}
|x|^{2w}\,\d\nu_{\sigma,s}(x)
=
2^{2w+1}
\frac{\Gamma(w+1)}{\Gamma(-w)}
\mathcal R_{M,w}[H_{\sigma,s}],
\qquad
-1<\operatorname{Re}w<0.
\end{equation}
By Proposition~\ref{lem:Reg-Riesz-Hsigma},
\begin{equation*}
w\longmapsto
\mathcal R_{M,w}[H_{\sigma,s}]
\end{equation*}
is holomorphic on \(\Omega_M\). Therefore
\begin{equation}\label{analyticextensionfunction}
w\longmapsto
2^{2w+1}
\frac{\Gamma(w+1)}{\Gamma(-w)}
\mathcal R_{M,w}[H_{\sigma,s}]
\end{equation}
is holomorphic on \(\Omega_M\). Hence the  function (\ref{analyticextensionfunction}) gives an analytic continuation of
$w\longmapsto (2\sigma-1)^{s^2}
\mathbb E\left[
|\zeta(\sigma,X)|^{2s}
|T_\sigma|^{2w}\right]$ from the initial strip $\{w\in \mathbb{C}: -1<\textnormal{Re}(w)<0\}$ to $\Omega_{M}$. Note that,
for any $w\in \mathbb{C}$ with $\textnormal{Re}(w)>-1$,
\begin{equation}\label{1023formula}
(2\sigma-1)^{s^2}
\mathbb E\left[
|\zeta(\sigma,X)|^{2s}
|T_\sigma|^{2w}
\right]=
\int_{\mathbb R^2}
|x|^{2w}\,
\d\nu_{\sigma,s}(x)
\end{equation}
by Lemmas \ref{lem:mixed-moment} and \ref{lem:weighted-pushforward-nu}. To prove the claim of this lemma, it remains only to show that
\begin{equation}
\label{anothergoal}
F_{\sigma,s}(w)
:=
\int_{\mathbb R^2}|x|^{2w}\,\d\nu_{\sigma,s}(x)
\end{equation}
is holomorphic on
\(\{\operatorname{Re}w>-1\}\). Indeed, once this is established, the connectedness of \(\Omega_M\),
together with the identity theorem applied to
\eqref{eq:moment-Riesz-initial-strip}, implies, in view of
\eqref{1023formula}, that
\begin{equation}
\label{1023formula4}
(2\sigma-1)^{s^2}
\mathbb E\left[
|\zeta(\sigma,X)|^{2s}
|T_\sigma|^{2w}
\right]
=
2^{2w+1}
\frac{\Gamma(w+1)}{\Gamma(-w)}
\mathcal R_{M,w}[H_{\sigma,s}],
\qquad
w\in\Omega_M.
\end{equation}
Equation~\eqref{eq:moment-Riesz-Hsigma} then follows by setting
\(w=s\) in \eqref{1023formula4}. 

It therefore remains to prove that \(F_{\sigma,s}(w)\) is holomorphic
on $\{w\in \mathbb{C}: \operatorname{Re}w>-1\}$. Let $L\subset\{w\in\mathbb C:\operatorname{Re}w>-1\}$ be compact, and set
\[
a:=\inf_{w\in L}\operatorname{Re}w>-1,
\qquad
b:=\sup_{w\in L}\operatorname{Re}w.
\]
For \(x\neq0\), the function
\[
w\longmapsto |x|^{2w}
=
\exp\bigl(2w\log|x|\bigr)
\]
is entire, and for every \(k\in\mathbb Z_{\geq0}\),
\[
\frac{\partial^k}{\partial w^k}|x|^{2w}
=
(2\log|x|)^k|x|^{2w}.
\]
We now establish an integrable majorant, uniformly for \(w\in L\).
For \(|x|\leq1\),
\[
|(2\log|x|)^k|x|^{2w}|
\leq
2^k|\log|x||^k|x|^{2a}.
\]
By (\ref{absoluteconvergence}),
\(d\nu_{\sigma,s}=h\,dx\) with
\(h\in L^\infty(\mathbb R^2)\), we have
\[
\begin{aligned}
\int_{|x|\leq1}
|\log|x||^k|x|^{2a}\,
\d|\nu_{\sigma,s}|(x)
&=
\int_{|x|\leq1}
|\log|x||^k|x|^{2a}
|h(x)|\,\d x\\
&\leq
2\pi\|h\|_\infty
\int_0^1
|\log r|^k r^{2a+1}\,\d r
<\infty,
\end{aligned}
\]
because \(a>-1\). For \(|x|\geq1\), choose \(q>\max\{2b,0\}\). Then
\[
(\log r)^k r^{2b}\ll_{L,k,q} r^q,
\qquad r\geq1,
\]
and therefore
\[
|(2\log|x|)^k|x|^{2w}|
\ll_{L,k,q}|x|^q,
\qquad |x|\geq1.
\]
Moreover,
\[
\int_{\mathbb R^2}|x|^q\,d|\nu_{\sigma,s}|(x)<\infty
\]
by Lemma \ref{lem:weighted-pushforward-nu} and (\ref{eq:mixed-moment-finite}). Thus, for every \(k\geq0\), the family
\[
\left\{
(2\log|x|)^k|x|^{2w}:w\in L
\right\}
\]
is dominated by a \(d|\nu_{\sigma,s}|\)-integrable function.
Hence differentiation under the integral sign is justified, and
\[
F_{\sigma,s}^{(k)}(w)
=
\int_{\mathbb R^2}
(2\log|x|)^k|x|^{2w}\,
\d \nu_{\sigma,s}(x),
\qquad w\in \{w\in \mathbb{C}: \operatorname{Re}w>-1\}.
\]
In particular, \(F_{\sigma,s}\) is holomorphic on $\{w\in \mathbb{C}: \operatorname{Re}w>-1\}$. 
\end{proof}

\subsubsection{Convergence of the regularized Fourier--Riesz functional}

We now turn to the limit $\sigma\downarrow1/2$ in the regularized
Fourier--Riesz functional
\[
\mathcal R_{M,s}[H_{\sigma,s}].
\]
We treat separately the contributions from a
neighbourhood of the origin in Fourier space,
the Fourier tail, and the finite correction terms
arising from Taylor subtraction. Near the origin, the factor $r^{-2s-1}$ is singular, so local
uniform convergence of $H_{\sigma,s}$ alone is not sufficient.
Instead, we use the local Fourier derivative bounds to subtract the
circular Taylor polynomial uniformly and obtain a remainder of order
$r^{2M+2}$.  This produces an integrable majorant provided
$\sup_{s\in K}\operatorname{Re}(s)<M+1$.  The precise estimate is
given in the following lemma.

\begin{lem}[Uniform control near the Fourier origin]
\label{lem:local-Riesz-origin-control}
Let \(K\) be a compact subset of the open right half-plane
\[
\{s\in\mathbb C:\operatorname{Re}s>0\}.
\]
and let \(M\in\mathbb Z_{\geq0}\) satisfy
\begin{equation}\label{restriction}
\sup_{s\in K}\operatorname{Re}s<M+1.
\end{equation}
Let  \(H_{\sigma,s}\) be defined in \eqref{eq:Hsigma-def-final}. Then there exist \(\sigma_0>1/2\) and \(C_{K,M}>0\) such that
\[
\Bigg|
\mathcal A H_{\sigma,s}(r)
-
\sum_{j=0}^{M}
\frac{r^{2j}}{4^j(j!)^2}
\Delta^jH_{\sigma,s}(0)
\Bigg|
\leq
C_{K,M}r^{2M+2}
\]
uniformly for
\[
\frac12<\sigma\leq\sigma_0,
\qquad
s\in K,
\qquad
0\leq r\leq1.
\]
Moreover, as \(\sigma\downarrow1/2\),
\[
\begin{aligned}
&\int_0^1
\Bigg[
\mathcal A H_{\sigma,s}(r)
-
\sum_{j=0}^{M}
\frac{r^{2j}}{4^j(j!)^2}
\Delta^jH_{\sigma,s}(0)
\Bigg]
r^{-2s-1}\,\d r\\
&\qquad\longrightarrow
\int_0^1
\Bigg[
\mathcal A H_s(r)
-
\sum_{j=0}^{M}
\frac{r^{2j}}{4^j(j!)^2}
\Delta^jH_s(0)
\Bigg]
r^{-2s-1}\,\d r
\end{aligned}
\]
uniformly for \(s\in K\).
Here, 
\begin{equation}\label{1023eq}
H_s(\lambda,u)
=
a_s
\mathrm e^{-\i s\lambda-(\lambda^2+u^2)/4}
\end{equation}
and 
$\mathcal{A}H$ is defined by (\ref{def:circular-mean2}).
\end{lem}
\begin{proof}
Taking \(R=1\) in
Lemma~\ref{lem:local-Fourier-derivatives-rigorous} and choosing a
common \(\sigma_0>1/2\) for the finitely many pairs
\[
(\alpha,\beta)
=
(j,2M+2-j),
\qquad
0\leq j\leq2M+2,
\]
the chain rule gives
\[
\sup_{\substack{
1/2<\sigma\leq\sigma_0\\
s\in K,\ 0\leq r\leq1\\
\theta\in[0,2\pi]}}
\left|
\frac{d^{2M+2}}{dr^{2M+2}}
H_{\sigma,s}(r\cos\theta,r\sin\theta)
\right|
\leq C_{K,M}<\infty.
\]
By Taylor's theorem at \(r=0\), with the remainder in integral form,
we obtain
\begin{align}
H_{\sigma,s}(r\cos\theta, r\sin\theta)
={}&
\sum_{k=0}^{2M+1}
\frac{r^k}{k!}
\left(
\cos\theta\,\partial_\lambda
+
\sin\theta\,\partial_u
\right)^k
H_{\sigma,s}(0)+
O_{K,M}(r^{2M+2}),
\label{eq:directional-Taylor-origin}
\end{align}
where the error term is uniform for
\[
\frac12<\sigma\leq\sigma_0,
\qquad
s\in K,
\qquad
0\leq r\leq1,
\qquad
\theta\in[0,2\pi].
\]
Averaging \eqref{eq:directional-Taylor-origin} in \(\theta\), all odd
terms vanish. For the even terms,
\begin{equation}
\label{eq:angular-differential-identity}
\frac1{2\pi}
\int_0^{2\pi}
\left(
\cos\theta\,\partial_\lambda
+
\sin\theta\,\partial_u
\right)^{2j}
\,d\theta
=
\frac{(2j)!}{4^j(j!)^2}\Delta^j,
\end{equation}
where
\[
\Delta=\partial_\lambda^2+\partial_u^2.
\]
Indeed, note that the differential operators \(\partial_\lambda\) and
\(\partial_u\) commute on sufficiently smooth functions, the operator on the left-hand side is a rotationally invariant
homogeneous polynomial of degree \(2j\) in
\(\partial_\lambda,\partial_u\), and is therefore a constant multiple
of \(\Delta^j\). The constant is determined by the coefficient of
\(\partial_\lambda^{2j}\):
\[
\frac1{2\pi}
\int_0^{2\pi}\cos^{2j}\theta\,\d \theta
=
\frac{(2j)!}{4^j(j!)^2}.
\]
It follows that
\begin{equation}
\label{eq:circular-Taylor-origin}
\mathcal A H_{\sigma,s}(r)
=
\sum_{j=0}^{M}
\frac{r^{2j}}{4^j(j!)^2}
\Delta^jH_{\sigma,s}(0)
+
O_{K,M}(r^{2M+2}),
\end{equation}
uniformly for
\[
\frac12<\sigma\leq\sigma_0,
\qquad
s\in K,
\qquad
0\leq r\leq1.
\]
The locally uniform convergence in
Lemma~\ref{lem:local-Fourier-derivatives-rigorous} also implies the
corresponding uniform derivative bounds for the limiting function
\(H_s\). Hence the same argument gives
\begin{equation}
\label{eq:circular-Taylor-origin-limit}
\mathcal A H_s(r)
=
\sum_{j=0}^{M}
\frac{r^{2j}}{4^j(j!)^2}
\Delta^jH_s(0)
+
O_{K,M}(r^{2M+2}),
\end{equation}
uniformly for
\[
s\in K,
\qquad
0\leq r\leq1.
\]
Define
\begin{equation}\label{defofRsigmas}
R_{\sigma,s}(r)
:=
\mathcal A H_{\sigma,s}(r)
-
\sum_{j=0}^{M}
\frac{r^{2j}}{4^j(j!)^2}
\Delta^jH_{\sigma,s}(0)
\end{equation}
and
\begin{equation}\label{defofRs}
R_s(r)
:=
\mathcal A H_s(r)
-
\sum_{j=0}^{M}
\frac{r^{2j}}{4^j(j!)^2}
\Delta^jH_s(0).
\end{equation}
Equations \eqref{eq:circular-Taylor-origin} and
\eqref{eq:circular-Taylor-origin-limit} imply
\begin{equation}
\label{eq:Riesz-origin-remainder-bound}
|R_{\sigma,s}(r)|+|R_s(r)|
\leq
C_{K,M}r^{2M+2},
\qquad
0\leq r\leq1.
\end{equation}
Moreover, Lemma~\ref{lem:local-Fourier-derivatives-rigorous} gives
\[
\sup_{\substack{s\in K\\0\leq r\leq1}}
\left|
\mathcal A H_{\sigma,s}(r)-\mathcal A H_s(r)
\right|
\longrightarrow0
\]
and, for every \(0\leq j\leq M\),
\[
\sup_{s\in K}
\left|
\Delta^jH_{\sigma,s}(0)-\Delta^jH_s(0)
\right|
\longrightarrow0.
\]
Therefore
\begin{equation}
\label{eq:Riesz-remainder-convergence}
\sup_{\substack{s\in K\\0\leq r\leq1}}
|R_{\sigma,s}(r)-R_s(r)|
\longrightarrow0
\end{equation}
as \(\sigma\downarrow1/2\).
It remains to prove the uniform convergence of the integrals. By \eqref{restriction} and \eqref{eq:Riesz-origin-remainder-bound},
\[
\begin{aligned}
\left|
\bigl(R_{\sigma,s}(r)-R_s(r)\bigr)r^{-2s-1}
\right|
&\leq
C_{K,M}
r^{2M+1-2\sup_{s\in K}\operatorname{Re}s}\in L^{1}(0,1].
\end{aligned}
\]
Hence, for $\epsilon>0$, choose \(\delta\in(0,1)\) such that
\[
\int_0^\delta
\left|
R_{\sigma,s}(r)-R_s(r)
\right|
r^{-2\operatorname{Re}s-1}\,dr
<
\frac{\epsilon}{2}.
\]
uniformly for \(s\in K\). On the interval \([\delta,1]\), one has
\[
\begin{aligned}
\sup_{s\in K}
\left|
\int_\delta^1
\bigl(R_{\sigma,s}(r)-R_s(r)\bigr)
r^{-2s-1}\,dr
\right|\leq
\delta^{-2\sup_{s\in K}\operatorname{Re}s-1}
\sup_{\substack{s\in K\\\delta\leq r\leq1}}
|R_{\sigma,s}(r)-R_s(r)|.
\end{aligned}
\]
By \eqref{eq:Riesz-remainder-convergence}, the right-hand side tends
to zero as \(\sigma\downarrow1/2\). Thus it is smaller than
\(\epsilon/2\) for \(\sigma\) sufficiently close to \(1/2\). Combining the estimates on \((0,\delta)\) and \([\delta,1]\), we
obtain
\[
\begin{aligned}
\sup_{s\in K}
\left|
\int_0^1
R_{\sigma,s}(r)r^{-2s-1}\,dr
-
\int_0^1
R_s(r)r^{-2s-1}\,dr
\right|
\longrightarrow0,
\end{aligned}
\]
which proves the desired uniform convergence.
\end{proof}

The
contribution near the Fourier origin is controlled by the local
Fourier derivative estimates and the resulting uniform Taylor
remainder, while the contribution from large Fourier variables is
controlled by the global decay of $H_{\sigma,s}$.  The finite
correction terms are handled by the local convergence of the
derivatives at the origin.  The following lemma makes these three
steps precise.

\begin{lem}[Convergence of the regularized Fourier--Riesz functional]
\label{lem:Riesz-functional-convergence}
Let
\[
s\in
\left\{
s\in\mathbb C:\operatorname{Re}s>0
\right\}
\setminus\mathbb Z_{>0}.
\]
Let \(H_{\sigma,s}\) and \(H_s\) be defined by
\eqref{eq:Hsigma-def-final} and \eqref{1023eq}, respectively. Then,
for every \(M\in\mathbb Z_{\geq0}\) satisfying
\[
\operatorname{Re}s<M+1,
\]
one has
\begin{equation}
\label{eq:Riesz-functional-convergence}
\mathcal R_{M,s}[H_{\sigma,s}]
\longrightarrow
\mathcal R_{M,s}[H_s],
\qquad
\text{as }\sigma\downarrow\frac12,
\end{equation}
where \(\mathcal R_{M,s}[\cdot]\) is defined by
\eqref{eq:Riesz-functional-abstract}. Moreover, if \(K\) is a compact subset of $\left\{
s\in\mathbb C:\operatorname{Re}s>0
\right\}
\setminus\mathbb Z_{>0}$ and
\(M\in\mathbb Z_{\geq0}\) satisfies
\[
\sup_{s\in K}\operatorname{Re}s<M+1,
\]
then
\[
\sup_{s\in K}
\left|
\mathcal R_{M,s}[H_{\sigma,s}]
-
\mathcal R_{M,s}[H_s]
\right|
\longrightarrow0,
\qquad
\text{as }\sigma\downarrow\frac12.
\]
\end{lem}
\begin{proof}
Let
\[
D
:=
\left\{
s\in\mathbb C:\operatorname{Re}s>0
\right\}
\setminus\mathbb Z_{>0}.
\]
We first prove the locally uniform statement. Let \(K\) be a compact
subset of \(D\), and let \(M\in\mathbb Z_{\geq0}\) satisfy
\[
\sup_{s\in K}\operatorname{Re}s
<
M+1.
\]
We shall prove that
\[
\sup_{s\in K}
\left|
\mathcal R_{M,s}[H_{\sigma,s}]
-
\mathcal R_{M,s}[H_s]
\right|
\longrightarrow0
\qquad
\text{as }\sigma\downarrow\frac12.
\]
By the definition of \(\mathcal R_{M,s}\), we have
\begin{align}
&\mathcal R_{M,s}[H_{\sigma,s}]
-
\mathcal R_{M,s}[H_s]
\nonumber\\
&=
\int_0^1
\left[
R_{\sigma,s}(r)-R_s(r)
\right]
r^{-2s-1}\,\d r
+
\int_1^\infty
\mathcal A
\left(
H_{\sigma,s}-H_s
\right)(r)
r^{-2s-1}\,\d r
\nonumber\\
&\quad+
\sum_{j=0}^{M}
\frac{
\Delta^jH_{\sigma,s}(0)
-
\Delta^jH_s(0)
}{
4^j(j!)^2(2j-2s)
},
\label{eq:Riesz-functional-difference}
\end{align}
where $R_{\sigma,s}(r)$ and $R_s(r)$ are defined by (\ref{defofRsigmas}) and (\ref{defofRs}), respectively. We treat the three terms on the right-hand side of
\eqref{eq:Riesz-functional-difference} separately. First, Lemma~\ref{lem:local-Riesz-origin-control} gives
\begin{equation}
\label{eq:Riesz-local-part-convergence}
\sup_{s\in K}
\left|
\int_0^1
\left[
R_{\sigma,s}(r)-R_s(r)
\right]
r^{-2s-1}\,\d r
\right|
\longrightarrow0
\end{equation}
as \(\sigma\downarrow1/2\). We next consider the finite sum. Since \(K\) is compact and
\[
K\cap\{0,1,\ldots,M\}=\varnothing,
\]
we have
\[
d_{K,M}
:=
\min_{\substack{s\in K\\0\leq j\leq M}}
|j-s|
>0.
\]
By Lemma~\ref{lem:local-Fourier-derivatives-rigorous}, for every
\(0\leq j\leq M\),
\[
\Delta^jH_{\sigma,s}(0)
\longrightarrow
\Delta^jH_s(0)
\]
uniformly for \(s\in K\). Therefore
\begin{align}
&\sup_{s\in K}
\left|
\sum_{j=0}^{M}
\frac{
\Delta^jH_{\sigma,s}(0)
-
\Delta^jH_s(0)
}{
2\,4^j(j!)^2(j-s)
}
\right|
\nonumber\\
&\qquad\leq
\sum_{j=0}^{M}
\frac{1}{
2\,4^j(j!)^2d_{K,M}
}
\sup_{s\in K}
\left|
\Delta^jH_{\sigma,s}(0)
-
\Delta^jH_s(0)
\right|
\longrightarrow0.
\label{eq:Riesz-finite-sum-convergence}
\end{align}
It remains to treat the integral over \(r\geq1\). By
Proposition~\ref{prop:global-Fourier-decay-rigorous}, after decreasing
\(\sigma_0>1/2\) if necessary, for every \(N>0\) there exists
\(C_{K,N}>0\) such that
\[
|H_{\sigma,s}(\lambda,u)|
\leq
C_{K,N}
(1+\lambda^2+u^2)^{-N}
\]
for
\[
\frac12<\sigma\leq\sigma_0,
\qquad
s\in K,
\qquad
(\lambda,u)\in\mathbb R^2.
\]
The explicit formula
\[
H_s(\lambda,u)
=
a_s
\mathrm e^{-\i s\lambda-(\lambda^2+u^2)/4}
\]
implies the same type of rapid-decay bound for \(H_s\), uniformly for
\(s\in K\). Consequently, combining this with the definition of \(\mathcal A H\) in
\eqref{def:circular-mean2},
\begin{equation}
\label{eq:circular-tail-uniform-bound}
\sup_{\substack{
1/2<\sigma\leq\sigma_0\\
s\in K}}
\left|
\mathcal A
\left(
H_{\sigma,s}-H_s
\right)(r)
\right|
\leq
C_{K,N}(1+r^2)^{-N},
\qquad r\geq1.
\end{equation}
Fix \(R>1\). It follows that
\begin{align}
&\sup_{s\in K}
\left|
\int_1^\infty
\mathcal A
\left(
H_{\sigma,s}-H_s
\right)(r)
r^{-2s-1}\,\d r
\right|
\nonumber\\
&\quad\leq
\int_1^R
\sup_{s\in K}
\left|
\mathcal A
\left(
H_{\sigma,s}-H_s
\right)(r)
\right|
r^{-2\inf_{s\in K}\operatorname{Re}s-1}\,\d r
\nonumber\\
&\qquad+
C_{K,N}
\int_R^\infty
(1+r^2)^{-N}r^{-2\inf_{s\in K}\operatorname{Re}s-1}\,\d r.
\label{eq:Riesz-tail-split}
\end{align}
For fixed \(R\), Lemma~\ref{lem:local-Fourier-derivatives-rigorous},
with \(\alpha=\beta=0\), gives
\[
H_{\sigma,s}(\lambda,u)
\longrightarrow
H_s(\lambda,u)
\]
uniformly for
\[
s\in K,
\qquad
\lambda^2+u^2\leq R^2.
\]
Hence
\[
\sup_{\substack{s\in K\\1\leq r\leq R}}
\left|
\mathcal A
\left(
H_{\sigma,s}-H_s
\right)(r)
\right|
\longrightarrow0.
\]
Thus the first integral on the right-hand side of
\eqref{eq:Riesz-tail-split} tends to zero as
\(\sigma\downarrow1/2\). On the other hand, the second integral in
\eqref{eq:Riesz-tail-split} tends to zero as \(R\to\infty\),
uniformly in \(\sigma\) and \(s\in K\). Therefore,
\begin{equation}
\label{eq:Riesz-tail-convergence}
\sup_{s\in K}
\left|
\int_1^\infty
\mathcal A
\left(
H_{\sigma,s}-H_s
\right)(r)
r^{-2s-1}\,\d r
\right|
\longrightarrow0
\end{equation}
as \(\sigma\downarrow1/2\). Combining
\eqref{eq:Riesz-local-part-convergence},
\eqref{eq:Riesz-finite-sum-convergence}, and
\eqref{eq:Riesz-tail-convergence} in
\eqref{eq:Riesz-functional-difference}, we obtain
\[
\sup_{s\in K}
\left|
\mathcal R_{M,s}[H_{\sigma,s}]
-
\mathcal R_{M,s}[H_s]
\right|
\longrightarrow0
\qquad
\text{as }\sigma\downarrow\frac12.
\]
This proves the stated local uniform convergence. The pointwise convergence follows immediately from the compact-uniform
convergence established above. This completes the proof.
\end{proof}
We have now shown that the regularized Fourier--Riesz functional
associated with $H_{\sigma,s}$ converges, as $\sigma\downarrow1/2$,
to the corresponding functional associated with the limiting
Fourier transform
\[
H_s(\lambda,u)
=
a_s \e^{-\i s\lambda-(\lambda^2+u^2)/4}.
\]
It therefore remains to evaluate
\[
\mathcal R_{M,s}[H_s]
\]
explicitly.  Since $H_s$ is the Fourier transform of a shifted
complex Gaussian density, the same Fourier--Riesz argument as above
reduces this evaluation to an explicit Gaussian integral.  The
result is given in the following lemma.

\subsection{Evaluation of the limiting Fourier--Riesz functional}
\begin{lem}
\label{lem:limiting-Riesz-evaluation}
Let
\[
s\in\mathbb C\setminus\mathbb Z_{>0},
\qquad
\operatorname{Re}s>0,
\]
and let \(a_s\) be the arithmetic factor defined in
\eqref{definitionofthearithmeticfactor}. Then, for any
\(M\in\mathbb Z_{\geq0}\) satisfying
\[
\operatorname{Re}s<M+1,
\]
one has
\begin{equation}
\label{eq:limiting-Riesz-evaluation}
2^{2s+1}
\frac{\Gamma(s+1)}{\Gamma(-s)}
\mathcal R_{M,s}[H_s]
=
a_s \e^{-s^2}\Gamma(s+1)
{}_1F_1(s+1;1;s^2),
\end{equation}
where
\begin{equation}\label{1023eq}
H_s(\lambda,u)
=
a_s
\mathrm e^{-\i s\lambda-(\lambda^2+u^2)/4}.
\end{equation}
\end{lem}

\begin{proof}
Fix
\[
s\in\mathbb C\setminus\mathbb Z_{>0},
\qquad
\operatorname{Re}s>0,
\]
and let \(M\in\mathbb Z_{\geq0}\) satisfy
\[
\operatorname{Re}s<M+1.
\]
Then \(s\in\Omega_M\). For \(t>0\), set
\[
f_{s,t}(x):=|x|^{2s}\mathrm e^{-t|x|^2}.
\]
By the definition of \(\mathcal M_t(s)\) in
\eqref{eq:fixed-t-smoothed-limit}, together with
\eqref{eq:gs-Fourier-for-Mt}, we have
\begin{equation}
\label{eq:Mt-pairing-Hs}
\mathcal M_t(s)
=
\frac{1}{(2\pi)^2}
\left\langle
\widehat{f_{s,t}},H_s
\right\rangle,
\end{equation}
where $H_{s}(\lambda,u)$ is defined in (\ref{1023eq}).
In particular,
$
H_s\in\mathcal S(\mathbb R^2).
$ By \eqref{eq:fw-Sprime-limit},
\[
f_{s,t}
\longrightarrow
|x|^{2s}
\qquad\text{in }\mathcal S'(\mathbb R^2)
\]
as \(t\downarrow0\). Since the Fourier transform is a continuous
automorphism of \(\mathcal S'(\mathbb R^2)\), it follows that
\[
\widehat{f_{s,t}}
\longrightarrow
\widehat{|x|^{2s}}
\qquad\text{in }\mathcal S'(\mathbb R^2).
\]
Since \(H_s\in\mathcal S(\mathbb R^2)\), passing to the limit in
\eqref{eq:Mt-pairing-Hs} gives
\begin{equation}
\label{eq:Mt-limit-distribution}
\lim_{t\downarrow0}\mathcal M_t(s)
=
\frac{1}{(2\pi)^2}
\left\langle
\widehat{|x|^{2s}},H_s
\right\rangle.
\end{equation}
By \eqref{1023evening}, with \(w=s\),
\[
\widehat{|x|^{2s}}(\xi)
=
2^{2s+2}\pi
\frac{\Gamma(s+1)}{\Gamma(-s)}
\,\operatorname{Reg}|\xi|^{-2s-2}.
\]
Hence \eqref{eq:Mt-limit-distribution} becomes
\[
\lim_{t\downarrow0}\mathcal M_t(s)
=
\frac{2^{2s+2}\pi}{(2\pi)^2}
\frac{\Gamma(s+1)}{\Gamma(-s)}
\left\langle
\operatorname{Reg}|\xi|^{-2s-2},
H_s
\right\rangle.
\]
Since \(H_s\in\mathcal S(\mathbb R^2)\), the usual
\(\mathcal S'-\mathcal S\) pairing agrees with the extended pairing
of Proposition~\ref{lem:Reg-Riesz-Hsigma}. Since \(s\in\Omega_M\),
that proposition therefore gives
\[
\left\langle
\operatorname{Reg}|\xi|^{-2s-2},
H_s
\right\rangle
=
2\pi\mathcal R_{M,s}[H_s].
\]
Consequently,
\begin{equation}
\label{eq:Mt-limit-Riesz}
\lim_{t\downarrow0}\mathcal M_t(s)
=
2^{2s+1}
\frac{\Gamma(s+1)}{\Gamma(-s)}
\mathcal R_{M,s}[H_s].
\end{equation}
On the other hand, Lemma~\ref{lem:evaluation-fixed-t-limit} gives
\[
\mathcal M_t(s)
=
a_s\mathrm e^{-s^2}\Gamma(s+1)
(1+t)^{-s-1}
{}_1F_1\left(
s+1;1;\frac{s^2}{1+t}
\right).
\]
Letting \(t\downarrow0\), we obtain
\begin{equation}
\label{eq:Mt-limit-explicit}
\lim_{t\downarrow0}\mathcal M_t(s)
=
a_s\mathrm e^{-s^2}\Gamma(s+1)
{}_1F_1(s+1;1;s^2).
\end{equation}
Comparing \eqref{eq:Mt-limit-Riesz} and
\eqref{eq:Mt-limit-explicit} yields
\eqref{eq:limiting-Riesz-evaluation}.
\end{proof}

\subsection{Proof of Theorem~\ref{thm:complex-random-derivative-moments}}\label{thm41}

We now combine the preceding results to prove
Theorem~\ref{thm:complex-random-derivative-moments}.

\begin{proof}[Proof of Theorem~\ref{thm:complex-random-derivative-moments}]
We first consider
\[
s\in\mathbb C\setminus\mathbb Z_{>0},
\qquad
\operatorname{Re}s>0.
\]
Choose \(M\in\mathbb Z_{\geq0}\) such that
\[
\operatorname{Re}s<M+1.
\]
Recall that
\[
T_\sigma=(2\sigma-1)L'(\sigma,X)
\]
and
\[
\zeta'(\sigma,X)
=
\zeta(\sigma,X)L'(\sigma,X).
\]
Consequently,
\[
\begin{aligned}
(2\sigma-1)^{s^2+2s}
\mathbb E\left[
|\zeta'(\sigma,X)|^{2s}
\right]=
(2\sigma-1)^{s^2}
\mathbb E\left[
|\zeta(\sigma,X)|^{2s}
|T_\sigma|^{2s}
\right].
\end{aligned}
\]
By Lemma~\ref{lem:moment-Riesz-Hsigma},
for all $\sigma>1/2$ sufficiently close to $1/2$,
\begin{equation}
\label{eq:final-proof-Riesz-representation}
\begin{aligned}
(2\sigma-1)^{s^2+2s}
\mathbb E\left[
|\zeta'(\sigma,X)|^{2s}
\right]=
2^{2s+1}
\frac{\Gamma(s+1)}{\Gamma(-s)}
\mathcal R_{M,s}[H_{\sigma,s}].
\end{aligned}
\end{equation}
By Lemma~\ref{lem:Riesz-functional-convergence},
\[
\mathcal R_{M,s}[H_{\sigma,s}]
\longrightarrow
\mathcal R_{M,s}[H_s]
\qquad
\text{as }\sigma\downarrow\frac12.
\]
Combining this convergence with
Lemma~\ref{lem:limiting-Riesz-evaluation},
we obtain
\begin{equation}
\label{finalconclusion2}
\begin{aligned}
\lim_{\sigma\downarrow1/2}
(2\sigma-1)^{s^2+2s}
\mathbb E\left[
|\zeta'(\sigma,X)|^{2s}
\right]=
a_s\mathrm e^{-s^2}\Gamma(s+1)
{}_1F_1(s+1;1;s^2).
\end{aligned}
\end{equation}
Moreover, the convergence is locally uniform
for
\[
s\in D:=
\{s\in\mathbb C:\operatorname{Re}s>0\}
\setminus\mathbb Z_{>0},
\]
by the locally uniform convergence established
in Lemma~\ref{lem:Riesz-functional-convergence}. It remains to extend the result to
$s\in\mathbb Z_{>0}$.
For $\sigma>1/2$, define
\[
F_\sigma(s)
:=
(2\sigma-1)^{s^2+2s}
\mathbb E\left[
|\zeta'(\sigma,X)|^{2s}
\right],
\]
and let
\[
F(s)
:=
a_s\mathrm e^{-s^2}\Gamma(s+1)
{}_1F_1(s+1;1;s^2).
\]
By \eqref{eq:fixed-sigma-all-moments}, for every
compact set $K\subset\{s:\operatorname{Re}s>0\}$,
there exists $A>0$ such that
\[
\sup_{s\in K}
\left||\zeta'(\sigma,X)|^{2s}\right|
\leq 1+|\zeta'(\sigma,X)|^{2A},
\]
whose expectation is finite.
Since $|\zeta'(\sigma,X)|^{2s}$ is holomorphic
in $s$ for $\operatorname{Re}s>0$,
dominated convergence and Morera's theorem
imply that
$\mathbb E[|\zeta'(\sigma,X)|^{2s}]$
is holomorphic in this half-plane.
Consequently, $F_\sigma$ is holomorphic
there.
Consequently, $F_\sigma(s)$ is holomorphic
in this half-plane. The function $F$ is also holomorphic in
the half-plane $\operatorname{Re}s>0$.
Indeed, the Euler product defining $a_s$
converges locally uniformly in $s$, since
each Euler factor equals
$1+O_K(p^{-2})$ uniformly on every compact
set $K\subset\mathbb C$.
Moreover, $\Gamma(s+1)$ is holomorphic
for $\operatorname{Re}s>0$, while
$\mathrm e^{-s^2}$ and
${}_1F_1(s+1;1;s^2)$ are entire functions
of $s$. Fix $m\in\mathbb Z_{>0}$ and choose
$0<r<1/2$.
The closed disc
\[
\overline{B(m,r)}
=
\{s\in\mathbb C:|s-m|\leq r\}
\]
lies in the half-plane $\operatorname{Re}s>0$,
and its boundary contains no positive integers.
Since $F_\sigma-F$ is holomorphic on
a neighbourhood of $\overline{B(m,r)}$,
the maximum modulus principle gives
\[
\sup_{|s-m|\leq r}
|F_\sigma(s)-F(s)|
\leq
\sup_{|s-m|=r}
|F_\sigma(s)-F(s)|.
\]
The right-hand side tends to zero
as $\sigma\downarrow1/2$, by the locally
uniform convergence already established on $D$.
Hence
\[
\sup_{|s-m|\leq r}
|F_\sigma(s)-F(s)|
\longrightarrow0.
\]
In particular,
\[
\lim_{\sigma\downarrow1/2}
F_\sigma(m)=F(m).
\]
Since $m\in\mathbb Z_{>0}$ was arbitrary,
the asymptotic formula holds for every
$s\in\mathbb C$ with $\operatorname{Re}s>0$.
The locally uniform convergence on $D$,
together with the uniform convergence
on the discs around positive integers,
also establishes locally uniform convergence
throughout the half-plane $\operatorname{Re}s>0$.
This completes the proof.
\end{proof}

\bibliographystyle{plain}
\bibliography{main}

\end{document}